\documentclass[12pt]{amsart}

\usepackage{times}

\usepackage{amsmath}
\usepackage{amssymb}
\usepackage{amsthm}
\usepackage{amscd}
\usepackage{xypic}

\headsep=1cm \numberwithin{equation}{section}
\newtheorem{theorem}{Theorem}[section]
\newtheorem{lemma}[theorem]{Lemma}
\newtheorem{proposition}[theorem]{Proposition}
\newtheorem{corollary}[theorem]{Corollary}

\theoremstyle{definition}
\newtheorem{definition}[theorem]{Definition}
\newtheorem{remark}[theorem]{Remark}
\newtheorem{convention}[theorem]{Convention}
\newtheorem{reminder}[theorem]{Reminder}
\newtheorem{notation}[theorem]{Notation}
\newtheorem{example}[theorem]{Example}

\newcommand\Proj{\operatorname{Proj}}
\newcommand\Spec{\operatorname{Spec}}
\newcommand\Hom{\operatorname{Hom}}
\newcommand\Ext{\operatorname{Ext}}
\newcommand\Tor{\operatorname{Tor}}
\newcommand\height{\operatorname{height}}
\newcommand\depth{\operatorname{depth}}
\newcommand\codim{\operatorname{codim}}
\newcommand\beg{\operatorname{beg}}
\newcommand\reg{\operatorname{reg}}
\newcommand\Ker{\operatorname{\Ker}}
\newcommand\Coker{\operatorname{Coker}}

\newcommand\Supp{\operatorname{Supp}}
\newcommand\Ass{\operatorname{Ass}}
\newcommand\Ann{\operatorname{Ann}}
\newcommand\Rad{\operatorname{Rad}}
\newcommand\Nor{\operatorname{Nor}}
\newcommand\Sec{\operatorname{Sec}}
\newcommand\NZD{\operatorname{NZD}}
\newcommand\Sing{\operatorname{Sing}}
\newcommand\rank{\operatorname{rank}}

\begin{document}

\title[ALMOST MINIMAL DEGREE]
         {ARITHMETIC PROPERTIES OF PROJECTIVE VARIETIES OF ALMOST MINIMAL
         DEGREE}

\author{Markus BRODMANN and Peter SCHENZEL}

\date{Z\"urich and Halle, June 2005}


\thanks{The second author was partially supported by Swiss National Science
           Foundation (Projects No. 20 - 66980-01 \& No. 20 - 103491 / 1)}


\begin{abstract}
We study the arithmetic properties of projective varieties of
almost minimal degree, that is of non-degenerate irreducible
projective varieties whose degree exceeds the codimension by
precisely $2$. We notably show, that such a variety $X \subset
{\mathbb P}^r$ is either arithmetically normal (and arithmetically
Gorenstein) or a projection of a variety of minimal degree $\tilde
{X} \subset {\mathbb P}^{r + 1}$ from an appropriate point $p \in
{\mathbb P}^{r + 1} \setminus \tilde {X}$. We focus on the latter
situation and study $X$ by means of the projection $\tilde {X}
\rightarrow X$.

If $X$ is not arithmetically Cohen-Macaulay, the homogeneous
coordinate ring $B$ of the projecting variety $\tilde {X}$ is the
endomorphism ring of the canonical module $K(A)$ of the
homogeneous coordinate ring $A$ of $X.$ If $X$ is non-normal and
is maximally Del Pezzo, that is arithmetically Cohen-Macaulay but
not arithmetically normal $B$ is just the graded integral closure
of $A.$ It turns out, that the geometry of the projection $\tilde
{X} \rightarrow X$ is governed by the arithmetic depth of $X$ in
any case.

We study in particular the case in which the projecting variety
$\tilde {X} \subset {\mathbb P}^{r + 1}$ is a (cone over a)
rational normal scroll. In this case $X$ is contained in a variety
of minimal degree $Y \subset {\mathbb P}^r$ such that $\codim_Y(X)
= 1$. We use this to approximate the Betti numbers of $X$.

In addition we present several examples to illustrate our results
and we draw some of the links to Fujita's classification of
polarized varieties of $\Delta $-genus $1$.
\end{abstract}

   \maketitle
   \thispagestyle{empty}

\section{Introduction} \label{1. Introduction}
Let ${\mathbb P}^r_k$ denote the projective $r$-space over an
algebraically closed field $k.$ Let $X \subset {\mathbb P}^r_k$ be
an irreducible non-degenerate projective variety of dimension $d$.
The degree $\deg X$ of $X$ is defined as the number of points of
$X \cap {\mathbb L},$ where $\mathbb L$ is a linear subspace
defined by generically chosen linear forms $\ell _1,\ldots , \ell
_d$. It is well known that
   \[ \deg X \geq \codim X + 1,
   \]
(cf e.g. \cite{EH}), where $\codim X = r - d$ is used to denote
the codimension of $X$. In case equality holds, $X$ is called a
variety of minimal degree. Varieties of minimal degree are
classified and well understood. A variety $X$ of minimal degree is
either a quadric hypersurface, a (cone over a) Veronese surface in
$\mathbb P^5_k,$ or a (cone over a smooth) rational normal scroll
(cf \cite[Theorem 19.9]{H}). In particular these varieties are
arithmetically Cohen-Macaulay and arithmetically normal.

The main subject of the present paper is to investigate varieties
of almost minimal degree, that is irreducible, non-degenerate
projective varieties $X \subset \mathbb P^r_k$ with $\deg X =
\codim X + 2$. From the point of view of polarized varieties,
Fujita \cite{Fu1}, \cite{Fu2}, \cite{Fu3} has studied extensively
such varieties in the framework of varieties of $\Delta $-genus
$1$. Nevertheless, in our investigation we take a purely
arithmetic point of view and study our varieties together with a
fixed embedding in a projective space.

A natural approach to understand a variety $X \subset {\mathbb
P}^r_k$ of almost minimal degree is to view it (if possible) as a
birational projection of a variety of minimal degree $\tilde {X}
\subset {\mathbb P}^{r + 1} _k$ from a point $p \in {\mathbb P}^{r
+ 1}_k \setminus \tilde {X}$. If sufficiently many varieties of
almost minimal degree can be obtained by such projections, we may
apply to them the program of: ``classifying by projections of
classified varieties''. It turns out, that this classification
scheme can indeed be applied to an interesting class of varieties
of almost minimal degree $X \subset {\mathbb P}^r_k$, namely
those, which are not arithmetically normal or equivalently, to all
those which are not simultaneously normal and arithmetically
Gorenstein. More precisely, we shall prove the following result,
in which $\Sec_p(\tilde {X})$ is used to denote secant cone of
$\tilde {X}$ with respect to $p$:

\begin{theorem}
\label{1.1 Theorem} Let $X \subset {\mathbb P}^r_k$ be a
non-degenerate irreducible projective variety and let $t \in \{ 1,
2, \cdots , \dim X + 1 \} $. Then, the following conditions are
equivalent:

\begin{itemize}
\item[{\rm (i)}] $X$ is of almost minimal degree, of arithmetic
depth $t$ and not arithmetically normal.

\item[{\rm (ii)}] $X$ is of almost minimal degree and of
arithmetic depth $t$, where either $t \leq \dim X$ or else $t =
\dim X + 1$ and $X$ is not normal.

\item[{\rm (iii)}] $X$ is of almost minimal degree and of
arithmetic depth $t$, where either $X$ is not normal and $t > 1$
or else $X$ is normal and $t = 1$.

\item[{\rm (iv)}] $X$ is a (birational) projection of a variety
$\tilde {X} \subset {\mathbb P}^{r + 1}_k$ of minimal degree from
a point $p \in {\mathbb P}^{r + 1}_k \setminus \tilde {X}$ such
that $\dim \Sec_p(\tilde {X}) = t - 1$.
\end{itemize}
\end{theorem}

For the proof of this result see Theorem~\ref{5.6 Theorem} (if $t
\leq \dim X$) resp. Theorem~\ref{6.9 Theorem} (if $t = \dim X +
1$). In the spirit of Fujita \cite{Fu1} we say that a variety of
almost minimal degree is {\it maximally Del Pezzo} if it is
arithmetically Cohen-Macaulay (or -- equivalently --
arithmetically Gorenstein). Then, as a consequence of
Theorem~\ref{1.1 Theorem} we have:

\begin{theorem}
\label{1.2 Theorem} A variety $X \subset {\mathbb P}^r_k$ of
almost minimal degree is either maximally Del Pezzo and normal or
a (birational) projection of a variety $\tilde {X} \subset
{\mathbb P}^{r + 1}_k$ of minimal degree from a point $p \in
{\mathbb P}^{r + 1}_k \setminus \tilde {X}$.
\end{theorem}

In this paper, our interest is focussed on those varieties $X
\subset {\mathbb P}^r_k$ of almost minimal degree which are
birational projections of varieties of minimal degree. As already
indicated by Theorem~\ref{1.1 Theorem} and in accordance with our
arithmetic point of view, the arithmetic depth of $X$ is the key
invariant of our investigation. It turns out, that this arithmetic
invariant is in fact closely related to the geometric nature of
our varieties. Namely, the picture sketched in Theorem~\ref{1.1
Theorem} can be completed as follows:

\begin{theorem}
\label{1.3 Theorem} Let $X \subset {\mathbb P}^r_k$ be a variety
of almost minimal degree and of arithmetic depth $t$, such that $X
= \varrho (\tilde {X})$, where $\tilde {X} \subset {\mathbb P}^{r
+ 1}_k$ is a variety of minimal degree and $\varrho :
{\mathbb P}^{r + 1}_k \setminus \{ p \} \rightarrow {\mathbb
P}^r_k$ is a birational projection from a point $p \in {\mathbb
P}^{r + 1}_k \setminus  \tilde{X}. $ Then:

\begin{itemize}
\item[{\rm (a)}] $\nu := \varrho \upharpoonright : \tilde {X}
\rightarrow X$ is the normalization of $X$.

\item[{\rm (b)}] The secant cone $\Sec_p(\tilde {X}) \subset
{\mathbb P}^{r + 1}_k$ is a projective subspace ${\mathbb P}^{t -
1} _k \subset {\mathbb P}^{r + 1}_k$.

\item[{\rm (c)}] The singular locus $\Sing(\nu ) = \varrho
(\Sec_p(\tilde {X}) \setminus \{ p \} ) \subset X$ of $\nu $ is a
projective subspace ${\mathbb P}^{t - 2}_k \subset {\mathbb
P}^r_k$ and coincides with the non-normal locus of $X$.

\item[{\rm (d)}] If $t \leq \dim X, \Sing(\nu )$ coincides with
the non $S_2$-locus and the non-Cohen-Macaulay locus of $X$ and
the generic point of $\Sing(\nu )$ in $X$ is of Goto-type.

\item[{\rm (e)}] The singular fibre $\nu ^{-1}(\Sing(\nu )) =
\Sec_p(\tilde {X}) \cap \tilde {X}$ is a quadric in ${\mathbb
P}^{t - 1}_k = \Sec_p(\tilde {X})$.
\end{itemize}
\end{theorem}

For the proves of these statements see Theorem~\ref{5.6 Theorem}
and Corollary~\ref{6.10 Corollary}.

Clearly, the projecting variety $\tilde {X} \subset {\mathbb P}
^{r + 1}_k$ of minimal degree plays a crucial r\^ole for $X$. We
thus may distinguish the {\it exceptional case} in which $X$ is a
cone over the Veronese surface and the {\it general case} in which
$\tilde {X}$ is a cone over a rational normal scroll. In this
latter case, we have the following crucial result, in which we use
the convention $\dim \emptyset = -1:$

\begin{theorem}
\label{1.4 Theorem} Let $X \subset {\mathbb P}^r_k$ be a variety
of almost minimal degree which is a birational projection of a
(cone over a) rational normal scroll $\tilde {X} \subset {\mathbb
P}^{r + 1} _k$ from a point $p \in {\mathbb P}^{r + 1}_k \setminus
\tilde {X}$. Then, there is a (cone over a) rational normal scroll
$Y \subset {\mathbb P}^r_k$ such that $X \subset Y$ and
$\codim_X(Y) = 1$. Moreover, if the vertex of $\tilde {X}$ has
dimension $h$, the dimension $l $ of the vertex of $Y$ satisfies
$h \leq l \leq h + 3.$ In addition, the arithmetic depth $t$
satisfies $t \leq h +5.$
\end{theorem}

For a proof of this result see Theorem~\ref{Theorem 8.3} and
Corollary~\ref{Corollary 8.4}. It should be noticed, that there
are varieties of almost minimal degree, which cannot occur as a
$1$-codimensional subvariety of a variety of minimal degree (cf
Example~\ref{6.4 Example} and Remark \ref{6.5 Remark}).

Our paper is built up following the idea, that the arithmetic
depth $t:= \depth A$ of a variety $X \subset {\mathbb P} ^r_k =
\Proj(S), S = k [x_0, \cdots , x_r],$ of almost minimal degree
with homogeneous coordinate ring $A = A_X$ is a key invariant. In
Section~\ref{2. Preliminaries} we present a few preliminaries and
discuss the special case where $X$ is a curve.

In Section~\ref{3. The case ``Arithmetic Depth $= 1$''} we
consider the case where $t = 1$. We show that the total ring of
global sections $\oplus _{n \in {\mathbb Z}} H^0(X, {\mathcal
O}_X(n))$ of $X$ -- that is the $S_+$-transform $D(A)$ of $A$ --
is the homogeneous coordinate ring of a variety $\tilde {X}
\subset {\mathbb P}^{r + 1}_k$ of minimal degree. In geometric
terms: $X$ is isomorphic to $\tilde {X}$ by means of a projection
from a generic point $p \in {\mathbb P}^{r + 1}_k$ and hence
normal but not arithmetically normal (cf Propositions~\ref{3.1
Proposition} and \ref{3.4 Proposition}).

In Section~\ref{4. The non-$arithmetically Cohen-Macaulay$ case}
we begin to investigate the case $(1 <) \ t \leq \dim X$, that is
the case in which $X$ is not arithmetically Cohen-Macaulay. First,
we prove some vanishing statements for the cohomology of $X$ and
describe the structure of the $t$-th deviation module $K^t(A)$ of
$A$. Moreover we determine the Hilbert series of $A$ and the
number of defining quadrics of $X$ (cf Theorem~\ref{4.2 Theorem}
and Corollary~\ref{4.4 Corollary}).

In Section~\ref{5. Endomorphism Rings of Canonical Modules} we aim
to describe $X$ as a projection if $t \leq \dim X$. As a
substitute for the $S_+$-transform $D(A)$ of the homogeneous
coordinate ring $A$ (which turned out to be useful in the case $t
= 1$) we now consider the endomorphism ring $B:=
\mbox{End}_A(K(A))$ of the canonical module of $A$ (cf
Theorem~\ref{5.3 Theorem}). It turns out that $B$ is the
homogeneous coordinate ring of variety $\tilde {X} \subset
{\mathbb P}^{r + 1}_k$ of minimal degree, and this allows to
describe $X$ as a projection of $\tilde {X}$ (cf Theorem~\ref{5.6
Theorem}). Endomorphism rings of canonical modules have been
studied extensively in a purely algebraic setting (cf \cite{AGo},
\cite{S}). The striking point is the concrete geometric meaning of
these rings in the case of varieties of almost minimal degree.

In Section~\ref{6. Del Pezzo Varieties and Fujita's
Classification} we study the case where $t = \dim X + 1$, that is
the case where $X$ is arithmetically Cohen-Macaulay. Now, $X$ is a
Del Pezzo variety in the sense of Fujita \cite{Fu3}. According to
our arithmetic point of view we shall speak of maximal Del Pezzo
varieties in order to distinguish them within the larger class of
polarized Del Pezzo varieties. We shall give several equivalent
characterizations of these varieties (cf Theorem~\ref{6.2
Theorem}). We shall in addition introduce the notion of Del Pezzo
variety and show among other things that this notion coincides
with Fujita's definition for the polarized pair $(X, {\mathcal
O}_X(1))$ (cf Theorem~\ref{6.8 Theorem}). Finally we shall prove
that the graded integral closure $B$ of the homogeneous coordinate
ring $A$ of a non-normal maximal Del Pezzo variety $X \subset
{\mathbb P}^r_k$ is the homogeneous coordinate ring of a variety
of minimal degree $\tilde {X} \subset {\mathbb P}^{r + 1} _k$ (cf
Theorem~\ref{6.9 Theorem}) and describe $X$ as a projection of
$\tilde {X}$ (cf Corollary~\ref{6.10 Corollary}). Contrary to the
case in which $t \leq \dim X$, we now cannot characterize $B$ as
the endomorphism ring of the canonical module $K(A)$, simply as
$A$ is a Gorenstein ring. We therefore study $B$ by geometric
arguments, which rely essentially on the fact that we know already
that the non-normal locus of $X$ is a linear subspace (cf
Proposition~\ref{5.8 Proposition}). It should be noticed that on
turn these geometric arguments seem to fail if $t \leq \dim X$.

In Section~\ref{7. Varieties of minimal degree that are
projections} we assume that $X$ is a (birational) projection of a
(cone over a) rational normal scroll $\tilde {X} \subset {\mathbb
P}^{r + 1}_k$. We then prove what is claimed by the previous
Theorem~\ref{1.4 Theorem}. Here, we extensively use the
determinantal description of rational normal scrolls (cf
\cite{H}). As an application we give some constraints on the
arithmetic depth $t$ of $X$ (cf Corollary~\ref{Corollary 8.4} and
Corollary~\ref{Corollary 8.5}).

In Section~\ref{8. Betti numbers} we study the Betti numbers of
the homogeneous coordinate ring $A$ of our variety of almost
minimal degree $X \subset {\mathbb P}^r_k.$ We focus on those
cases, which after all merit a particular interest, that is the
situation where $t \leq \dim X$ and $X$ is a projection of a
rational normal scroll. Using what has been shown in
Section~\ref{7. Varieties of minimal degree that are projections},
we get a fairly good and detailed view on the behaviour of the
requested Betti numbers.

Finally, in Section~\ref{6. Examples} we present various examples
that illustrate the results proven in the previous sections. In
several cases we calculated the Betti numbers of the vanishing
ideal of the occuring varieties on use of the computer algebra
system {\sc Singular} \cite{GrP}.

\section{Preliminaries}
\label{2. Preliminaries}

We first fix a few notation, which we use throughout this paper.
By ${\mathbb N}_0$ (resp. ${\mathbb N})$ we denote the set of
non-negative (resp. positive) integers.

\begin{notation} \label{2.1 Notation}
A) Let $k$ be an algebraically closed field, let $S:= k[x_0,
\cdots , x_r]$ be a polynomial ring, where $r \geq 2$ is an
integer. Let $X \subset {\mathbb P}^r_k = \Proj (S)$ be a reduced
irreducible projective variety of positive dimensions $d$.
Moreover, let ${\mathcal J} = {\mathcal J}_X \subset {\mathcal
O}_{{\mathbb P}^r_k}$ denote the sheaf of vanishing ideals of $X$,
let $I = I_X = \oplus _{n \in {\mathbb Z}} H^0 ({\mathbb P}^r_k,
{\mathcal J}(n)) \subset S$ denote the vanishing ideal of $X$ and
let $A = A_X:= S / I$ denote the homogeneous coordinate ring
of $X$.\\
B) If $M$ is a finitely generated graded $S$-module and if $i \in
{\mathbb Z}$, we use $H^i(M) = H^i_{S_+}(M)$ to denote the $i$-th
{\it local cohomology module} of $M$ with respect to the
irrelevant ideal $S_+ = \oplus _{n \in {\mathbb N}} S_n$ of $S$.
Let $D(M) = D_{S_+}(M)$ denote {\it the $S_+$-transform} $
{\varinjlim} \Hom_S(S^n_+, M)$ of $M$. Moreover, let us introduce
the $i$-th {\it deficiency module} of $M$:
\[
K^i(M) = K^i_S(M):= \Ext^{r + 1 - i}_S( M, S(- r - 1)).
\]
The $S$-modules $D(M), H^i(M)$ and $K^i(M)$ are always furnished
with their natural gradings.
\end{notation}

\begin{reminder} \label{2.2 Reminder}
A) Let $i \in {\mathbb Z}$. If $U = \oplus_{n \in {\mathbb Z}} \
U_n$ is a graded $S$-module, we denote by ${^\ast \Hom}_k(U,k)$
the graded $S$-module $\oplus _{n \in {\mathbb Z}} \
\Hom_k(U_{-n}, k)$. If $M$ is a finitely generated graded
$S$-module, by graded local duality, we have isomorphisms of
graded $S$-modules
\begin{align}
      \label{2.1} K^i(M) &\simeq {^\ast \Hom}_k ( H^i(M), k) \mbox{ and } \\
      \label{2.2} H^i(M) &\simeq {^\ast \Hom}_k ( K^i(M), k) \simeq
      \Hom_S ( K^i(M), E),
\end{align}
where $E$ denotes the graded injective envelope of the $S$-module
$k = S / S_+$.\\
B) By $\depth M$ we denote the depth of the finitely generated
graded $S$-module $M$ (with respect to the irrelevant ideal $S_+$
of $S$), so that
\begin{align}
      \begin{split} \label{2.3}
        \depth M &= \inf \{ i \in {\mathbb Z}
                     \big\arrowvert H^i(M) \not= 0 \} \\
                &= \inf \{ i \in {\mathbb Z} \big\arrowvert K^i
                     (M) \not= 0 \} ,
      \end{split}
\end{align}
(with the usual convention that $\inf \emptyset = \infty $). Here
$\depth A$ is called the {\it arithmetic depth of the variety} $X
\subset {\mathbb P}^r_k$. If we denote the Krull dimension of $M$
by $\dim M$, we have
\begin{align}
      \begin{split} \label{2.4}
        \dim M &= \sup \{ i \in {\mathbb Z}
                     \big\arrowvert H^i(M) \not= 0 \}\\
                &= \sup \{ i \in {\mathbb Z} \big\arrowvert K^i
                     (M) \not= 0 \} ,
      \end{split}
\end{align}
(with the conventions that $\sup \emptyset = - \infty $ and $\dim
0 = - \infty $). \\
C) For a graded $S$-module $U = \oplus_{n \in {\mathbb Z}} \ U_n$,
let $\text{end } U:= \sup \{ n \in {\mathbb Z} \mid U_n \not= 0 \}
$ and $\beg U:= \inf \{ n \in {\mathbb Z} \mid U_n \not= 0 \} $
denote the {\it end } resp. the {\it beginning} of $U$. In these
notation, the {\it Castelnuovo-Mumford regularity} of the finitely
generated graded $S$-module $M$ is defined by
\begin{equation} \label{2.5}
        \reg M = \sup \{ \text{ end } H^i(M) + i \big\arrowvert
                        i \in {\mathbb Z}\}
                 = \inf \{ - \beg K^i(M)) + i \big\arrowvert
                       i \in {\mathbb Z}\} .
\end{equation}
Keep in mind that the Castelnuovo-Mumford regularity of the
variety $X \subset {\mathbb P}^r_k$ is defined as
\begin{equation} \label{2.6}
    \reg X = \reg I = \reg A + 1 .
\end{equation}
\end{reminder}

We are particularly interested in the {\it canonical module} of
$A$, that is in the graded $A$-module
\begin{equation} \label{2.7}
K(A):= K^{\dim (A)}(A) = K^{d + 1}(A) .
\end{equation}

\begin{remark} \label{2.3 Remark}
A) Let $0 < i < \dim (A) = d + 1$ and let $\mathfrak p \in \Spec
S$ with $\dim S / \mathfrak p  = i$. Then, the $S_{\mathfrak p}
$-module $A_{\mathfrak p} $ has positive depth and hence vanishes
or is of projective dimension $< \dim S_{\mathfrak p} = r + 1 -
i$. Therefore $K^i(A)_{\mathfrak p} \simeq \Ext^{r + 1 -
i}(A_{\mathfrak p} ,S_{\mathfrak p} ) = 0$. So
\begin{equation} \label{2.8}
\dim K^i(A) < i \, \mbox { for }\,  0 < i < \dim (A) = d + 1.
\end{equation}\\
B) Let $n \in {\mathbb N}$ and let $f \in A_n \backslash \{ 0 \}.$
Then, $f$ is $A$-regular and the short exact sequence $0 \to A(-n)
\overset {f}{\longrightarrow } A \to A / fA \to 0$ yields an
epimorphism of graded $A$-modules $f : H^{d + 1}(A)(-n)
\twoheadrightarrow H^{d + 1}(A)$. So, by the isomorphisms
\eqref{2.1} of Reminder \ref{2.2 Reminder}, the multiplication map
$f : K^{d + 1}(A) \to K^{d + 1}(A)(n)$ is injective. Moreover,
localizing at the prime ideal $I \subset S$ we get
\[
K^{d + 1}(A) \otimes _A \ \mbox{Quot}(A) \simeq K^{d + 1}(A)_I
\simeq \Ext^{r - d}_{S_I} (S_I / IS_I , S_I) \simeq S_I / IS_I
 = \mbox{Quot }(A).
\]
So, we may resume:
\begin{equation} \label{2.9}
\mbox{The canonical module $K(A)$ of $A$ is torsion free and of
rank $1$.}
\end{equation}\\
C) Let $\ell \in S_1 \backslash \{ 0 \} $ be a linear form. We
write $T:= S / \ell S$ and consider $T$ as a polynomial ring in
$r$ indeterminates. For the $T$-deficiency modules $K^i_T(A / \ell
A)$ of $A / \ell A$, the isomorphisms \eqref{2.1} of Reminder
\ref{2.2 Reminder} together with the base ring independence of
local cohomology furnish the following isomorphisms of graded $A /
\ell A$-modules
\begin{equation*}
K^i_T(A / \ell A) \simeq {^\ast \Hom}_k
                        \left( H^i_{T_+}(A / \ell A), k \right)
                     \simeq {^\ast \Hom}_k \left( K^i_{S_+}(A / \ell
                        A), k \right) \simeq K^i_S(A / \ell A)  .
\end{equation*}
So for all $i \in {\mathbb Z}$ we obtain
\begin{equation} \label{2.10}
K^i_T(A / \ell A) \simeq K^i_S(A / \ell A) = \Ext^{r + 1 - i}_S
\left( A / \ell A , S(- r - 1) \right) .
\end{equation}\\
D) Let $\ell $ be as above. If we apply $\Ext^{r + 1 -
i}_S(\bullet , S(- r - 1))$ to the short exact sequence $0 \to
A(-1) \overset {\ell }{\longrightarrow } A \to A / \ell A \to 0$
and keep in mind the isomorphisms \eqref{2.10}, we get for each $i
\in {\mathbb Z}$ an exact sequence of graded $A / \ell A$-modules
\begin{equation} \label{2.11}
0 \to \left( K^{i + 1}_S(A) / \ell K^{i + 1}_S(A) \right)(1) \to
K^i_T(A / \ell A) \to 0 :_{K^i_S(A)} \ell  \to 0 .
\end{equation}
Correspondingly, applying local cohomology, we get for each $i \in
{\mathbb Z}$ an exact sequence of graded $A / \ell A$-modules
\begin{equation} \label{2.12}
0 \to  H^i_{S_+}(A) / \ell H^i_{S_+}(A) \to H^i_{T_+}(A / \ell A)
\to (0 :_{H^{i + 1}_{S_+}(A)} \ell )(-1) \to 0 .
\end{equation}\\
E) We keep the above notation. In addition, we assume that $\ell
\in S_1 \backslash \{ 0 \} $ is chosen generically. Then,
according to Bertini's Theorem (cf \cite{J}) the hyperplane
section
\begin{align*}
    Y:&= X \cap \Proj (T) = \Proj (T / IT)
            \simeq \Proj (A / \ell A) \subset \Proj (T)
      = {\mathbb P}^{r - 1}_k
\end{align*}
is reduced and irreducible if $\dim A > 2.$ The homogeneous
coordinate ring of $Y$ is
\[
A' = A / (\ell A)^{\mbox{sat}} \simeq T / (I T)^{\mbox{sat}} ,
\]
where $\bullet ^{\mbox{sat}}$ is used to denote the saturation of
a graded ideal in a homogeneous $k$-algebra. Observe that we have
the following isomorphisms of graded $A / \ell A$-modules (cf
\eqref{2.1}, \eqref{2.10}).
\begin{align} \label{2.13}
H^i_{S_+}(A / \ell A) &\simeq H^i_{T_+} (A / \ell A)
                       \simeq H^i_T(A') \,\mbox{ for all }\, i > 0 ; \\
\label{2.14}  K^i_{S}(A / \ell A) &\simeq K^i_{T} (A / \ell A)
                       \simeq K^i_T(A') \,\mbox{ for all }\,i > 0 .
\end{align}
On use of \eqref{2.12} and \eqref{2.13} we now easily get
\begin{equation} \label{2.15}
H^i_{T_+}(A')_{\geq m} = 0 \Rightarrow H^{i + 1}_{S_+} (A)_{\geq m
- 1} = 0 \mbox { for all } i > 0 \mbox { and all } m \in {\mathbb
Z} ,
\end{equation}
where, for a graded $S$-module $U = \oplus _{n \in {\mathbb Z}}
U_n$, we use $U_{\geq n}$ to denote the $m$-th left truncation
$\oplus_{n \geq m} U_n$ of $U$. Finally, if $\depth A > 1$, we
have $A' = A / \ell A$. If $\depth A = 1$, we know that
$H^1_{S_+}(A)$ is a finitely generated non-zero $A$-module so
that, by Nakayama, $\ell H^1_{S_+}(A) \not= H^1_{S_+}(A)$ and
hence $H^1_{T_+}(A') \not= 0$ (cf \eqref{2.12} and \eqref{2.13}).
So, the arithmetic depth of $Y$ behaves as follows
\begin{equation} \label{2.16}
\depth  A' =
    \begin{cases} \depth A - 1 , \ &\mbox{if
                                   depth } A > 1, \\
                                1 , \ &\mbox{if depth } A = 1.
    \end{cases}
\end{equation}
\end{remark}

The aim of the present paper is to investigate the case in which
the degree of $X$ exceeds the codimension of $X$ by $2$. Keep in
mind, that the degree of $X$ always exceeds the codimension of $X$
by 1. Therefore, we make the following convention.

\begin{convention} \label{2.4 Convention}
We write $\dim X, \codim X$ and $\deg X$ for the dimension, the
codimension and the degree of $X$ respectively, so that $d = \dim
X = \dim A - 1, \codim X = \height I = r - \dim X = r - d$. Keep
in mind that
\[
\deg X \geq \codim X + 1
\]
(cf e.g. \cite{EH}). We say that $X$ is of {\it almost minimal
degree}, if $\deg X = \codim X + 2 = r - d + 2.$ Note that $X$ is
called of {\it minimal degree} (cf \cite{EH}) whenever $\deg X =
\codim X + 1.$
\end{convention}

We now discuss the case in which $X$ is a curve of almost minimal
degree.

\begin{remark} \label{2.5 Remark}
A) We keep the hypotheses and notations of Remark \ref{2.3 Remark}
and assume that $\dim X = 1$ and that $\deg X = \codim X + 2 = r +
1$. Then, for a generic linear form $\ell \in S_1 \backslash \{ 0
\} $ and in the notation of part E) of Remark \ref{2.3 Remark},
the generic hyperplane section
\begin{align*}
Y:&= \Proj (T / IT) \simeq \Proj (A / \ell A) = \Proj
                        (A') \subset \Proj (T)
     = {\mathbb P}^{r - 1}_k
\end{align*}
is a scheme of $r + 1$ points in semi-uniform position in
${\mathbb P}^{r - 1}_k$ (cf \cite{Ba}, \cite{N1}). Consequently,
by (cf \cite{GL}) we can say that $IT$ is generated by quadrics.
Therefore we may conclude: The homogeneous $T$-module
\begin{equation}
\label{2.17}
    H^0_{T_+} ( (IT)^{\mbox{sat}} / IT ) \simeq H^0_T(A / \ell A)
    \mbox{ is generated in degree $2$.}
\end{equation}
Moreover, (cf  \cite[(2.4) a)]{BS1})
\begin{equation} \label{2.18}
\dim_k H^1_{T_+} (A / \ell A)_n = \dim _k H^1_{T_+} (A')_n =
    \begin{cases}  r + 1 , \ &\mbox{ if } n < 0 \\
                                r , \ &\mbox{ if } n = 0 \\
                                \leq 1 , \ &\mbox{ if } n = 1 \\
                                0 , \ &\mbox{ if } n > 1
    \end{cases}
\end{equation}
So, by the exact sequences \eqref{2.12} and by statement
\eqref{2.15} we get
\begin{align} \label{2.19}
    H^1(A) / \ell H^1(A) &\subset k(-1) \mbox{ and } \\
     \label{2.20} \text{ end }\, H^2(A) &\leq 0 .
\end{align}\\
B) Assume first, that $A$ is a Cohen-Macaulay ring. Then $H^1(A) =
0$ and so, by \eqref{2.20}, the Hilbert polynomial $P_A(x) \in
{\mathbb Q}[x]$ of $A$ satisfies $P_A(n) = \dim_k A_n$ for all $n
> 0$. As $\dim _k A_1 = r + 1$ it follows $P_A(x) = (r + 1)x$ and
hence $H^2(A)_{-n} \simeq A_n$ for all $n \in {\mathbb Z}$. So, by
\eqref{2.1} $K(A)_n = K^2(A)_n \simeq A_n$ for all $n \in {\mathbb
Z}$. As $K(A)$ is torsion-free of rank $1$ (cf \eqref{2.9}), we
get an isomorphism of graded $A$-modules $K(A) \simeq A(0)$.
Therefore, $A$ is a Gorenstein ring.

If $A$ is normal, $X \subset {\mathbb P}^r_k$ is a smooth
non-degenerate curve of genus $\dim_k K(A)_0 = 1$ and of degree $r
+ 1$, hence an elliptic normal curve: we are in the case
$\bar{I}$ of \cite[(4.7) B)]{BS1}.\\
C) Yet assume that $A$ is a Cohen-Macaulay (and hence a
Gorenstein) ring. Assume that $A$ is not normal. Let $B$ denote
the graded normalization of $A$. Then, there is a short exact
sequence of graded $S$-modules $0 \to A \to B \overset {\pi }{\to
} C \to 0$ with $\dim C = 1$. As $H^0_{S_+}(B) = H^1_{S_+}(B) =
H^0_{S_+}(A) = H^1_{S_+}(A) = 0$ we get $H^0_{S_+}(C) = 0$ and an
exact sequence of graded $S$-modules
\[
0 \to H^1_{S_+}(C) \to H^2(A) \to H^2_{S_+}(B) \to 0 .
\]
As $\dim C = 1$ and $H^0_{S_+}(C) = 0$, there is some $c \in
{\mathbb N}$ such that
\[
\dim_k C_n + \dim _k H^1_{S_+}(C)_n = c \, \mbox{ for all }\,
        n \in {\mathbb Z}.
\]
By \eqref{2.20} and the above sequence $\dim _k C_n = c$ for all
$n > 0$. As $C_0 = 0$ and $\dim _k H^2(A)_0 = \dim _k A_0 = 1$ (cf
part B) ) it follows $c = 1$. As $H^0_{S_+}(C) = 0$, there exits a
$C$-regular element $h \in S_1 \setminus \{ 0 \},$ and choosing $t
\in C_1 \setminus \{ 0 \} $ we get
\begin{equation} \label{2.21}
C = k[h]t \simeq k[h](-1).
\end{equation}
Choose $\bar{y} \in B_1$ such that $\pi (\bar{y}) = t.$ Then we
get $B / A = (\bar{y} A + A) / A$ and hence $B = A[\bar{y}]$. So,
if $y$ is an indeterminate, there is a surjective homomorphism of
homogeneous $k$-algebras
\[
S[y] = k[x_0 , \cdots , x_r, y] \overset {\beta }
        {\twoheadrightarrow } B , \, y \mapsto \bar{y},
\]
which occurs in the commutative diagram
\[
\xymatrix@!R{ S[y] \ar[rr]^-{\beta} && B  \\
           S \ar@{_{(}->}[u] \ar[rr]^-{\alpha } &&A
           \ar@{_{(}->}[u]^-{\alpha }  }
\]
where $\alpha $ is the natural map. Thus, the normalization
$\tilde {X} := \Proj(B)$ of $X$ is a curve of degree $r + 1$ in
$\Proj (S[y]) = {\mathbb P}^{r + 1}_k$ -- a rational normal curve
-- and the normalization morphism $\nu : \tilde {X} \to X$ is
induced by a simple projection $\varrho : {\mathbb P}^{r + 1}
\backslash \{ p \} \twoheadrightarrow {\mathbb P}^r_k$ with center
$\{ p \} = | \Proj (S[y] / S_+ S[y]) |$.

Moreover, by \eqref{2.21} we have $\nu _\ast {\mathcal O}_{\tilde
{X}} / {\mathcal O}_X \simeq \tilde {C} \simeq k$, so that $\nu
_\ast {\mathcal O}_{\tilde {X}} / {\mathcal O}_X$ is supported in
a single point $q \in X$, -- the unique singularity of $X$ -- a
double point. That is, we are in the case $\bar{\mbox{III}}$ of
\cite[(4.7)
B)]{BS1}.\\
D) We keep the notations and hypotheses of part A). But contrary
to what we did in parts B) and C) we now assume that $A$ is not
Cohen-Macaulay, so that $H^1(A) \not= 0$. Then, by \eqref{2.19}
and by Nakayama it follows $H^1(A) / \ell H^1(A) = k(-1)$. In
particular $H^1(A)_1 \simeq k$ and the multiplication map $\ell :
H^1(A)_n \to H^1(A)_{n + 1}$ is surjective for all $n \geq 1$.

Now we claim that $H^1(A)_n = 0$ for all $n > 1$. Assuming the
opposite, we would have an isomorphism $H^1(A)_1 \overset{\ell }
{\to } H^1(A)_2$ and the exact sequence of graded $S$-modules $0
\to H^0_{S_+}(A / \ell A) \to H^1(A)(-1) \overset{\ell }{\to } H^1
(A)$ would imply that $H^0_{T_+}(A / \ell A)_2 \simeq H^0_{S_+} (A
/ \ell A)_2 = 0$ and hence $H^0_{T_+}(A / \ell A) = 0$ (cf
\eqref{2.17}). This would imply $\depth A > 1$, a contradiction.
This proves our claim and shows (cf \eqref{2.1})
\begin{equation} \label{2.22}
H^1(A) \simeq k(-1) \mbox { and } K^1(A) \simeq k(1) .
\end{equation}
By \eqref{2.18} (applied for $n = 1$) it follows that the natural
map $H^1(A)_1 \to H^1(A / \ell A)_1$ is an isomorphism. So
$H^2(A)_0 = 0$. In particular, we get
\begin{equation} \label{2.23}
P_A(x) = (r + 1) x + 1 , \ \text{ end } H^2(A) = - 1 ,
\end{equation}
where $P_A(x) \in {\mathbb Q}[x]$ is used to denote the Hilbert
polynomial of $A$.

As $K(A)$ is torsion-free over the $2$-dimensional domain $A$ (cf
\eqref{2.9}) and satisfies the second Serre property $S_2$ (cf
\cite[3.1.1]{S}), in view of the second statement of \eqref{2.23}
we get :
\begin{equation} \label{2.24}
K(A) \mbox { is a $CM$-module with beg $K(A) = 1$} .
\end{equation}
According to \eqref{2.22}, the $S_+$-transform $D(A)$ of $A$ is a
domain which appears in a short exact sequence $0 \to A \to D(A)
\to k(-1) \to 0$. Choosing $\bar{y} \in D(A)_1 \setminus A_1$ we
obtain $D(A) = A[\bar{y}]$. So, if $y$ is an indeterminate, there
is a surjective homomorphism of homogeneous $k$-algebras $S[y] =
k[x_0, \cdots , x_r, y] \overset{\gamma }{\twoheadrightarrow }
D(A)$, sending $y$ to $\bar{y}$ and extending the natural map
$\alpha : S \to A$ (cf part C) ). In particular $\tilde{X}:= \Proj
(D(A))$ is a curve of degree $r + 1$ in $\Proj (S[y]) = {\mathbb
P}^{r + 1}_k$ -- a rational normal curve. Moreover, the natural
morphism $\varepsilon : \tilde{X} \to X$ is an isomorphism induced
by the simple projection $\varrho : {\mathbb P}^{r + 1}_k
\backslash \{ p \} \twoheadrightarrow {\mathbb P}^r_k$ with $\{ p
\} = |\Proj (S[y]|S_+ S[y]) | \ (\not\in \tilde{X}).$ That is, we
are in the case $\bar{\mbox{II}}$ of \cite[(4.7) B)]{BS1}.
\end{remark}

\section{The case ``Arithmetic Depth $= 1$''}
\label{3. The case ``Arithmetic Depth $= 1$''}

In this section we study varieties of almost minimal degree and
arithmetic depth one. In particular, we shall extend the results
of part D) of Remark~\ref{2.5 Remark} from curves to higher
dimensions.

\begin{proposition} \label{3.1 Proposition} Let $X \subset {\mathbb P}
^r_k$ be a projective variety of almost minimal degree  such that
$\depth  A = 1$ and $\dim X = d$. Then
\begin{itemize}
\item[{\rm (a)}] $H^i(A) = K^i(A) = 0$ for all $i \not= 1, d + 1$;

\item[{\rm (b)}] ${\rm end }\, H^{d + 1}(A) = - \beg K(A)
  = - d $;

\item[{\rm (c)}] $H^1(A) \simeq k(-1) , \ K^1(A) \simeq k(1)$;

\item[{\rm (d)}] $K(A)$ is a torsion-free $CM$-module of rank one;

\item[{\rm (e)}] $D(A)$ is a homogeneous $CM$ integral domain with
$\reg D(A) = 1$ and $\dim_k D(A)_1 = r + 2$.
\end{itemize}
\end{proposition}

\begin{proof} ({\it Induction on} $d = \dim X$). For $d = 1$ all
our claims are clear by the results of part D) of Remark~\ref{2.5
Remark}.

So, let $d > 1$. Let $\ell \in S_1 \backslash \{ 0 \} $ be
generic. Then in the notation of part E) of Remark~\ref{2.3
Remark} we have $\dim A' = d$ and $\depth A' = 1$ (cf
\eqref{2.16}). By induction  $H^i_{T_+}(A') = 0$ for all $i \not=
1, d$. So, by \eqref{2.12} we obtain
\begin{equation}
\label{3.1} H^i(A) = 0 , \mbox { for all } i \not= 1, 2, d + 1 .
\end{equation}
Moreover, by induction and in view of \eqref{2.13} we get
$H^1_{T_+} (A / \ell A) \simeq k(-1)$. As $H^1(A)$ is a non-zero
and finitely generated graded $S$-module, we have $\ell H^1(A)
\not= H^1(A)$. So, by \eqref{2.12} we obtain
\begin{equation}
\label{3.2} H^1(A) / \ell H^1(A) \simeq k(-1)
\end{equation}
and
\begin{equation}
\label{3.3} H^2(A) = 0 .
\end{equation}
Combining \eqref{3.1}, \eqref{3.3} and \eqref{2.1}, we get claim
(a). By induction
\[
\text{ end } H^d(A / \ell A) = \text{ end }
H^d_{T_+} (A') = - d + 1.
\]
As $H^d(A) = 0$, \eqref{2.12} gives
$\text{ end } H^{d + 1}(A) = - d$. In view of \eqref{2.1} we get
claim (b). Also, by induction $\depth K^d_T(A') = d$. As $d > 0$,
we have $H^{d}_{T_+}(A') \simeq H^d_{T_+}(A / \ell A)$ and hence
$K^d _T(A / \ell A) \simeq K^d_T(A')$, (cf \eqref{2.1}). As
$K^d(A) = 0$, \eqref{2.9} and \eqref{2.11} prove statement (d).
Moreover $D(A)$ is a positively graded finite integral extension
domain of $A$ such that $H^1_{S_+}(D(A)) = 0$ and $H^i_{S_+}(D(A))
\simeq H^i(A)$ for all $i > 1$, it follows from statements (a) and
(b), that $D(A)$ is a $CM$-ring with $\reg D(A) = 1$. In view of
\eqref{3.2} and the natural exact sequence
\begin{equation}
\label{3.4} 0 \to A \overset{\eta }{\to } D(A) \overset {\xi }{\to
} H^1(A) \to 0
\end{equation}
there is some $\delta \in D(A)_1 \backslash A$ such that $D(A) = A
+ \delta A$. In particular we have $D(A) = A[\delta ]$ and $D(A)_1
\simeq A \oplus k.$ Therefore statement (e) is proved.

It remains to show statement (c). In view of \eqref{2.1} it
suffices to show that $H^1(A) \simeq k(-1)$. By \eqref{3.2} and as
$H^1(A)_n = 0$ for all $n \gg 0$, there is an isomorphism of
graded $S$-modules $H^1(A) \simeq S / {\mathfrak q}(-1)$, where
${\mathfrak q} \subset S$ is a graded $S_+$-primary ideal. We have
to show that ${\mathfrak q} = S_+$. There is a minimal epimorphism
of graded $S$-modules
\[
\pi : S \oplus S(-1) \to D(A) \to 0
\]
such that $\pi \upharpoonright _S$ coincides with the natural map
$\alpha : S \twoheadrightarrow A$ and $\pi (S(-1)) = \delta A =
\delta S$. As $\beg(\mbox{Ker }(\alpha ) = I) \geq 2$ and $\pi
(S_1) \cap \pi (S(-1)_1) = A_1 \cap \delta k = 0$, it follows
$\beg \mbox{Ker }(\pi ) \geq 2$. Moreover, by statement (e) we
have $\reg D(A) = 1$. Therefore a minimal free presentation of
$D(A)$ has the form
\begin{equation}
\label{3.5} S^\beta (-2) \to S \oplus S(-1) \stackrel
                 {\pi }{\to } D(A) \to 0
\end{equation}
with $\beta \in {\mathbb N}_0$. It follows $\Tor^S_1(k, D(A))
\simeq k^\beta (-2)$. As $1 = \eta (1)$ is a minimal generator of
the $S$-module $D(A)$, the sequence \eqref{3.4} induces an
epimorphism of graded $S$-modules
\[
\Tor^S_1(k, D(A)) \to \Tor^S_1(k, H^1(A)) \to 0 .
\]
Therefore
\[
({\mathfrak q} / S_+ {\mathfrak q})(-1) \simeq \Tor^S_0 ( k,
{\mathfrak q}(-1) ) \simeq \Tor^S_1 ( k,(S / {\mathfrak q}) (-1) )
\simeq \Tor^S_1 ( k, H^1(A) )
\]
is concentrated in degree $2$. So, by Nakayama, ${\mathfrak q}$ is
generated in degree one, thus ${\mathfrak q} = S_+$.
\end{proof}

Varieties of almost minimal degree and arithmetic depth one can be
characterized as simple generic  projections from varieties of
minimal degree.

\begin{reminder}
\label{3.2 Reminder} A) Recall that an irreducible reduced
non-degenerate projective variety $\tilde{X} \subset {\mathbb
P}^s_k$ is said to be of {\it minimal degree} if
$\mbox{deg } \tilde{X} = \mbox{codim } \tilde{X} + 1$.\\
B) Projective varieties of minimal degree are rather well
understood, namely (cf e.g. \cite[Theorem 19.9]{H}): A projective
variety $\tilde{X} \subset {\mathbb P}^s_k$ of minimal degree is
either
\begin{align}
\label{3.6} &\mbox{a quadric hypersurface }, \\
\label{3.7} &\mbox{a (cone over a) Veronese surface in ${\mathbb
                    P}^5_k$ or } \\
\label{3.8} &\mbox{a (cone over a) rational normal scroll}.
\end{align}

C) In particular, a variety $\tilde{X} \subset {\mathbb P}^s_k$ of
minimal degree is arithmetically Cohen-Macaulay and arithmetically
normal.
\end{reminder}

\begin{remark}
\label{3.3 Remark} A) Let $\tilde{X} \subset {\mathbb P}^s_k$ be
an irreducible reduced projective variety, let $p \in {\mathbb P}
^s_k \backslash \tilde{X}$, let $\varrho : {\mathbb P}^s_k
\backslash \{ p \} \twoheadrightarrow {\mathbb P}^{s - 1}_k$ be a
projection with center $p$ and let $X:= \varrho (\tilde{X})
\subset {\mathbb P}^{s - 1}_k$. Then, the induced morphism
$\varrho \upharpoonright : \tilde{X} \twoheadrightarrow X$ is
finite. Moreover, we have $\deg  X = \deg  \tilde{X}$ if and only
if $\varrho \upharpoonright $ is birational, hence if and only if
there is a line $\bar{\ell} \subset {\mathbb P}^s_k$ with $p \in
\bar{\ell }$ and such that the scheme $\bar{\ell } \cap \tilde{X}$
is non-empty, reduced and irreducible. It is equivalent to say
that there are lines $\bar{\ell } \subset {\mathbb P}^s_k$ which
join $p$ and $\tilde{X}$ and are not secant lines of $\tilde{X}$.

But this means precisely that the join $\text{Join} (p, \tilde{X})$
of $p$ and $\tilde{X}$ is not contained in the secant cone
$\Sec_p(\tilde{X})$ of $\tilde{X}$ of $p$. Observe that here
$\Sec_p(\tilde{X})$ is understood as the union of $p$ with all lines
$\bar{\ell } \subset {\mathbb P}^s_k$ such that $p \in \bar{\ell}$ 
and  $\bar{\ell } \cap \tilde{X}$ is a scheme of dimension $0$ and 
of degree $> 1$.

Also, $\varrho \upharpoonright $ is an isomorphism if and only if
for any line $\bar{\ell } \subset {\mathbb P}^s_k$ with $p \in
\bar{\ell }$, the scheme $\bar{\ell } \cap \tilde{X}$ is either
empty or reduced and irreducible. It is equivalent to say that $p
\notin \Sec(\tilde{X})$, where the secant variety
$\Sec(\tilde{X})$ of $\tilde{X}$ is understood as the union of all
lines $\bar{\ell } \subset {\mathbb P}^s_k$ such that $\bar{\ell }
\cap \tilde{X}$ is a scheme of dimension $0$ and of degree $> 1$,
or else $\bar{\ell } \subset \tilde{X}$.

B) Assume now in addition that $\tilde{X} \subset {\mathbb P}^s_k$
is of minimal degree. Then by the above observation we can say:
\begin{equation}
     \label{3.9} X \subset {\mathbb P}^{s - 1}_k \mbox { is of almost
                 minimal degree if and only if $\text{Join}(p,
                 \tilde{X}) \not= \Sec_p(\tilde
                 {X})$.}
\end{equation}
\begin{align}
\begin{split}
\label{3.10} &X \subset {\mathbb P}^{s - 1}_k \mbox { is of
                     almost minimal degree and } \varrho
                     \upharpoonright : \tilde{X} \to
                     X \mbox{an isomorphism }\\
                  &\mbox{if and only if } p \notin
                     \Sec (\tilde{X}), \mbox { thus if and
                     only if } \Sec_p(\tilde{X}) =
                     \{ p \} .
\end{split}
\end{align}
Moreover, by the classification of Reminder \ref{3.2 Reminder} B) 
it follows that $\Sec( \tilde{X}) = \mathbb P^s_k$ whenever $\tilde{X}$ 
is not smooth.
\end{remark}

Now, we can give the announced geometric characterization of
varieties of almost minimal degree and arithmetic depth one.

\begin{proposition}
\label{3.4 Proposition} The following statements are equivalent:
\begin{itemize}
\item[{\rm (i)}] $X$ is of almost minimal degree and of arithmetic
depth $1$.

\item[{\rm (ii)}] $X$ is the projection $\varrho (\tilde{X})$ of a
variety $\tilde{X} \subset {\mathbb P}^{r + 1}_k$ of minimal degree
from a point $p \in {\mathbb P}^{r + 1}_k \backslash  \Sec
(\tilde{X})$.
\end{itemize}
\end{proposition}

\begin{proof} (i) $\Longrightarrow $ (ii): Assume that $X$ is of
almost minimal degree with $\depth A = 1$. Then, by statement e)
Proposition~\ref{3.1 Proposition}, there is some $\bar{y} \in
D(A)_1 \backslash A_1$ such that $D(A) = A[\bar{y}]$. Now, as in
the last paragraph of part D) in Remark~\ref{2.5 Remark}, we may
view $\tilde{X}:= \mbox{ Proj}(D(A))$ as a non-degenerate
irreducible projective variety in ${\mathbb P}^{r + 1}_k$ such
that $\deg(\tilde{X}) = \deg(X)$ and a projection $\varrho :
{\mathbb P}^{r + 1}_k \backslash \{ p \} \twoheadrightarrow
{\mathbb P}^r_k$ from an appropriate point $p \in {\mathbb P}^{r +
1} _k \backslash \tilde{X}$ induces an isomorphism $\varrho
\upharpoonright : \tilde{X} \to X$. In view of Remark \ref{3.3 Remark} B) this
proves statement (ii).

(ii) $\Longrightarrow $ (i): Assume that there is a variety
$\tilde{X} \subset {\mathbb P}^{r + 1}_k$ of minimal degree and a
projection $\varrho : {\mathbb P}^{r + 1}_k \backslash \{ p \}
\twoheadrightarrow {\mathbb P}^r_k$ from a point $p \notin \mbox{
Sec}(\tilde{X})$ with $X = \varrho (\tilde{X})$. Then, by
\eqref{3.10}, $X$ is of almost minimal degree and $\varrho
\upharpoonright : \tilde{X} \to X$ is an isomorphism. It remains to
show that $\depth A = 1$, hence that $H^1 (A) \not= 0$. So, let
$B$ denote the homogeneous coordinate ring of $\tilde{X} \subset
{\mathbb P}^{r + 1}_k$. Then, the isomorphism $\varrho
\upharpoonright : \tilde{X} \to X$ leads to an injective
homomorphism of graded integral domains $A \hookrightarrow B$ such
that $B / A$ is an $A$-module of finite length. Therefore $B
\subset \bigcup _{n \in {\mathbb N}} \ (A :_B S^n _+) = D(A)$. As
$\tilde{X} \subset {\mathbb P}^{r + 1}_k$ is non-degenerate, we have
$\dim_k B_1 = r + 2 > \dim_k A_1$, hence $A \subsetneqq B \subset
D(A)$ so that $H^1(A) \simeq D(A) / A \not= 0$.
\end{proof}

\section{The non-arithmetically Cohen-Macaulay case}
\label{4. The non-$arithmetically Cohen-Macaulay$ case}

In this section we study projective varieties of almost minimal
degree which are not arithmetically Cohen-Macaulay. So, we are
interested in the case where $\deg  X = \codim  X + 2$ and $1 \leq
\depth A \leq \dim X$.

Our first aim is to generalize Proposition~\ref{3.1 Proposition}.
In order to do so, we prove the following auxiliary result, in
which $\NZD_S(M)$ is used to denote the set of non-zero divisors
in $S$ with respect to the $S$-module $M$.

\begin{lemma}
\label{4.1 Lemma} Let $M$ be a finitely generated graded
$S$-module, let $m \in \{ 0, \cdots , r \} $ and $n \in {\mathbb
Z}$. Let $z_m, \cdots , z_r \in S_1$ be linearly independent over
$k$ such that $z_m \in \NZD_s(M)$ and such that there is an
isomorphism of graded $S$-modules $M / z_m M \simeq (S / (z_m,
\cdots , z_r))(n)$.

Then, there are linearly independent elements $y_{m + 1}, \cdots ,
y_r \in S_1$ such that there is an isomorphism of graded
$S$-modules $M \simeq ( S / (y_{m + 1} , \cdots , y_r) ) (n). $
\end{lemma}

\begin{proof} By Nakayama there is an isomorphism of graded
$S$-modules $M \simeq (S / {\mathfrak q})(n)$, where ${\mathfrak
q} \subset S$ is a homogeneous ideal. In particular
\begin{equation}
\label{4.1} z_m \in \NZD_S(S / {\mathfrak q}) .
\end{equation}
As $S / (z_m , {\mathfrak q}) \simeq (M / z_m M)(-n) \simeq S /
(z_m, \cdots , z_r)$, we have
\[
(z_m, {\mathfrak q}) = (z_m, z_{m + 1}, \cdots ,z_r).
\]
Also, by \eqref{4.1}, we have $z_m \notin {\mathfrak q}_1$, so
that ${\mathfrak q}_1$ becomes a $k$-vector space of dimension $r
- m$. Let $y_{m + 1}, \cdots , y_r \in S_1$ form a $k$-basis of
${\mathfrak q}_1$. As
\[
(y_{m + 1}, \cdots ,y_r) \subset {\mathfrak q} \subset (z_m, y_{m
+ 1}, \cdots , y_r)
\]
and in view of \eqref{4.1} we obtain
\[
{\mathfrak q} = (y_{m + 1}, \cdots , y_r) + {\mathfrak q} \cap z_m
S = (y_{m + 1} , \cdots , y_r) + z_m {\mathfrak q} .
\]
So, by Nakayama ${\mathfrak q} = (y_{m + 1}, \cdots , y_r)$.
\end{proof}

Now, we are ready to prove the first main result of this section,
which recover results of \cite{HSV} written down in the context of
the modules of deficiency.

\begin{theorem} \label{4.2 Theorem}
Assume that $X \subset {\mathbb P}^r_k$ is of almost minimal
degree and that $t:= \depth  A \leq \dim X =: d$. Then
\begin{itemize}
\item[(a)] $H^i(A) = K^i(A) = 0$ for all $i \not= t, d + 1$;

\item[(b)] $\mbox{\em end } H^{d + 1}(A) = - \beg K(A) = - d$;

\item[(c)] There are linearly independent forms $y_{t - 1} ,
\cdots , y_r \in S_1$ such that there is an isomorphism of graded
$S$-modules
\[
K^t(A) \simeq ( S / (y_{t - 1} , \cdots , y_r) )(2 - t);
\]

\item[(d)] $K(A)$ is a torsion-free CM-module of rank one.
\end{itemize}
\end{theorem}

\begin{proof} ({\it Induction on $t$}). The case $t = 1$ is clear by
Proposition~\ref{3.1 Proposition}.

So, let $t > 1$ and $\ell \in S_1 \backslash \{ 0 \} $ be generic.
Then, in the notation of Remark~\ref{2.3 Remark} E) we have $A' =
A / \ell A$, $\dim A' = \dim Y + 1 = d$ and $\depth A' = t - 1
\leq d - 1 = \dim Y$ (cf \eqref{2.16}).

(a): By induction, $H^j_{T_+}(A') = 0$ for all $j \not= t - 1, d$.
So, \eqref{2.15} gives $H^i(A) = 0$ for all $i \not= 0, 1, t, d +
1$. As $t > 1$, we have $H^0(A) = H^1(A) = 0$. In view of
\eqref{2.1} this proves statement (a).

(b): By induction $\text{end } H^d_{T_+} (A / \ell A) = - d + 1$.
As $H^d(A) = 0$, \eqref{2.12} implies that $\text{end } $ $H^{d +
1}$ $(A) = - d$ and \eqref{2.1} gives our claim.

(c): By induction there are forms $z_{t - 1}, \cdots , z_r \in
S_1$ whose images $\bar{z}_{t - 1}, \cdots , \bar{z}_r \in T_1$
are linearly independent over $k$ and such that there is an
isomorphism of graded $T$-modules $K^{t - 1}(A / \ell A) = (T /
(\bar{z}_{t - 1}, \cdots , \bar{z}_r)) (2 - (t - 1))$. Let $z_{t -
2}:= \ell $. Then $z_{t - 2}, \cdots , z_r \in S_1$ are linearly
independent and
\[
K^{t - 1}_T(A / \ell A) \simeq ( S / (z_{t - 2}, z_{t - 1}, \cdots
, z_r) ) (3 - t) .
\]
By statement (a) we have $K^{t - 1}(A) = 0.$ So, the sequence
\eqref{2.11} gives an isomorphism of graded $S$-modules
\begin{equation}
\label{4.2} K^t(A) / z_{t - 2} K^t(A) \simeq ( S / (z_{t - 2},
                \cdots , z_r) ) (2 - t) .
\end{equation}
Assume first, that $t < d$. By induction $K^t_T(A / z_{t - 2}A) =
K ^t_T(A')$ vanishes and hence \eqref{2.11} yields $0 :_{K^t(A)}
z_{t - 2} = 0$, thus $z_{t - 2} \in \mbox{NZD}_S(K^t(A))$. So,
\eqref{4.2} and Lemma~\ref{4.1 Lemma} imply statement (c).

Now, let $t = d$. Then \eqref{4.2} implies $\dim K^d(A) / z_{t -
2} K^d (A) = d - 2$. Our first aim is to show that $\dim K^d(A) = d - 1$. 
If $d > 2,$ this follows by the genericity of $z_{t - 2} = \ell .$ 
So, let $t = d = 2.$ Then $A/\ell A$ is a domain of depth 
$1$ which is the coordinate ring of a curve $Y \subset \mathbb P^{r-1}_k$ 
of almost minimal degree (cf Remark \ref{2.5 Remark} A)). So, 
according to Remark \ref{2.5 Remark} D) we have $H^1(A/\ell A) \simeq 
k(-1)$ and $H^2(A/\ell A)_n = 0$ for all $n \geq 0.$ If we apply 
cohomology to the exact sequence $0 \to A(-1) \stackrel{\ell}{\to} A 
\to A/\ell A \to 0$ and keep in mind that $H^1(A) = 0$ we thus get 
$\ell : H^2(A)_{-1} \stackrel{\simeq}{\longrightarrow} H^2(A)_0 \simeq k.$ 
According to Remark \ref{2.5 Remark} D) there is an isomorphism 
$\tilde{Y} \stackrel{\simeq}{\longrightarrow} Y,$ where $\tilde{Y} \subset \mathbb P^{r-1}_k$ 
is a rational normal curve, so that $Y \simeq \mathbb P^1_k$ is smooth. 
As $Y$ is a hyperplane section of $X$ it follows that the non-singular 
locus of $X$ is finite. So, if we apply \cite[Proposition 5.2]{AB} 
to the ample sheaf of $\mathcal O_X$-modules $\mathcal L := \mathcal 
O_X(1)$ and observe that $H^2(A)_n \simeq H^1(X,\mathcal L^{\otimes n})$ 
for all $n \in \mathbb Z,$ we get that $H^2(A)_n \simeq k$ for all 
$n \leq 0.$ Consequently, $K^2(A)_n \not= 0$ for all $n \geq 0,$ 
hence $\dim K^2(A) > 0 = \dim K^2(A)/z_0 K^2(A).$ Therefore 
$\dim K^2(A) = 1,$ which concludes the case $t = d.$ 

According to \eqref{4.2} the
$S$-module $K^d(A) / z_{t - 2} K^d(A)$ is generated by a single
homogeneous element of degree $d - 2$. By Nakayama, $K^d (A)$ has
the same property. So, there is a graded ideal ${\mathfrak q}
\subset S$ with $K^d(A) \simeq (S / {\mathfrak q}) (2 - d)$. In
particular we have $\dim S / {\mathfrak q} = d - 1$.

Now, another use of \eqref{4.2} yields
\begin{multline*}
S / ({\mathfrak q}, z_{d-2}) \simeq (S / {\mathfrak q}) / z_{d -
2}(S / {\mathfrak q}) \simeq K^t(A)(d - 2) / z_{d - 2} K^t(A)(d - 2) \\
\simeq (K^t(A) / z_{d - 2}K^t(A))(d - 2) \simeq S / (z_{t - 2},
\cdots, z_r),
\end{multline*}
so that $({\mathfrak q}, z_{d - 2}) = (z_{t - 2}, \cdots , z_r)$
is a prime ideal. As
\[
\dim S / {\mathfrak q} = d -1 > \dim (S / (z_{t - 2} , \cdots ,
z_r))
\]
it follows, that ${\mathfrak q}$ is a prime ideal. Moreover, as
$z_{d - 2} \notin {\mathfrak q}$, we obtain $z_{d - 2} \in
\mbox{NZD}_S(S / {\mathfrak q}) = \mbox{NZD}_S(K^d(A))$.

Now, our claim follows from \eqref{4.2} and Lemma~\ref{4.1 Lemma}.

(d): In view \eqref{2.9} it remains to show that $\depth K(A) = d
+ 1$. By \eqref{2.9} and by induction we have
\begin{equation}
\label{4.3} \ell \in \NZD_S ( K(A) )
                \mbox { and } \depth K^d_T(A / \ell A) = d .
\end{equation}
So, by the sequence \eqref{2.11}, applied with $i = d$, it
suffices to show that $0 :_{K^d(A)} \ell = 0$. If $t < d$, this
last equality follows from statement (a). If $t = d$, statement
(c) yields $\depth K^d(A) = d - 1 > 0$ and by the genericity of
$\ell $ we get $\ell \in \mbox{NZD}_S(K^d(A))$.
\end{proof}

\begin{remark} \label{4.3 Remark} Keep the notations and hypotheses of
Theorem~\ref{4.2 Theorem}. Then, by statement (c) of Theorem
\ref{4.2 Theorem} and in view of \eqref{2.1} we get $\text{end }
H^t(A) = 2 - t$. So, by statements (a) and (b) of Theorem \ref{4.2
Theorem} we obtain
\begin{equation}
\label{4.5} \reg(A) = 2 \, \text{ and } \, \dim _k A_n = P_A(n) ,
                \, \text{ for all } n > 2 - t .
\end{equation}
\end{remark}

\begin{corollary} \label{4.4 Corollary}
Let $X \subset {\mathbb P}^r_k$ be of almost minimal degree with
$\dim X = d$ and $\depth A = t$. Then:
\begin{itemize}
\item[(a)] The Hilbert series of $A$ is given by
\[
F(\lambda , A) = \frac {1 + (r + 1 - d) \lambda }{(1 -
\lambda)^{d+1}}
       - \frac {\lambda }{(1 - \lambda )^{t - 1}} .
\]

\item[(b)] The Hilbert polynomial of $A$ is given by
\[
P_A(n) = (r - d + 2) \binom{n+d-1}{d}
       + \binom{n + d - 1}{d - 1}
       - \binom{n + t - 2}{t - 2}.
\]

\item[(c)] The number of independent quadrics in $I$ is given by
\[
\dim _k (I_2) = t +\binom{r + 1 - d}{2} - d - 2 .
\]
\end{itemize}
\end{corollary}

\begin{proof} (a): ({\it Induction on} $t$). If $t = 1$, $D(A)$
is a $\mbox{CM}$-module of regularity $1$ (cf Proposition~\ref{3.1
Proposition} (e) ) and of multiplicity $\deg  X = r - d + 2$.
Therefore
\begin{equation}
\label{4.6} F ( \lambda , D(A) ) = \frac {1 + (r + 1 - d)\lambda
}{(1 - \lambda )^{d + 1}} .
\end{equation}
In view of statement c) of Proposition~\ref{3.1 Proposition} we
thus get
\[
F(\lambda , A) = F ( \lambda , D(A) ) - \lambda = \frac {1 + (r +
1 - d) \lambda }{(1 - \lambda )^{d + 1}} - \lambda
\]
and hence our claim.

So, let $t > 1$. Then, as $t' = t -1$ (cf \eqref{2.16}) and $A' =A
/\ell A$ we get by induction
\begin{align*}
F(\lambda , A) =  \frac {F(\lambda , A')} {1 - \lambda } \ &= \
 [\frac {1 + ((r - 1) + 1 - (d - 1))\lambda }{(1 - \lambda )
^d} \ - \ \frac {\lambda }{(1 - \lambda )^{t - 2}}] (1 -
\lambda )^{-1} \\
&= \ \frac {1 + (r + 1 - d)\lambda }{(1 - \lambda )^{d + 1}} -
\frac {\lambda }{(1 - \lambda)^{t - 1}} .
\end{align*}

(b), (c): These are purely arithmetical consequences of statement
(a).
\end{proof}

\begin{remark}
\label{4.7 Remark} Observe that Corollary~\ref{4.4 Corollary} also
holds if $X$ is arithmetically Cohen-Macaulay. In this case, the
shape of the Hilbert series $F(\lambda , A)$ (cf statement (a) )
yields that $A$ is a Gorenstein ring (cf  \cite{St}) which says
that a projective variety of almost minimal degree which is
arithmetically Cohen-Macaulay is already arithmetically
Gorenstein. For $\dim X = 0$ this may be found in \cite{HSV}.

Finally, by Remark~\ref{2.5 Remark} B), by statement \eqref{2.9}
and the exact sequence \eqref{2.12} it follows immediately by
induction on $d = \dim X$ that $K(A) \cong A(1 - d)$ if $X$ is is
arithmetically Cohen-Macaulay. This shows again that $X$ is
arithmetically Gorenstein.
\end{remark}

\section{Endomorphism Rings of Canonical Modules}
\label{5. Endomorphism Rings of Canonical Modules}

Or next aim is to extend the geometric characterization of
Proposition~\ref{3.4 Proposition} to arbitrary non-arithmetically
Cohen-Macaulay varieties of almost minimal degree.

We attack this problem via an analysis of the properties of the
endomorphism ring of the canonical module $K(A)$ of $A$, which in
the local case has been studied already in \cite{AGo}. The crucial
point is, that this ring has a geometric meaning in the context of
varieties of almost minimal degree.

\begin{notation}
\label{5.1 Notation} We write $B$ for the endomorphism ring of the
canonical module of $A$, thus
\[
B:= \Hom_S( K(A), K(A)) .
\]
Observe that $B$ is a finitely generated graded $A$-module and
that
\begin{equation}
\label{5.1} B = \Hom_A( K(A), K(A) ) .
\end{equation}
In addition we have a homomorphism of graded $A$-modules
\begin{equation}
\label{5.2} \varepsilon : A \to B , \; a \mapsto a
                \mbox{ id}_{K(A)} .
\end{equation}
Keep in mind, that $B$ carries a natural structure of (not
necessarily commutative) ring and that $\varepsilon $ is a
homomorphism of rings.
\end{notation}

The homomorphism $\varepsilon : A \to B$ occurs to be of genuine
interest for its own. So we give a few properties of it.

\begin{proposition}
\label{5.2 Proposition} Let $d:= \dim X \geq 1$. Then
\begin{itemize}
\item[(a)] $B = k \oplus B_1 \oplus B_2 \oplus \cdots $ is a
positively graded commutative integral domain of finite type over
$B_0 = k$.

\item[(b)] $\varepsilon : A \to B$ is a finite injective
birational homomorphism of graded rings.

\item[(c)] There is a (unique) injective homomorphism $\tilde
\varepsilon $ of graded rings, which occurs in the commutative
diagram
\[
\xymatrix{A \ar[dr]_-{\varepsilon } \ar@{^{(}->}[rr] && D(A)
       \ar[dl]^-{\tilde \varepsilon } \\
       &B}
\]

\item[(d)] If $\mathfrak p \in \Spec(A)$, the ring $A_{\mathfrak
p} $ has the second Serre property $S_2$ if and only if the
localized map $\varepsilon_{\mathfrak p} : A_{\mathfrak p} \to
B_{\mathfrak p} $ is an isomorphism.

\item[(e)] $\varepsilon : A \to B$ is an isomorphism if and only
if $A$ satisfies $S_2$.

\item[(f)] $\tilde {\varepsilon } : D(A) \to B$ is an isomorphism
if and only if $X$ satisfies $S_2$.

\item[(g)] $B$ satisfies $S_2$ (as an $A$-module and as a ring).

\item[(h)] If the $A$-module $K(A)$ is Cohen-Macaulay, then $B$ is
Cohen-Macaulay (as an $A$-module and as a ring).
\end{itemize}
\end{proposition}

\begin{proof} (a), (b): By \eqref{2.9} (cf Remark~\ref{2.3 Remark})
we know that $K(A)$ is torsion-free and of rank one. From this it
follows easily that $B$ is a commutative integral domain. Also the
map $\varepsilon : A \to B$ is a homomorphism of $A$-modules, and
so becomes injective by the torsion-freeness of the $A$-module
$K(A)$. The intrinsic $A$-module structure on $B$ and the
$A$-module structure induced by $\varepsilon $ are the same. As
$B$ is finitely generated as an $A$-module it follows that
$\varepsilon $ is a finite homomorphism of rings.

It is easy to verify that the natural grading of the $A$-module
$B$ respects the ring structure on $B$ and thus turns $B$ into a
graded ring. In particular $\varepsilon $ becomes a homomorphism
of graded rings. As $A$ is positively graded, $\varepsilon $ is
finite and $B$ is a domain, it follows that $B$ is finite. As $k$
is algebraically closed and $B_0$ is a domain, we get $B_0 \simeq
k$. As $A$ is of finite type over $k$ and $\varepsilon $ is
finite, $B$ is of finite type over $k$, too.

(c): As $\dim A > 1$ we know that $H^1_{S_+}(A)$ is of finite
length. Therefore
\[
\Ext^j_S(H^1_{S_+}(A), S) = 0 \; \mbox{ for all } j \not= r + 1.
\]
So, the short exact sequence $0 \to A \to D(A) \to H^1(A) \to 0$
yields an isomorphism of graded $A$-modules
\[
K(A) = \Ext^{r - d}_S ( A, S(- r - 1))
       \simeq \Ext^{r - d}_S ( D(A), S(- r - 1) ) .
\]
Therefore, $K(A)$ carries a natural structure of graded
$D(A)$-module. As $D(A)$ is a birational extension ring of $A$, we
can write
\[
B = \Hom_A ( K(A), K(A)) = \Hom
       _{D(A)} ( K(A), K(A) )
\]
and hence consider $B$ as a graded $D(A)$-module in a natural way.
In particular, there is a homomorphism of rings
\[
\tilde {\varepsilon }: D(A) \to B , \;
       c \mapsto c \mbox{ id}_{K(A)} ,
\]
the unique homomorphism of rings $\tilde \varepsilon $ which
appears in the commutative diagram
\[
\xymatrix{A \ar[dr]^-{\varepsilon } \ar@{^{(}->}[rr] && D(A)
       \ar[dl]_-{\tilde \varepsilon } \\
       &B}
\]
As $A, D(A)$ and $B$ are domains and as $\varepsilon $ is
injective, $\tilde {\varepsilon }$ is injective, too. Clearly
$\tilde \varepsilon $ is finite, and  respects gradings.

(d): Let $\mathfrak p \in \Spec(A)$. Then, by the chain condition
in $\Spec(S)$, $K(A)_{\mathfrak p} $ is nothing else than the
canonical module $K_{A_{\mathfrak p} }$ of the local domain
$A_{\mathfrak p} $. In particular we may identify $B_{\mathfrak p}
= \Hom_A(K(A), K(A))_{\mathfrak p} \simeq \Hom_{A_{\mathfrak p} }
(K _{A_{\mathfrak p} }, K_{A_{\mathfrak p} })$. Then, the natural
map $\varepsilon_{\mathfrak p} : A_{\mathfrak p} \to B_{\mathfrak
p} $ induced by $\varepsilon $ coincides with the natural map
\[
A_{\mathfrak p} \to \Hom_{A_{\mathfrak p}}(K_{A_{\mathfrak p} },
       K_{A_{\mathfrak p} }), \quad  b \mapsto b \mbox{ id}
       _{K_{A_{\mathfrak p} }} .
\]
But this latter map is an isomorphism if and only if $A_{\mathfrak
p} $ satisfies $S_2$ (cf \cite[3.5.2]{S}).

(e): Is clear by statement (d).

(f): By statement (d), $X$ satisfies $S_2$ if and only if
$\varepsilon _{\mathfrak p} : A_{\mathfrak p} \to B_{\mathfrak p}
$ is an isomorphism for all ${\mathfrak p} \in \Proj (A)$. But
this latter statement is equivalent to the fact that $B /
\varepsilon (A)$ has finite length, thus to the fact that $B
\subset \varepsilon (A) :_B A^n_+$ for some $n$, hence to $B
\subset \tilde {\varepsilon }(D(A))$.

(g): Let ${\mathfrak p} \in \Spec(A)$ of $\height \geq 2$. Then,
the canonical module $K_{A_{\mathfrak p} }$ is of $\depth \geq 2$
(cf \cite[3.1.1]{S}). So, there is a $K_{A_{\mathfrak p}
}$-regular sequence $x, y \in A_{\mathfrak p} $. By the
left-exactness of the functor $\Hom_{A_{\mathfrak p}
}(K_{A_{\mathfrak p} }, \cdot )$ it follows that $x, y$ is a
regular sequence with respect to $\Hom _{A_{\mathfrak p}
}(A_{\mathfrak p} , A_{\mathfrak p} ) = B_{\mathfrak p} $, so that
$\depth_{A_{\mathfrak p} } B_{\mathfrak p} \geq 2$. This shows
that the $A$-module $B$ satisfies $S_2$. As $B$ is finite over
$A$, it satisfies $S_2$ as a ring.

(h): Assume that $K(A)$ is a Cohen-Macaulay module. Consider the
exact sequence of graded $A$-modules $0 \to A \overset
{\varepsilon }{\to } B \to B / A \to 0$. By statement (d) we have
$(B / A) _{\mathfrak p} = 0$ as soon as ${\mathfrak p} \in
\Spec(A)$ is of height $\leq 1$. So $\dim B / A \leq d - 1$ and
$\varepsilon $ induces an isomorphism of graded $S$-modules $H^{d
+ 1}(A) \simeq H^{d+ 1}(B)$. By \eqref{2.1} we get an isomorphism
of graded $S$-modules $K ^{d + 1}(B) \simeq K(A)$. Therefore, the
$A$-module $K^{d + 1}(B)$ is Cohen-Macaulay. In view of
\cite[3.2.3]{S} we thus get $H^i_{A_+} (B) = 0$ for $i = 2, \cdots
, d$. By statement g) we have $H^i_A(B) = 0$ for $i = 0, 1$. So,
the $A$-module $B$ is Cohen-Macaulay. As $B$ is finite over $A$,
it becomes a Cohen-Macaulay ring.
\end{proof}

We now apply the previous result in the case of varieties of
almost minimal degree. We consider $B$ as a graded extension ring
of $A$ by means of $\varepsilon : A \to B$.

\begin{theorem}
\label{5.3 Theorem} Assume that $X \subset {\mathbb P}^r_k$ is of
almost minimal degree and that $t:= \depth A \leq \dim X =: d$.
Then
\begin{itemize}
\item[(a)] $B$ is a finite graded birational integral extension
domain of $A$ and $CM$.

\item[(b)] There are linearly independent linear forms $y_{t - 1},
\cdots , y_r \in S_1$ and an isomorphism of graded $S$-modules $B
/ A \simeq (S / (y_{t - 1}, \cdots , y_r))(-1)$.

\item[(c)] The Hilbert polynomial of $B$ is given by
\[
P_B(x) = (r - d + 2) \binom{x + d - 1}{d} +
       \binom{x + d - 1}{d - 1}.
\]

\item[(d)] If $t = 1$, then $B = D(A)$.
\end{itemize}
\end{theorem}

\begin{proof} (a): According to statement (d) of Theorem
\ref{4.2 Theorem} the $A$-module $K(A)$ is CM. So our claim
follows form statements (a), (b) and (h) of Proposition~\ref{5.2
Proposition}.

(d): Let $t = 1$. According to Proposition~\ref{3.4 Proposition}
and statement \eqref{3.10} of Remark~\ref{3.3 Remark}, $X$ is
isomorphic to a variety $\tilde{X} \subset {\mathbb P}^{r + 1} _k$
of minimal degree and thus CM (cf Reminder~\ref{3.2 Reminder},
part B) ). So, by statement (f) of Proposition~\ref{5.2
Proposition} we have $D(A) = B$.

(b): We proceed by induction on $t$. If $t = 1$, statement (d)
gives $B / A \simeq H^1(A)$ and so we may conclude by statement
(c) of Proposition~\ref{3.1 Proposition}. Let $t > 1$. We write $C
= B / A$ and consider the exact sequence of graded $S$-modules
\[
0 \to A \overset {\varepsilon }{\to }
       B \to C \to 0 .
\]
Let $\ell \in S_1 \backslash \{ 0 \} $ be generic. As $t > 1$ and
$B$ is $\mbox{CM}$ (as an $A$-module) we have $\depth_A(C) > 0$.
Therefore $\ell \in \mbox{NZD}(C)$. We write $A' = A / \ell A$ and
$T = S / \ell S$. Then $A'$ is a domain and $Y:= \Proj (A')
\subset \Proj (T) = {\mathbb P}^{r - 1}_k$ is a variety of almost
minimal degree (cf Remark~\ref{2.3 Remark}).

Let $K(A'):= K^d_T(A'), B':= \Hom_T(K(A'), K(A'))$ and $C' = B' /
A'$. By induction there are linearly independent linear forms
$\bar{z}_{t - 1}, \cdots , \bar{z}_r \in T_1$ and an isomorphism
of graded $T$-modules $B' / A' \simeq (T / (\bar{z}_{t - 1},
\cdots , \bar{z}_r)T)(-1)$. We write $\ell = z_{t - 2}$. Then,
there are linear forms $z_{t - 1}, \cdots , z_r \in S_1$ such that
$z_{t - 2}, z_{t - 1}, \cdots , z_r$ are linearly independent and
such that there is an isomorphism of graded $S$-modules $C' \simeq
(S / (z_{t - 2}, \cdots , z_r)S)(-1)$. As $z_{t - 2} = \ell \in
\mbox{NZD}(C)$, it suffices to show that there is an isomorphism
of graded $S$-modules $C / \ell C \simeq C'$ (cf Lemma~\ref{4.1
Lemma}). As $\ell \in \mbox{NZD}(C)$ there is an exact sequence of
graded $S$-modules
\[
0 \to A' \overset {\alpha }{\to }
       B / \ell B \to C / \ell C \to 0
\]
in which $\alpha $ is induced by $\varepsilon $. It thus suffices
to construct an isomorphism of graded $A$-modules $\bar{\gamma }$,
which occurs in the commutative diagram
\[
\xymatrix{& A'  \ar[dr]^{\varepsilon '} \ar[dl]_-{\alpha} \\
        B / \ell B \ar[rr]_-{\bar{\gamma }} &&B' }
\]
where $\varepsilon '$ is used to denote the natural map. As $A'$
is a domain and as the $A$-module $B$ is $\mbox{CM}$ and torsion
free of rank $1$ (by statement (a)), the $A'$-module $B / \ell B$
is torsion-free and again of rank $1$ (as $\ell \in S_1$ is
generic). By statement (a) the $A'$-module $B'$ is also
torsion-free of rank $1$. So, it suffices to find an epimorphism
$\gamma : B \to B'$, which occurs in the commutative diagram
\[
\xymatrix{ A \ar[rr]^-{\mbox{can}_0} \ar[d]_-{\varepsilon }
                       && A' \ar[d]^-{\varepsilon }   \\
                   B \ar[rr]^-{\gamma } && B' \ar[rr] && 0 }
\]
and hence such that $\gamma (1_B) = 1_{B'}$.

By our choice of $\ell $ and in view of statements (a), (c) and
(d) of Theorem~\ref{4.2 Theorem} we have $\ell \in
\mbox{NZD}(K^d(A)) \cap \mbox{NZD}(K(A))$. So, by \eqref{2.11} of
Remark~\ref{2.3 Remark} we get an exact sequence of graded
$S$-modules
\[
0 \to K(A)(-1) \overset{\ell }{\to } K(A)
       \overset{\pi }{\to } K(A')(1) \to 0 .
\]
Let $U:= \Ext^1_S(K(A), K(A))(-1)$. If we apply the functors
\[
\Hom_S(K(A), \cdot ) \mbox{ and } \Hom_S(\cdot , K(A'))(1)
\]
to the above exact sequence, we get the following diagram of
graded $S$-modules with exact rows and columns
\begin{equation}
    \label{5.3}  \xymatrix{&0 \ar[d] \\
                 &B' = \Hom_S(K(A'), K(A')) \ar[d]^-{\nu } \\
                 \Hom_S(K(A), K(A)) \ar@{=}[d] \ar[r]^-{\mu }
                    & \Hom_S(K(A), K(A'))(-1) \ar[d]^-{0} \ar[r]
                    &  0 :_U   \ell  \ar[r] & 0 \\
                 B  & \Hom_S(K(A), K(A'))}
\end{equation}
where $\mu := \Hom_S(id _{K(A)}, \pi ), \nu := \Hom_S (\pi ,
id_{K(A')})$. With  $\gamma := \nu ^{-1 } \circ \mu $, it follows
\begin{align*}
\begin{split}
\gamma (1_B) = \gamma (id_{K(A)}) = \nu ^{-1}(\mu (id_{K(A)})) =
\nu^{-1}(\pi \circ id_{K(A)}) = \nu ^{-1}(\pi ) \\ = \nu
^{-1}(id_{K(A')(-1)} \circ \pi ) = \nu ^{-1}(\nu (id_{K(A')})) =
id_{K(A')} = 1_{B'}.
\end{split}
\end{align*}

So, it remains to show that $\mu $ is surjective. It suffices to
show that $(0 :_U \ell) = 0$. Assume to the contrary, that $0 :_U
\ell \not= 0$. Then $\ell $ belongs to some associated prime ideal
$\mathfrak p \in \Ass_S U.$ As $\ell$ is generic, this means that
$S_1 \subset {\mathfrak p} $ so that $S_+ \subset {\mathfrak p} $
and hence $0 :_U S_+ \not= 0$, thus $0 :_{(0 :_U \ell )} S_+ \not=
0$. Therefore $\depth_S(0 :_U \ell ) = 0$. In view of statement
(a) we have $\depth_S(B / \ell B) = \depth_S(B') = d$. Moreover
the above diagram \eqref{5.3} yields an exact sequence of graded
$S$-modules
\[
0 \to B / \ell B \to B' \to 0 :_U \ell  \to 0 ,
\]
which shows that $\depth_S(0 :_U \ell) \geq d - 1 > 0$, a
contradiction.

(c): By statement (b) we have $P_{B / A}(x) = \binom{x + t - 3}{t
- 2}$ if $t > 1$ and $P_{B / A}(x) = 0$ if $t = 1$. In view of
statement (a) of Corollary~\ref{4.4 Corollary} we get our claim.
\end{proof}

Now, we are ready to draw a few conclusions about the geometric
aspect.

\begin{notation}
\label{5.4 Notation} A) We convene that ${\mathbb P}^{-1}_k =
\emptyset $ and we use CM$(X)$, $S_2(X)$ and $\Nor(X)$ to denote
respectively the locus of Cohen-Macaulay points, $S_2$-points and
normal points of $X$.

B) If $\nu : \tilde {X} \to X$ is a morphism of schemes, we denote
by $\Sing(\nu )$ the set
\[
\{ x \in X \mid \nu ^\sharp _x : {\mathcal O}_{X,x}
\overset{\not\simeq }{\longrightarrow } (\nu _\ast {\mathcal
O}_{\tilde X})_x \}
\]
of all points $x \in X$ over which $\nu $ is singular.
\end{notation}

\begin{definition}
\label{5.5 Definition} We say that $x \in X$ is a {\it Goto} or
{\it $G$-point}, if the local ring ${\mathcal O}_{X,x}$ is of
``Goto type'' (cf \cite{Go}) thus if $\dim {\mathcal O}_{X,x} > 1$
and
\[
H^i_{{\mathfrak m}_{X,x}}({\mathcal O}_{X,x}) =
       \begin{cases} 0 , &\mbox{if } i \not= 1 , \dim {\mathcal O}
                        _{X,x}, \\
                     \kappa (x), &\mbox{if } i = 1.
       \end{cases}
\]
\end{definition}

\begin{theorem}
\label{5.6 Theorem} Assume that $X \subset {\mathbb P}^r_k$ is of
almost minimal degree and that $t:= \depth A \leq \dim X =: d$.
Then
\begin{itemize}
\item[(a)] $B$ is the homogeneous coordinate ring of a
$d$-dimensional variety $\tilde {X} \subset {\mathbb P}^{r + 1}
_k$ of minimal degree.

\item[(b)] $B$ is the normalization of $A$.

\item[(c)] The normalization $\nu : \tilde X \to X$ given by the
inclusion $\varepsilon : A \to B$ is induced by a projection
$\varrho : {\mathbb P}^{r + 1}_k \setminus \{ p \} \to {\mathbb
P}^r_k$ from a point $p \in {\mathbb P}^{r + 1} _k \setminus
\tilde {X}$.

\item[(d)] The secant cone $\Sec_p(\tilde {X}) \subset {\mathbb
P}^{r + 1}_k$ is a projective subspace of dimension $t - 1$ and
$\Sing  (\nu ) = \varrho (\Sec_p(\tilde {X}) \setminus \{ p \})
\subset X$ is a projective subspace ${\mathbb P}^{t - 2}_k \subset
{\mathbb P}^r_k$.

\item[(e)] The generic point $x \in X$ of $\mbox{ \rm Sing}(\nu )$
is a $G$-point.

\item[(f)] $\Nor(X) = S_2(X) = \mbox{\rm CM}(X) = X \setminus
\mbox{\rm Sing }(\nu )$.
\end{itemize}
\end{theorem}

\begin{proof} (a): By Proposition~\ref{5.2 Proposition} (a), (b)
we see that $B$ is an integral, positively graded $k$-algebra of
finite type and with $\dim B = \dim A = d + 1$. By
Theorem~\ref{5.3 Theorem} (b) and on use of Nakayama we have in
addition $B = k[B_1]$ with $\dim _k B_1 = \dim _k A_1 + 1 = r +
2$. So, $B$ is the homogeneous coordinate ring of a non-degenerate
projective variety $\tilde {X} \subset {\mathbb P}^{r + 1}_k$ of
dimension $d$. By Theorem~\ref{5.3 Theorem} (c) we have $\deg
\tilde {X} = r + 1$.

(b): By statement (a) the ring $B$ is normal (cf Reminder \ref{3.2
Reminder} C) ). In addition, $B$ is a birational integral
extension ring of $A$.

(c): This follows immediately from the fact that $\dim _k B_1 =
\dim_k A_1 + 1 = r + 2$ (cf part C) of Remark~\ref{2.5 Remark}).

(e): According to Theorem~\ref{5.3 Theorem} (b) we have an exact
sequence of graded $S$-modules
\[
0 \to A \to B \to (S / P)(-1) \to 0 ,
\]
where $P:= (y_{t - 1}, \cdots , y_r)S$ with appropriate
independent linear forms $y_{t - 1}, \cdots , y_r \in S_1$. In
particular $0 :_S B / A = P$ and hence $I \subset P$. It follows
that
\[
\Sing(\nu ) = \Supp \Coker(\nu : {\mathcal O}_X \to \nu _\ast
{\mathcal O}_{\tilde {X}}) = \Supp((B / A)^\sim ) = \Proj (S / P)
\]
and that $x:= P / I \in \Proj (A) = X$ is the generic point of
$\Sing(\nu )$. Localizing the above sequence at $x$ we get an
exact sequence of ${\mathcal O}_{X,x}$-modules
\begin{equation}
\label{5.4}
0 \to {\mathcal O}_{X,x} \to
                (\nu _\ast {\mathcal O}_{\tilde {X}})_x
                \to \kappa (x) \to 0
\end{equation}
in which ${\mathcal O}_{X,x}$ has dimension $d - t + 1 > 1$. As
$B$ is a CM-module over $A$ (cf Theorem~\ref{5.3 Theorem} (a) ),
$(\nu _\ast {\mathcal O}_{\tilde X})_x \simeq \tilde {B}_x$ is a
CM-module over the local domain ${\mathcal O}_{X,x}$. So, the
sequence \eqref{5.4} shows that $H^1_{{\mathfrak
m}_{X,x}}({\mathcal O}_{X,x}) \simeq \kappa (x)$ and
$H^i_{{\mathfrak m}_{X,x}}({\mathcal O}_{X,x}) = 0$ for all $i
\not= 1, \dim ({\mathcal O}_{X,x})$.

(d): Let $P \subset S$ be as above. We already know that $\Sing
(\nu ) = \Proj (S / P) = {\mathbb P}^{t - 2}_k \subset {\mathbb
P}^r_k$. Moreover, a closed point $q \in X$ belongs to $\Sing(\nu
)$ if and only if the line $\varrho ^{-1}(q) \subset {\mathbb
P}^{r + 1}_k$ is a secant line of $\tilde {X}$. So
\[
\Sec_p(\tilde {X}) = \{ p \} \cup \varrho ^{-1} (\Sing(\nu )) = \{
p \} \cup \varrho ^{-1} ({\mathbb P}^{t - 1}_k) = {\mathbb P}^{t -
1}_k \subset {\mathbb P}^{r + 1}_k
\]
and $\Sing(\nu ) = \varrho (\Sec_p(\tilde {X}) \setminus \{ p
\})$.

(f): $\nu \upharpoonright : \tilde {X} \setminus \nu ^{-1}
(\Sing(\nu )) \to X \setminus \Sing(\nu )$ is an isomorphism. As
$\tilde {X}$ is a normal CM-variety (cf Reminder~\ref{3.2
Reminder} C) ) it follows that $\Nor(X), \mbox{CM}(X) \supseteq X
\setminus \Sing(\nu )$. As $S_2(X) \supseteq \Nor(X) \cup
\mbox{CM}(X)$ it remains to show that $S_2(X) \subset X \setminus
\Sing(\nu )$. As the generic point $x$ of $\Sing(\nu )$ is not an
$S_2$-point (cf statement (e) ), our claim follows.
\end{proof}

Now, we may extend Proposition~\ref{3.4 Proposition} to arbitrary
non arithmetically Cohen-Macaulay varieties.

\begin{corollary}
\label{5.7 Corollary} Let $1 \leq t \leq \dim (X)$. Then, the
following statements are equivalent:
\begin{itemize}
\item[(i)] $X$ is of almost minimal degree and of arithmetic
 depth $t$.

\item[(ii)] $X$ is the projection $\varrho (\tilde {X})$ of a
variety $\tilde {X} \subset {\mathbb P}^{r + 1}_k$ of minimal
degree from a point $p \in {\mathbb P}^{r + 1}_k \setminus \tilde
{X}$ such that $\dim \Sec_p(\tilde {X}) = t - 1$.
\end{itemize}
\end{corollary}

\begin{proof} (i) $\Longrightarrow $ (ii): Clear by
Theorem~\ref{5.6 Theorem} (d).

(ii) $\Longrightarrow $ (i): Let $\tilde {X} \subset {\mathbb
P}^{r + 1}_k$, $p$ and $\varrho : {\mathbb P}^{r + 1}_k \setminus
\tilde {X} \to X$ be as in statement (ii). Let $\tilde {A}$ be the
homogeneous coordinate ring of $\tilde {X}$. Then $\tilde {A}$ is
normal (cf Reminder~\ref{3.2 Reminder} C) ). Observe that $\varrho
\upharpoonright : \tilde {X} \twoheadrightarrow X$ is a finite
morphism (cf Remark~\ref{3.3 Remark} A) ) so that $\dim \tilde {X}
= d$.

As $\dim \Sec_p(\tilde {X}) = t - 1 < d = \dim \tilde {X} < \dim
\text{Join}(p, \tilde {X})$ there are lines joining $p$ and
$\tilde {X}$ which are not secant lines of $\tilde {X}$. So, in
view of Remark~\ref{3.3 Remark} A) and statement \eqref{3.9} of
Remark~\ref{3.3 Remark} B), the morphism $\varrho \upharpoonright
: \tilde {X} \twoheadrightarrow X$ is birational and $X \subset
{\mathbb P}^r_k$ is of almost minimal degree.

Moreover, the finite birational morphism $\varrho \upharpoonright
: \tilde {X} \twoheadrightarrow X$ is induced by a finite
injective birational homomorphism $\delta : A \hookrightarrow
\tilde {A}$ of graded rings. Thus, by Theorem~\ref{5.6 Theorem}
(b) we get an isomorphism of graded rings $\iota ,$ which occurs
in the commutative diagram
\[
\xymatrix{& A  \ar[dr]^{\varepsilon '} \ar[dl]_-{\delta} \\
        B \ar[rr]^-{\iota }_-{\simeq } &&\tilde {A} }
\]
Now Theorem~\ref{5.6 Theorem} d) shows that $\depth A = \dim
\Sec_p(\tilde{X}) + 1 = t$.
\end{proof}

As a further application of Theorem~\ref{5.6 Theorem} we now have
a glance at arithmetically Cohen-Macaulay varieties of almost
minimal degree and show that their non-normal locus is either
empty or a linear space. More precisely

\begin{proposition}
\label{5.8 Proposition} Assume that $X \subset {\mathbb P}^r_k$ is
of almost minimal degree, $S_2$ and not normal. Let $\dim X =: d$.
Then $X$ is arithmetically Cohen-Macaulay and the non-normal locus
$X \setminus \Nor (X)$ is a linear space ${\mathbb P}^{d - 1}_k
\subset {\mathbb P}^r_k$.
\end{proposition}

\begin{proof} ({\it Induction on $d$}) Let $d = 1$. Then $X$ is
a curve of degree $r + 1$ and thus may have at most one singular
point (cf \cite[(4.7) (B)]{BS1}). So, let $d > 1$. Then $\Nor (X)
\subset X = S_2(X)$ shows that $X$ is arithmetically
Cohen-Macaulay (cf Theorem~\ref{5.6 Theorem} f) ). Moreover, as
$X$ is $S_2$ and not normal, the Serre criterion for normal points
shows that the non-normal locus $X \setminus \Nor(X)$ of $X$ is of
pure codimension $1$.

Let $Z \subset X$ be the reduced and purely $(d -1)$-dimensional
closed subscheme supported by $X \setminus \Nor(X)$. It suffices
to show that $\deg Z = 1$. Now, let $\ell \in S_1$ be a generic
linear form and consider the hyperplane ${\mathbb P}^{r - 1}_k:=
\Proj (S / \ell S)$. Then, the hyperplane section $X':= {\mathbb
P}^{r - 1}_k \cap X = \Proj (A / \ell A)$ is an arithmetically
Cohen-Macaulay-variety of almost minimal degree in ${\mathbb P}^{r
- 1}_k$ with $\dim X' = d - 1$, and $Z':= {\mathbb P}^{r - 1}_k
\cap Z \subset X'$ is a reduced purely $1$-codimensional subscheme
with $\deg Z' = \deg Z$. Now, let $z'$ be one of the generic
points of $Z'$. Then $z' \in Z$ shows that ${\mathcal O}_{X, z'}$
is not normal and hence not regular. As ${\mathcal O}_{X', z'}$ is
a hypersurface ring of ${\mathcal O}_{X,z}$, the ring ${\mathcal
O}_{X', x'}$ is not regular either. As $\dim {\mathcal O}_{X', z'}
= 1$ it follows that ${\mathcal O}_{X', z'}$ is not normal and
hence $z' \in X' \setminus \Nor(X')$. By induction we have $X'
\setminus \Nor(X') = {\mathbb P}^{d - 2}_k$ for some linear
subspace ${\mathbb P}^{d - 2}_k \subset {\mathbb P}^r_k$. It
follows that $\{ \bar{z'} \} = {\mathbb P}^{d - 2}_k$. This shows
that the closed reduced subschemes $Z'$ and ${\mathbb P}^{d -
2}_k$ of ${\mathbb P}^r_k$ coincide, hence $\deg Z = \deg Z' = 1$.
\end{proof}

\section{Del Pezzo Varieties and Fujita's Classification}
\label{6. Del Pezzo Varieties and Fujita's Classification}

In this section we shall treat projective varieties of almost
minimal degree which are arithmetically Cohen-Macaulay.
We call these varieties maximal Del Pezzo varieties and make
sure that this is in coincidence with Fujita's notion of Del
Pezzo variety \cite{Fu3}. We also briefly discuss the link
with Fujita's classification of varieties of $\Delta $-genus
$1$.

\begin{remark}
\label{6.1 Remark} A) Let $d:= \dim (X) > 0$ and let $\omega _X =
K(A)^\sim $ denote the dualizing sheaf of $X$. Keep in mind that a
finitely generated graded $A$-module of $\depth
> 1$ is determined (up to a graded isomorphism) by the sheaf
of ${\mathcal O}_X$-modules $\widetilde {M}$ induced by $M$. So,
as $K(A)$ satisfies the second Serre property $S_2$ (cf
\cite[3.1.1]{S}), we have for each $r \in {\mathbb Z}$:
   \begin{align}
   \label{6.1} &\omega _X \simeq {\mathcal O}_X(r) \mbox{ if and only
               if }
               K(A) \simeq D(A)(r) . \\
   \label{6.2} &\mbox{If depth } A > 1, \mbox{ then } \omega _X
               \simeq  {\mathcal O}_X(r) \mbox{ if and only if }
               K(A) \simeq  A(r) .
   \end{align}

\noindent B) $X \subset {\mathbb P}^r_k$ is said to be {\it
linearly complete} if the inclusion morphism $X \hookrightarrow
{\mathbb P} ^r_k$ is induced by the complete linear system $|
{\mathcal O}_X (1)|$. It is equivalent to say that the natural
monomorphism $\eta : A_1 \rightarrow H^0(X, {\mathcal O}_X(1)) =
D(A)_1$ is an isomorphism hence -- equivalently -- that $H^1(A)_1
= 0$.
\end{remark}

\begin{theorem}
\label{6.2 Theorem} The following statements are equivalent:
\begin{itemize}
\item[{\rm (i)}] $X$ is arithmetically Gorenstein and of almost
minimal degree.
\item[{\rm (ii)}] $X$ is arithmetically Cohen-Macaulay and of almost
minimal degree.
\item[{\rm (iii)}] $X$ is $S_2$, linearly complete and of almost
minimal degree.
\item[{\rm (iv)}] $\omega _X \simeq  {\mathcal O}_X(1 - d)$ and $X$
is of almost minimal degree.
\item[{\rm (v)}] $K(A) \simeq  A(1 - d)$ and $X$ is arithmetically
Cohen-Macaulay.
\item[{\rm (vi)}] $\omega _X \simeq  {\mathcal O}_X(1 - d)$ and
$X$ is arithmetically Cohen-Macaulay.
\item[{\rm (vii)}] $H^{d + 1}(A)_{1 - d} \simeq  k$ and $H^1(A)_1 =
H^i(A)_n = 0$ for $2 \leq i \leq d$ and $1 - d \leq n \leq 1$.
\end{itemize}
\end{theorem}

\begin{proof} The implications (i) $\Longrightarrow $ (ii)
$\Longrightarrow $ (iii) and (v) $\Longrightarrow $ (vi) are
obvious. The implication (ii) $\Longrightarrow $ (i) follows by
Remark \ref{4.7 Remark}, the implication (vi) $\Longrightarrow $
(v) by \eqref{6.2}. It remains to prove the implications (v)
$\Longrightarrow $ (vii) $\Longrightarrow $ (ii) and (ii)
$\Longrightarrow $ (v) $\Longrightarrow $ (iv) $\Longrightarrow $
(iii) $\Longrightarrow $ (ii). The implication (v)
$\Longrightarrow $ (vii) is easy.

(vii) $\Longrightarrow $ (ii): $H^{d + 1}(A)_{1 - d} \simeq  k$
implies that $H^{d + 1}(A)_n = 0$ for all $n > 1 - d$. As $P_A(n)
= \dim _k A_k - \sum ^{d + 1}_{i = 1} (-1)^i \dim _k H^i(A)_n$,
statement (vii) implies that
\[
P_A(n) = \begin{cases} (-1)^d , &\mbox{if } n = 1 - d , \\
      0 , &\mbox{if } 1 - d < n < 0 , \\
      1 , &\mbox{if } n = 0 , \\
      r + 1 , &\mbox{if } n = 1 .  \end {cases}
\]
So, we may write $P_A(n) = (r - d + 2) \binom {n + d - 1}{d} +
\binom {n + d - 2}{d - 2}$. In particular, $X$ is of almost minimal
degree. But then, by Corollary~\ref{4.4 Corollary} b) we see that $X$
must be arithmetically Cohen-Macaulay.

(ii) $\Longrightarrow $ (v): Assume that $X$ is arithmetically Cohen-Macaulay
and of almost
minimal degree. Then, the shape of the Hilbert series given in
Corollary~\ref{4.4 Corollary} (a) allows to conclude that $K(A) \simeq
A(1 - d)$.

(v) $\Longrightarrow $ (iv): We know that (v) $\Longrightarrow$
(vii) and (vii) $\Longrightarrow$ (ii). So, (v) implies that $X$
is of almost minimal degree and hence induces (iv).

(iv) $\Longrightarrow $ (iii): By \eqref{6.2}, statement (iv)
implies $K(A) \simeq  D(A)(1 - d)$ and hence by Theorem~\ref{4.2 Theorem}
(b) that $X$ is arithmetically Cohen-Macaulay.

(iii) $\Longrightarrow $ (ii): Assume that $X$ is of almost
minimal degree, $S_2$ and linearly complete. We have to show that
$X$ is arithmetically Cohen-Macaulay. Assume to the contrary that
$t = \depth A \leq d$. As $H^1(A)_1 = 0$, Proposition~\ref{3.1
Proposition} (c) yields $t > 1$. But now, by Theorem~\ref{5.6
Theorem} (d) and (e), $X$ contains a $G$-point and thus cannot be
$S_2$ -- a contradiction.
\end{proof}

\begin{definition}
\label{6.3 Definition} A) $X \subset {\mathbb P}^r_k$ is called a
{\it maximal Del Pezzo variety} if it satisfies the equivalent
conditions (i) -- (vii) of Theorem~\ref{6.2 Theorem}.\\
B) $X\subset {\mathbb P}^r_k$ is called a {\it Del Pezzo variety},
if there is an integer $r' \geq r$, a maximal Del Pezzo variety
$X' \subset {\mathbb P}^{r'}_k$ and a linear projection $\pi :
{\mathbb P}^{r'}_k \setminus {\mathbb P}^{r' - r - 1}_k
\rightarrow {\mathbb P}^r_k$ with $X' \cap {\mathbb P}^{r' - r -
1} = \emptyset $ and such that $\pi $ gives rise to an isomorphism
$\pi \kern- 3.5 pt \upharpoonright : X' \overset {\simeq
}{\longrightarrow } X$. So, $X \subset {\mathbb P}^r_k$ is Del
Pezzo if and only if it is obtained by a non-singular projection
of a maximal Del Pezzo variety.
\end{definition}

\begin{remark}
\label{6.4 Remark} A) As a linearly complete variety $X \subset
{\mathbb P}^r_k$ cannot be obtained by a proper non-singular
projection of a non-degenerate variety $X' \subset {\mathbb
P}^{r'}_k$ we can say that $X \subset {\mathbb P}^r_k$ is
maximally Del Pezzo if and only if it is Del Pezzo and linearly
complete.\\
B) Keep the previous notation and let $r' = h^0 (X, {\mathcal O}_X
(1)) - 1$, whence
\[
r' = \dim _k D(A)_1 - 1 = r + \dim _k H^1(A)_k.
\]
Moreover let $\varphi : X \rightarrow {\mathbb P}^{r'}_k$ be the
closed immersion defined by the complete linear system $|
{\mathcal O} _X(1)|$ and set $X':= \varphi (X)$. Then $X' \subset
{\mathbb P} ^{r'}_k$ is linearly complete with homogeneous
coordinate ring $A' = k[D(A)_1] \subset D(A)$, whereas the
isomorphism $\varphi : X \overset {\simeq  }{\longrightarrow } X'$
is induced by the inclusion $A \hookrightarrow A'$ and inverse to
an isomorphism ${\pi \kern- 3.9 pt \upharpoonright } : X' \overset
{\simeq  } {\longrightarrow }X$ which is the restriction of a
linear projection $\pi : {\mathbb P}^{r'}_k \setminus {\mathbb
P}^{r' - r - 1}_k \rightarrow {\mathbb P}^r_k$ with $X' \cap
{\mathbb P}^{r' - r - 1}_k = \emptyset $. Now, clearly $X' \subset
{\mathbb P}^{r'}_k$ is linearly complete and moreover
   \begin{align*} D(A') &= D(A), H^i(A') \simeq  H^i(A) 
                    , i \not= 1 ,
                  K(A')\simeq  K(A) \mbox{ and } P_{A'}(x) = P_A(x).
   \end{align*}
Note that $X' \subset {\mathbb P}^{r'}_k$ is called the {\it
linear completion of} $X \subset {\mathbb P}^r_k$. \\
C) Observe that the linear completion of $X \subset {\mathbb P}^r
_k$ is the maximal non-degenerate projective variety $X' \subset
{\mathbb P}^{r'}_k$ which can be projected non-singularly onto
$X.$ More precisely: If $\tilde{X} \subset {\mathbb P}^{\overline
r} _k$ is a non-degenerate projective variety and ${\bar{\pi
}\kern- 3.9 pt \upharpoonright } : \tilde{X} \overset {\simeq
}{\longrightarrow } X$ is an isomorphism induced by a linear
projection $\bar{\pi } : {\mathbb P}^{\overline r}_k
\longrightarrow  {\mathbb P}^r_k$, then $\bar{r} \leq r'$ and the
isomorphism ${\bar{\pi }\kern- 3.9 pt \upharpoonright } ^{-1}
\circ {\pi \kern- 3.9 pt \upharpoonright } : X' \overset {\simeq
}{\longrightarrow } \tilde{X}$ comes from a linear projection
$\varrho : {\mathbb P}^{r'}_k { \longrightarrow } {\mathbb P}
^{\bar{r}}_k$. In particular, if $\tilde{X} \subset {\mathbb
P}^{\bar{r}}_k$ is linearly complete we have $r' = \bar{r}$ and
$\varrho $ becomes an isomorphism so that we may identify $X'$
with $\tilde{X}$. Consequently, by what we said in part A) it
follows that $X \subset {\mathbb P}^r_k$ is Del Pezzo if and only
if its linear completion $X' \subset {\mathbb P}^{r'}_k$ is
(maximally) Del Pezzo.
\end{remark}

We now shall tie the link to Fujita's classification of polarized
Del Pezzo varieties.

\begin{remark}
\label{6.5 Remark} (see \cite{Fu3}). A) A {\it polarized variety over}
$k$ is a pair $(V, {\mathcal L})$ consisting of a reduced irreducible
projective variety $V$ over $k$ and an ample invertible sheaf of
${\mathcal O}_V$-modules ${\mathcal L}$.

B) Let $(V, {\mathcal L})$ be a polarized $k$-variety. For a
coherent sheaf of ${\mathcal O}_V$-modules ${\mathcal F}$ and $i
\in {\mathbb N}_0$ let $h^i(V,{\mathcal F})$ denote the
$k$-dimension of the $i$-th Serre cohomology group
$H^i(V,{\mathcal F})$ of $V$ with coefficients in ${\mathcal F}$.
Then, the function
\[
n \mapsto \chi _{(V,{\mathcal L})} (n):= \sum ^{\dim V}_{i = 0}
      (-1)^i h^i(V, {\mathcal L}^{\otimes n})
\]
is a polynomial of degree $\dim V$, the so called {\it Hilbert
polynomial of the polarized variety} $(V,{\mathcal L})$.

C) Let $(V, {\mathcal L})$ be a polarized variety of dimension $d$.
Then, there are uniquely determined integers $\chi _i(V, {\mathcal L}),
i = 0, \ldots , d$ such that
\[
\chi _{(V,{\mathcal L})}(n) = \sum ^d_{i = 0} \chi _i(V,
      {\mathcal L}) \binom{ n + i - 1}{ i} .
\]
Clearly $\chi _d(V,{\mathcal L}) > 0$. The {\it degree}, the
$\Delta $-{\it genus} and the {\it sectional genus} of the polarized
variety $(V,{\mathcal L})$ are defined respectively by
\begin{align*}
\deg (V,{\mathcal L}):= \ &\chi _d(V,{\mathcal L}) ; \\
\Delta (V,{\mathcal L}):= \ &d + \deg (V,
                    {\mathcal L}) - h^0 (V,{\mathcal L}); \\
g_s(V,{\mathcal L}):= \ &1 - \chi _{d - 1}(V,{\mathcal L}) .
\end{align*}

D) According to Fujita (cf \cite{Fu3}) the polarized variety $(V,
{\mathcal L})$ is
called a {\it Del Pezzo variety}, if it satisfies the following
conditions
    \begin{align}
        \label{6.3} &\Delta (V,{\mathcal L}) = 1 , \\
        \label{6.4} &g_s(V,{\mathcal L}) = 1 , \\
    \label{6.5} &V \mbox{ has only Gorenstein singularities and }
            \omega _V \simeq  {\mathcal L}^{\otimes (1 - \dim V)}, \\
    \label{6.6} &\mbox{For all } i \not= 0, \dim V \mbox{ and all }
    n \in {\mathbb Z} \mbox{ it holds } H^i(V, {\mathcal L}^{\otimes n})
    = 0.
    \end{align}
\end{remark}

\begin{remark}
\label{6.6 Remark} A) We consider the polarized variety
$(X,{\mathcal O} _X(1))$. For all $n \in {\mathbb Z}$ we have $H^0
(X,{\mathcal O} _X(1)^{\otimes n}) = H^0(X, {\mathcal O}_X(n)) =
D(A)_n$. Thus:
   \begin{equation}
   \label{6.7} \chi _{(X,{\mathcal O}_X(1))}(n) = P_A(n),
   \end{equation}
where $P_A$ is the Hilbert polynomial of $A.$  Therefore
   \begin{align}
   \label{6.8} \deg (X,{\mathcal O}_X(1)) &= \deg X , \\
   \label{6.9} \Delta (X,{\mathcal O}_X(1)) &= \deg X - \codim X
                  - 1 - \dim _k H^1(A)_1 .
   \end{align}
As a consequence of the last equality we obtain
\begin{equation}
   \label{6.10}
   \Delta (X,{\mathcal O}_X(1)) \leq \deg X - \codim X + 1,
\end{equation}
with equality if and only if  $X \subset {\mathbb P}^r_k$ is
linearly complete.

B) Let $X' \subset {\mathbb P}^{r'}_k$ be the linear completion of
$X \subset {\mathbb P}^r_k$. Then $(X,{\mathcal O}_X(1))$ and
$(X', {\mathcal O}_{X'}(1))$ are isomorphic polarized varieties.
In particular $(X',{\mathcal O}_{X'}(1))$ is Del Pezzo in the
sense of Fujita if and only $(X, {\mathcal O}_X(1))$ is.
\end{remark}

\begin{lemma}
\label{6.7 Lemma} Let $X \subset {\mathbb P}^r_k$ be of almost
minimal degree. Then
\begin{itemize}
\item[(a)]
\[ \Delta (X,{\mathcal O}_X(1)) = \begin{cases} 0 , &\mbox{ if }
      \mbox{ \rm depth } A = 1 , \\
      1 , &\mbox{ if } \mbox{ \rm depth } A > 1 . \end{cases}
\]
\item[(b)]
\[
g_s (X,{\mathcal O}_X(1)) = \begin{cases} 0 , &\mbox{ if } \ X
      \mbox{ is not arithmetically Cohen-Mcaulay}, \\
      1 , &\mbox{ if } \ X \mbox{ is arithmetically Cohen-Macaulay}. \end{cases}
\]
\end{itemize}
\end{lemma}

\begin{proof} (a): This follows immediately from
Proposition~\ref{3.1 Proposition}
c), Theorem~\ref{4.2 Theorem} b) and by \eqref{6.9}.

(b): This is a consequence of Corollary~\ref{4.4 Corollary} b) and
\eqref{6.7}.
\end{proof}

\begin{theorem}
\label{6.8 Theorem} Let $X \subset {\mathbb P}^r_k$ be of
dimension $d > 0$. Then, the following statements are equivalent:
\begin{itemize}
\item[{\rm (i)}] $X$ is Del Pezzo in the sense of
Definition~\ref{6.3 Definition} B).

\item[{\rm (ii)}] $(X,{\mathcal O}_X (1))$ is Del Pezzo in the sense
of Fujita.

\item[{\rm (iii)}] $\Delta (X, {\mathcal O}_X(1)) = g_s(X,
{\mathcal O}_X(1)) = 1$.

\item[{\rm (iv)}] $\Delta (X,{\mathcal O}_X(1)) = 1$ and $H^i(X,
{\mathcal O}_X(n)) = 0$ for all $i \not= 0, d$ and all $n \in
{\mathbb Z}$.

\item[{\rm (v)}] $g_s(X,{\mathcal O}_X(1)) = 1$ and $H^i(X,{\mathcal
O}_X(n)) = 0$ for all $i \not= 0, d$ and all $n \in {\mathbb Z}$.

\item[{\rm (vi)}] $\Delta (X, {\mathcal O}_X(1)) = 1$ and $\omega _X
\simeq  {\mathcal O}_X (1 - d)$.

\item[{\rm (vii)}] $H^i(X,{\mathcal O}_X(n)) = 0$ for all $ i \not=
0, d$ and all $n \in {\mathbb Z}$, and $\omega _X \simeq  {\mathcal O}
_X(1 - d)$.

\item[{\rm (viii)}] $\Delta (X, {\mathcal O}_X(1)) = 1$ and $X$ is
$S_2$.
\end{itemize}
\end{theorem}

\begin{proof} First let us fix a few notation. Let $r' = h^0 (X,
{\mathcal O}_X(1)) = \dim _k D(A)_1$
and let $X' \subset {\mathbb P}^{r'}_k$ be the linear completion
of $X \subset {\mathbb P}^r_k$.

(i) $\Longrightarrow $ (ii): Let $X$ be Del Pezzo in the sense of
Definition~\ref{6.3 Definition} B). According to Remark~\ref{6.4
Remark} C) we get that $X' \subset {\mathbb P} ^{r'}_k$ is
maximally Del Pezzo. According to Theorem~\ref{6.2 Theorem}, the
pair $(X',{\mathcal O}_{X'}(1))$ thus satisfies the requirements
\eqref{6.3} - \eqref{6.6} of Remark~\ref{6.5 Remark} D) and hence
is Del Pezzo in the sense of Fujita. By Remark~\ref{6.6 Remark} B)
the same follows for $(X,{\mathcal O} _X(1))$.

Clearly statement (ii) implies each of the statements (iii) -
(viii). So, it remains to show that each of the statements (iii) -
(viii) implies statement (i). According to Remark~\ref{6.4 Remark}
C) we may replace $X$ by $X'$ in statement (i). As $(X', {\mathcal
O}_{X'}(1))$ and $(X, {\mathcal O}_X(1))$ are isomorphic polarized
varieties we may replace $X$ by $X'$ in each of the statements
(iii) - (viii). So, we may assume that $X \subset {\mathbb P}^r_k$
is linearly complete.

(iii) $\Longrightarrow $ (i): According to statement ~
\eqref{6.10} of Remark~\ref{6.6 Remark} A) the equality $\Delta
(X,{\mathcal O}_X(1)) = 1$ implies that $X \subset {\mathbb P}
^r_k$ is of almost minimal degree. But now by Lemma~\ref{6.7
Lemma} b) the equality $g_s(X,{\mathcal O}_X(1)) = 1$ implies that
$X$ is arithmetically Cohen-Macaulay.

(iv) $\Longrightarrow $ (i): By $\Delta (X,{\mathcal O}_X(1)) = 1$
we see again that $X \subset {\mathbb P}^r_k$ is of almost minimal
degree. As $H^1 (A)_1 = 0$ we have $\mbox{\rm depth } A > 1$
(Proposition~\ref{3.1 Proposition} C) ). As $H^{i + 1}(A)_n \simeq
H^i(X,{\mathcal O}_X(n)) = 0$ for all $i \not= 0, d$, it follows
that $\mbox{\rm depth } A = d + 1$, so that $X$ is arithmetically
Cohen-Macaulay.

(vi) $\Longrightarrow $ (i): As we have proved the implication
(iii) $\Longrightarrow $ (i) it suffices to prove that statement
(v) implies the equality $\Delta (X,{\mathcal O}_X(1)) = 1$. We
proceed by induction on $d$. Let $d = 1$. As $g_s(X, {\mathcal
O}_X(1)) = 1$ implies $\chi _0(X,{\mathcal O}_X(1)) = 0$ we get
$\chi_{(X,{\mathcal O}_X(1))}(n) = (\deg X)n$. As $H^1(X,{\mathcal
O}_X(1)) = 0$ it follows $r + 1 = h^0 (X, {\mathcal O}_X(1)) =
{\chi}_{(X,{\mathcal O}_X(1))}(1) = \deg X$ and hence by
\eqref{6.10} we get $\Delta (X,{\mathcal O} _X(1)) = 1$.

So, let $d > 1$, let $\ell \in S_1$ be a generic linear form and
consider the irreducible projective variety $Y:= \Proj(A / \ell A)
\subset {\mathbb P}^r_k:= \Proj(S / \ell S)$ of dimension $d - 1$
and with homogeneous coordinate ring $ A' = (A / \ell A) / H^0(A /
\ell A)$. As
\[
\chi _{(Y, {\mathcal O}_Y (1))}(n) = P _{A'}(n) = \Delta P_A(x) =
\Delta \chi _{(X,{\mathcal O}_X(1))}(n)
\]
it follows $\Delta (Y, {\mathcal O}_Y(1)) =
\Delta (X, {\mathcal O}_X(1))$ and $g_s(Y, {\mathcal O}_Y(1))
= g_s(X, {\mathcal O}_X(1)) = 1$. So, by induction it suffices
to show that
\[
H^i(Y, {\mathcal O}_Y(n)) = H^{i + 1}( A')_n = 0 \mbox{ for all }i
\not= 0, d - 1 \mbox{ and all } n \in {\mathbb Z},
\]
and that $Y \subset {\mathbb P}^{r - 1}_k$ is linearly complete,
hence that $H^1(A')_1 = 0$. As
\[
H^{i + 1}(A)_n = H^i(X, {\mathcal
O}_X(n)) = 0
\]
for all $i \not= 0, d$ and all $n \in {\mathbb Z}$
and as $H^1(A)_1 = 0$ this follows immediately if we apply
cohomology to the sequence
\[
0 \rightarrow A(1) \overset{\ell }{\longrightarrow }
A \rightarrow A / \ell A \rightarrow 0
\]
and observe that $H^j(A') \simeq  H^j (A / \ell A)$ for all $j >
0$.

(vi) $\Longrightarrow $ (i): This is immediate by \eqref{6.10} and
Theorem~\ref{6.2 Theorem}.

(vii) $\Longrightarrow $ (i): As $X$ is linearly complete,
$H^1(A)_1 = 0$. Moreover by our hypothesis $H^i(A) = 0$ for
all $i \not= 1, d + 1$. Finally, by \eqref{6.2} we have $K(A)
\simeq  D(A)(1 - d)$, hence $H^{d + 1}(A)_{1 - d} \simeq  K(A)
_{d - 1} \simeq  D(A)_0 \simeq  k$. So, statement (vii) of
Theorem~\ref{6.2 Theorem} is true.

(viii) $\Longrightarrow $ (i): In view of \eqref{6.10},
statement (viii) implies statement (iii) of Theorem~\ref{6.2
Theorem}.
\end{proof}

Our next aim is to extend Theorem~\ref{5.3 Theorem} to maximal
Del Pezzo varieties.

\begin{theorem}
\label{6.9 Theorem} Let $X \subset {\mathbb P}^r_k$ be a maximal
Del Pezzo variety of dimension $d$ which is non-normal. Let $B = k
\oplus B_1 \oplus B_2 \oplus \cdots $ be the graded normalization
of $A$. Then:

\begin{itemize}
\item[(a)] There are linearly independent linear forms $y_d, y_{d
+ 1}, \cdots , y_r \in S_1$ such that $B / A \cong (S / (y_d, y_{d
+ 1}, \cdots , y_r))(-1)$.

\item[(b)] $B$ is the homogeneous coordinate ring of a variety of
minimal degree $\tilde {X} \subset {\mathbb P}^{r + 1}_k$. In
particular, $B$ is a Cohen-Macaulay ring.
\end{itemize}
\end{theorem}

\begin{proof} We make induction on $d$. The case $d = 1$ is clear by
Proposition~\ref{5.2 Proposition} C). Therefore, let $d > 1$.
Statement (b) follows easily from statement (a). So, we only shall
prove this latter. According to Proposition~\ref{5.8 Proposition}
there are linearly independent linear forms $y_d, y_{d + 1},
\cdots , y_r \in S_1$ such that $I \subset (y_d, y_{d + 1}, \cdots
, y_r)$ and such that ${{\mathfrak s}}:= (y_d, y_{d + 1}, \cdots
,y_r) / I \subset S / I = A$ defines the non-normal locus $X
\setminus \Nor(X)$ of $X$. Observe that ${\mathfrak s} \subset A$
is a prime of height $1$ and that $A / {\mathfrak s} \simeq S /
(y_d, y_{d + 1}, \cdots , y_r)S$ is a polynomial ring in $d$
inderminates over $k$. Next, we consider the canonical exact
sequence
   \begin{equation}
   \label{6.11} 0 \rightarrow A \rightarrow B \rightarrow C
                \rightarrow 0
   \end{equation}
in which $C:= B / A$ is a finitely generated graded $A$-module
such that $C_{\leq 0} = 0$ and $\Rad \Ann_B C = {\mathfrak s}$.
Our aim is to show that $C \simeq (A / {\mathfrak s}) (-1)$. Let
$\ell \in S_1$ be a generic linear form.

Then, according to Bertini and as $\depth A > 1$, the ideal $\ell
A \subset A$ is prime. Moreover $(A \setminus \ell A) \cap
{\mathfrak s} \not= \emptyset $, so that $(A \setminus \ell
A)^{-1} B = A_{\ell A}$. It follows that $\ell B$ has a unique
minimal prime $\mathfrak p $ and that $\mathfrak p B_\mathfrak p =
\ell B_\mathfrak p $. As $B$ is $S_2$ we get $\ell B = \mathfrak p
$, so that $\ell B$ is a prime ideal of $B$. Therefore $B / \ell
B$ is a finite birational integral extension domain of $A':= A /
\ell A$ and hence a subring of the graded normalization $B'$ of
$A'$.

As $\depth A > 2$ and $\depth  B \geq 2,$ the short exact sequence
\eqref{6.11} yields $\depth  C \geq 2$. In particular $\ell $ is
$C$-regular. Hence we get the following commutative diagram with
exact rows and columns in which $U = \Coker \iota$ is a graded
$A$-module
\begin{align}
\begin{split}
    \label{6.12}  \xymatrix{&&0 \ar[d] &0 \ar[d] \\
                  0 \ar[r] &A / \ell A \ar@{=}[d] \ar[r] &B / \ell B
                    \ar[d]^\lambda \ar[r] &C / \ell C \ar[d]^\iota
                    \ar[r] &0 \\
                  0 \ar[r] &A' \ar[r] &B'
                    \ar[d] \ar[r] &C' \ar[d] \ar[r] &0 \\
                  &&U \ar[d] \ar[r]^-{\simeq} &U \ar[d] \\
                  &&0 &0}
   \end{split}
   \end{align}
Now $X':= \Proj(A') \subset \Proj(S / \ell S) = {\mathbb P}^{r -
1}_k$ is again a maximal Del Pezzo variety. Moreover $X'$ is
non-normal, since otherwise $B / \ell B = A / \ell A,$ hence $B =
A$. Let ${\mathfrak s}' \subset A'$ be the prime of height $1$
which defines the non-normal locus of $X'$ and keep in mind that
$A' / {\mathfrak s}'$ is a polynomial ring in $d - 1$ inderminates
over $k$. By induction $C' \simeq (A' / {\mathfrak s}')(-1).$

Our next aim is to show that $\dim U \leq 0$. As $C \not= 0$, we
have $C / \ell C \not= 0$. As $C'$ is a free $A' / {\mathfrak
s}'$-module of rank one, it follows $\dim U < \dim (A' /
{\mathfrak s}') = d - 1$. As $\dim (B / \ell B) = \dim (B') = d$
it follows that $\lambda $ is an isomorphism in codimension one.
As $B'$ is normal and hence satisfies the property $R_1$, it
follows that $B / \ell B$ satisfies $R_1$, too.

Let $s \in {\mathfrak s} \setminus \{ 0 \}.$ As $B$ is normal, it
satisfies the Serre property $S_2$ so that $B / sB$ satisfies
$S_1.$ Therefore the set of generic points of the (closed)
non-$S_2$-locus of the $B$-module $B / sB$ is given by
\[
\mathcal Q := \{ {\mathfrak q} \in V(sB) : \depth (B_{\mathfrak q}
/ s B_{\mathfrak q}) = 1 < \height {\mathfrak q} - 1 \}.
\]
In particular $\mathcal Q$ is finite, and hence $\ell $ avoids all
members of ${\mathcal P}:= (\Ass_B(B / sB) \cup {\mathcal Q}) \cap
\Proj B.$

Now, let $\mathfrak r \in \Proj B \cap V(\ell B)$ such that
$\height \mathfrak r > 2.$ If $s \notin \mathfrak r,$ the equality
$A_s = B_s$ yields that $B_{\mathfrak r}$ is a Cohen-Macaulay
ring, so that $\depth (B_{\mathfrak r} / \ell B_{\mathfrak r}) >
1.$ If $s \in \mathfrak r$ the fact that $\ell $ avoids all
members of ${\mathcal P}$ implies that $s, \ell $ is a
$B_{\mathfrak r}$-sequence and $\depth B_{\mathfrak r} / s
B_{\mathfrak r}) > 1.$ It follows again that $\depth (B_{\mathfrak
r} / \ell B_{\mathfrak r}) > 1.$ This proves, that the scheme
$\Proj (B / \ell B)$ is $S_2.$ As $B / \ell B$ satisfies $R_1$ it
follows that $\Proj (B / \ell B)$ is a normal scheme, hence that
$\Proj (B / \ell B) = \Proj B'.$ As a consequence, we get indeed
that $\dim U \leq 0,$ that is $U$ is a graded $A$-module of finite
length.

Now, let ${\mathfrak t} \subset A$ be the preimage of ${\mathfrak
s}'$ under the canonical map $A \rightarrow A'$. Then ${\mathfrak
t}$ and ${\mathfrak s} + \ell A$ are primes of height  2 in $A$
and so
\[
\mathfrak s + \ell A = \Rad ((\Ann_A C) + \ell A) = \Rad (\Ann_A
C/\ell C) \supseteq \Ann_A C' = \Ann_A( A' / {\mathfrak s}'(-1)) =
{\mathfrak t}
\]
implies that ${\mathfrak s} + \ell A = \Ann_A (C / \ell C) =
{\mathfrak t}.$ As a consequence we get ${\mathfrak s}C \subseteq
\ell C$ and hence, by the genericity of $\ell $, that ${\mathfrak
s} C = 0$. It follows ${\mathfrak s}B \subseteq {\mathfrak s}$ and
${\mathfrak s}$ becomes an ideal of $B$. Now, let $a \in A$ and $c
\in C \setminus \{ 0 \}$ such that $ac = 0$. By the genericity of
$\ell $ we may assume that $c \notin \ell C$ so that $\iota (c +
\ell C) \not= 0$ and $a\iota (c + \ell C) = 0.$ It follows $a \in
{\mathfrak t} = {\mathfrak s} + \ell A$ and hence, by genericity,
that $a \in {\mathfrak s}.$ This shows that $C$ is a torsion-free
$A / {\mathfrak s}$-module and hence that $B / {\mathfrak s}$ is a
torsion-free $A / {\mathfrak s}$-module.

As $\rank_{A / {\mathfrak s}} C = e_0(C) = e_0(C / \ell C) = e_0
(C') = \rank_{A' / {\mathfrak s}'} C' = 1$ we get an exact
sequence of graded $A / {\mathfrak s}$-modules
   \begin{equation}
   \label{6.13} 0 \rightarrow C \rightarrow (A / {\mathfrak s})(-m)
                \rightarrow W \rightarrow 0
   \end{equation}
with $m \in {\mathbb Z}$ and $\dim  W < \dim  A / {\mathfrak s} =
d.$ We choose $m$ maximally. Then, there is no homogeneous element
$f \in A / {\mathfrak s}$ of positive degree with $C(m) \subseteq
f(A / {\mathfrak s})$, so that $\dim W = \dim(W(m)) = \dim((A /
{\mathfrak s}) / C(m)) < d - 1.$

As $\depth  C \geq 2$ we have $\depth W \geq 1.$ As $\ell $ is
generic we thus get an exact sequence
   \[ 0 \rightarrow C / \ell C \rightarrow (A' / {\mathfrak s}')(-m)
      \rightarrow W / \ell W \rightarrow 0
   \]
with $\dim W / \ell W = \dim  W - 1 < d - 2 = \dim A' / {\mathfrak
s}' - 1.$ Comparing Hilbert coefficients we get $e_1(C / \ell C) =
e_1((A' / {\mathfrak s}')(-m)) = m.$

The diagram \eqref{6.12} contains the short exact sequence $0
\rightarrow C / \ell C \rightarrow (A' / {\mathfrak s}')(-1)
\rightarrow U \rightarrow 0$ with $\dim U \leq 0.$ Assume first,
that $d \geq 3$ so that $\dim  U < \dim A' / {\mathfrak s}' - 1.$
Then, we may again compare Hilbert coefficients and get $e_1(C /
\ell C) = e_1((A' / {\mathfrak s}')(-1)) = 1$, hence $m = 1.$ It
follows $\dim W / \ell W = \dim U \leq 0,$ thus $\dim W \leq 1.$
Suppose that $\dim  W = 1.$ Then $H^1(W)_n \not= 0$ for all $n \ll
0$ and the sequence \eqref{6.11} yields $H^2(C)_n \not= 0$ for all
$n \ll 0.$ As $H^3(A) = 0$, the sequence \eqref{6.11} induces that
$H^2(B)_n \not= 0$ for all $n \ll 0;$ but this contradicts the
fact that $B$ is $S_2$. So, we have $\dim W \leq 0.$ Now, by the
sequence \eqref{6.13} we get $W = H^0(W) \simeq H^1(C),$ and
$\depth C > 1$ implies $W = 0.$ Therefore $C \simeq A / {\mathfrak
s}(-1)$.

It remains to treat the case $d = 2$. Now $A' / {\mathfrak s}'
\simeq k [x]$ so that the $A' / {\mathfrak s}'$-submodule $C /
\ell C$ of $(A' / {\mathfrak s}')(-1)$ is generated by a single
homogeneous element of degree $m \geq 1.$ By Nakayama it follows
$C \simeq (A / {\mathfrak s})(-m).$ It remains to show that $m =
1.$

By the case $d = 1$ we know that $B'$ is the homogeneous
coordinate ring of a rational normal curve, so that $H^2(B')_0 =
0.$ The middle column of the diagram \eqref{6.12} now implies
$H^2(B / \ell B)_0 = 0$. Applying cohomology to the exact sequence
$0 \rightarrow B(-1) \overset {\ell }{\rightarrow } B \rightarrow
B / \ell B \rightarrow 0$ we thus get an isomorphism $H^3(B)_{-1}
\simeq H^3(B)_0$ so that $H^3 (B)_{-1} = 0.$ If we apply
cohomology to the sequence \eqref{6.11} we thus get an exact
sequence
\[
0 \rightarrow H^2(B)_{-1} \rightarrow H^2(A / {\mathfrak
s})_{-1-m} \rightarrow H^3(A)_{-1} \rightarrow 0.
\]
By Theorem~\ref{6.2 Theorem} (cf statement (vii)) we have $H^3(A)
_{-1} \simeq k$. As $A / {\mathfrak s}$ is a polynomial ring in
two indeterminates over $k$ we have $H^2(A / {\mathfrak s})_{-1 -
m} \simeq k^m.$ It follows $H^2(B)_{-1} \simeq k^{m-1}.$

Moreover $\tilde {X}:= \Proj B$ is a projective normal surface and
the natural morphism $\nu : \tilde {X} \rightarrow X$ is a
normalization of $X.$ In particular ${\mathcal L}:= B(1)^\sim =
\nu ^\ast {\mathcal O}_X(1)$ is an ample invertible sheaf of
${\mathcal O} _{\tilde {X}}$-modules and ${\mathcal L}^{\otimes n}
= \nu ^\ast {\mathcal O}_X (n) = B(n)^\sim $ for all $n \in
{\mathbb Z}.$ In addition we have $H^2(B)_{-1} \simeq H^1(\tilde
{X}, {\mathcal L} ^{\otimes - 1})$ and $B_1 \simeq H^0(\tilde {X},
{\mathcal L}).$ Moreover $Y:= \Proj B / \ell B \simeq \Proj B'$ is
the effective divisor on $\tilde {X}$ defined by the global
section $\ell \in H^0(\tilde {X}, {\mathcal L}) \backslash \{ 0
\}. $ As $H^1(Y,{\mathcal O}_Y) = H^2(B / \ell B) = 0,$ the
sectional genus $g_s(\tilde {X},{\mathcal L})$ vanishes (cf
\cite[(5.3) B)]{AB}). As $\dim_k H^0(\tilde {X},{\mathcal L}) =
\dim_k B_1 \geq r + 1 > 1$ it follows $H^1(\tilde {X}, {\mathcal
L}^{\otimes - 1}) = 0$ (cf \cite[Proposition (5.4)]{AB}).
Therefore $k^{m - 1} \simeq H^2 (B)_{-1} = 0$, hence $m = 1.$
\end{proof}

Now, we may extend Theorem~\ref{5.6 Theorem} as follows

\begin{corollary}
\label{6.10 Corollary} Let $X \subset {\mathbb P}^r_k$ be of
almost minimal degree. Assume that either $t:= \depth  A \leq \dim
 X =: d$ or that $X$ is maximally Del Pezzo (that is $t = d + 1)$
and non-normal. Then, there is a $d$-dimensional variety $\tilde
{X} \subset {\mathbb P}^{r + 1}_k$ of minimal degree, a point $p
\in {\mathbb P}^{r + 1}_k \setminus \tilde {X}$ and a projection
$\varrho : {\mathbb P}^{r + 1}_k \setminus \{ p \} \rightarrow
{\mathbb P}^r_k$ from $p$ such that:

\begin{itemize}
\item[(a)] $\varrho (\tilde {X}) = X$ and $\nu : = \varrho
\upharpoonright : \tilde {X} \rightarrow X$ is the normalization
of $X.$

\item[(b)] The secant cone $\Sec_p(\tilde {X}) \subset {\mathbb
P}^{r + 1}_k$ is a projective subspace ${\mathbb P}^{t - 1} _k
\subset {\mathbb P}^{r + 1}_k$ and
\[
X \setminus \Nor X = \Sing (\nu ) = \varrho (\Sec_p(\tilde {X}
\setminus \{ p \} ) \subset X
\] is a projective subspace ${\mathbb P}^{t - 2}_k \subset
{\mathbb P}^r_k.$

\item[(c)] The singular fibre $\nu ^{-1}(\Sing (\nu )) =
\Sec_p(\tilde{X}) \cap \tilde{X} \subset \tilde {X}$ is  a quadric
in ${\mathbb P}^{t - 1}_k = \Sec_p(\tilde {X}).$
\end{itemize}
\end{corollary}

\begin{proof} Let $B$ be the graded normalization of $A$. Then,
according to Theorem~\ref{5.6 Theorem} resp. Theorem~\ref{6.9
Theorem} we see that $B$ is the homogeneous coordinate ring of a
variety $\tilde {X} \subset {\mathbb P}^{r + 1}_k$ of minimal
degree. Moreover, by Theorem~\ref{5.3 Theorem} resp.
Theorem~\ref{6.9 Theorem} there are linearly independent linear
forms $y_{t - 1}, y_t, \cdots , y_r \in S_1$ such that $I
\subseteq (y_{t - 1}, y_t, \cdots , y_r)$ and $B / A \simeq (S /
(y_{t - 1}, y_t, \cdots , y_r))(-1)$. Now, all claims except
statement (c) follow as in Theorem~\ref{5.6 Theorem}.

To prove statement (c), we consider the prime ${\mathfrak s} :=
(y_{t - 1}, y_t, \ldots , y_r) / I \subset A.$ Then $A /
{\mathfrak s}$ is a polynomial ring in $t - 1$ indeterminates over
$k$ and we have $B / A \simeq (A / {\mathfrak s})(-1).$ In
particular ${\mathfrak s} \subset B$ is an ideal and we have an
exact sequence $0 \rightarrow A / {\mathfrak s} \rightarrow B /
{\mathfrak s} \rightarrow (A / {\mathfrak s})(-1) \rightarrow 0$.
Therefore $B / {\mathfrak s} \simeq (A / {\mathfrak s})[z] / (f)$
for some polynomial $f = z^2 + uz + v$ with $u \in (A / {\mathfrak
s})_1$ and $v \in (A / {\mathfrak s})_2.$ As
\[
\nu ^{-1}(\Sing (\nu) ) = \Proj (B / {\mathfrak s}) \subset
\Sec_p(\tilde {X}) = {\mathbb P}^{t - 1}_k = \Proj ((A /
{\mathfrak s})[z])
\]
the claims of (c) follows.
\end{proof}

\section{Varieties of almost minimal degree that are projections}
\label{7. Varieties of minimal degree that are projections}

We now
wish to give a more detailed insight in the nature of those
varieties of almost minimal degree which are projections of cones
over rational normal scrolls. Let us first recall a few facts on
such scrolls.

\begin{remark} \label{Remark 8.1} (cf \cite[pp 94, 97,
108-110]{H} and \cite{ES}) A) Let $n \in \mathbb N$ and let $l,
a_1,\ldots,a_{l-1} \in \mathbb N.$  Let $a_0 = -1, a_l = n$ and
assume that $a_i - a_{i-1} > 1$ for $i = 1,\ldots,l.$ Then, up to
projective equivalence, the numbers $a_1, a_2, \ldots, a_{l-1}$
define a unique {\sl rational normal $l$-fold scroll}
$S_{a_1\cdots a_{l-1}} \subset \mathbb P^n_k.$ Keep in mind that
$S_{a_1 \cdots a_{l-1}}$ is smooth, rational,
arithmetically Cohen-Macaulay and of dimension $l.$ \\
B) Keep the notation of part A). After an appropriate linear
coordinate transformation we may assume that the vanishing ideal
of $S_{a_1 \cdots a_{l-1}}$ in the polynomial ring
$k[x_0,\ldots,x_r]$ is the ideal generated by the $2 \times
2$-minors of the $2 \times (n-l+1)$-matrix
\[
M_{a_1 \cdots a_{l-1}} := \left(%
\begin{array}{lllll}
  x_0 \; \cdots \; x_{a_1-1} & \Big\arrowvert & x_{a_1+1} \; \cdots \; x_{a_2-1}
  & \Big\arrowvert \; \cdots \; \Big\arrowvert & x_{a_{l-1}+1} \; \cdots \; x_{n-1} \\
  x_1 \; \cdots \; x_{a_1}   & \Big\arrowvert & x_{a_1+2} \; \cdots \; x_{a_2}
  & \Big\arrowvert \; \cdots \; \Big\arrowvert & x_{a_{l-1}+2} \; \cdots \; x_n
\end{array}%
\right) .
\]
\\
C) Let $V \in GL_2(k)$ and $W \in GL_{n-l}(k).$ Then, the
$2\times2$-minors of the {\sl conjugate matrix } $VM_{a_1 \cdots
a_{l-1}}W^{-1}$ generate the same ideal as the $2\times2$-minors of the matrix
of $M_{a_1 \cdots a_{l-1}}.$ So, if we subject $M_{a_1 \cdots
a_{l-1}}$ to regular $k$-linear row and column transformations,
the $2\times 2$-minors of the resulting $2 \times (n-l+1)$-matrix
still generate the vanishing ideal of the same scroll $S_{a_1
\cdots a_{l-1}}.$ \\
D) A $2 \times (n-l+1)$-matrix $N$ whose entries are linear forms
in $k[x_0,\ldots,x_r]$ is said to be {\sl $1$-generic,} if no
conjugate of $N$ has a zero entry. Observe that the property of
being $1$-generic is preserved under conjugation. Moreover, if
$N'$ is obtained by deleting some columns from the $1$-generic
matrix $N,$ then $N'$ is again $1$-generic.

Finally, let $N$ be a $1$-generic $2 \times (n-l+1)$-matrix whose
entries are linear forms in $k[x_0,\ldots,x_r].$  Let
$y_0,\ldots,y_m$ be a basis of the $k$-vector space $L \subset
k[x_0,\ldots,x_r]_1$ generated by the entries of $N.$ Then, $m >
n-l+1,$ thus $h := m-n+l > 1.$

Moreover, there are integers $b_1, \ldots, b_{h-1} \in \mathbb N$
such that with  $b_0 = -1$ and $b_h = m$ we have $b_i - b_{i-1} >
1$ for $i = 1, \ldots, h,$ and such that $N$ is conjugate to the
$2\times (m-l+1)$-matrix
\[
N' :=
\left(%
\begin{array}{lllll}
  y_0 \; \cdots \; y_{b_1-1} & \Big\arrowvert & y_{b_1+1} \; \cdots \;
  y_{b_2-1}
  & \Big\arrowvert \; \cdots \; \Big\arrowvert & y_{b_{h-1}+1} \; \cdots \;
  y_{m-1} \\
  y_1 \; \cdots \; y_{b_1}   & \Big\arrowvert & y_{b_1+2} \; \cdots \;
  y_{b_2}
  & \Big\arrowvert \; \cdots \; \Big\arrowvert & y_{b_{h-1}+2} \; \cdots \;
  y_m
\end{array}%
\right) .
\]
So, by parts B) and C) the $2 \times 2$-minors of $N$ generate the
vanishing ideal of a rational normal $h$-fold scroll in $\mathbb
P^m_k = \Proj(k[y_0,\ldots,y_m]).$
\end{remark}

\begin{remark} \label{Remark 8.2} A) Let $s \in \mathbb N$ and let
$\tilde{X} \subset \mathbb P^s_k = \Proj(R), R = k[x_0, \ldots,
x_s],$ be a cone over a rational normal scroll $S_{a_1 \cdots
a_{l-1}} \subset \mathbb P^n_k$ with $n \in \{1,2,\ldots,s\}.$
According to the previous remark we may assume that the vanishing
ideal of $\tilde{X}$ in $R$ is generated by the $2 \times
2$-minors of the $2 \times (n-l+1)$-matrix $M_{a_1 \cdots
a_{l-1}}.$ In this case (and with the convention that $\mathbb
P^{-1}_k = \emptyset$ and $\dim \emptyset = -1$), the vertex
$\Sing (\tilde{X})$ of $\tilde{X}$ is given by $\mathbb
P^{s-n-1}_k = \Proj( R/(x_0, \ldots, x_n)R)$ and so $\dim \Sing(
\tilde{X})
= s-n-1$ and $\dim \tilde{X} = l + s-n.$ \\
B) According to \ref{Remark 8.1} C) the $2 \times 2$-minors of any
matrix obtained from $M_{a_1 \cdots a_{l-1}}$ by $k$-linear row
and column operations generate the vanishing ideal of
$\tilde{X}$ in $R.$ \\
C) Let $N$ be a $1$-generic $2 \times (n-l+1)$-matrix whose
entries are linear forms in $k[x_0,\ldots,x_n].$ Let $y_0, \ldots,
y_m$ be a basis of the $k$-space spanned by the entries of $N$ and
let $h := m-n+l.$ Then, by parts A) and B) and by Remark
\ref{Remark 8.1} C), the $2 \times 2$-minors of $N$ generate the
vanishing ideal of a cone $Y \subset \mathbb P^s_k$ over a
rational normal $h$-fold scroll $Z \subset \mathbb P^m_k.$
In particular $\dim Y = h+s-m = l+s-n$ and $\dim
\Sing(Y) = s-m-1.$
\end{remark}

We now prove the result which shall be crucial in the rest of this
chapter.

\begin{theorem} \label{Theorem 8.3} Let $\tilde{X} \subset \mathbb
P^{r+1}_k$ be a (cone over a) rational normal scroll and let $
\varrho : \mathbb P^{r+1}_k \setminus \{p\} \to \mathbb P^r_k$ be
a linear projection from a point $p \in \mathbb P^{r+1}_k
\setminus \tilde{X}.$ Then, there is a (cone over a) rational
normal scroll $Y \subset \mathbb P^r_k$ such that $Y \supset
\varrho(\tilde{X}), \dim Y = \dim \tilde{X} + 1$ and $ \dim \Sing
(\tilde{X}) \leq \dim \Sing (Y) \leq \dim \Sing(\tilde{X}) + 3.$
\end{theorem}

\begin{proof} According to Remark \ref{Remark 8.2} A) we may
assume that the vanishing ideal of $\tilde{X}$ in $S' =
k[x_0,\ldots,x_{r+1}]$ is generated by the $2 \times 2$-minors of
the $2 \times (n-l+1)$-matrix
\[
M := \left(%
\begin{array}{lllll}
  x_0 \; \cdots \; x_{a_1-1} & \Big\arrowvert & x_{a_1+1} \; \cdots \;
  x_{a_2-1} & \Big\arrowvert \; \cdots \; \Big\arrowvert & x_{a_{l-1}+1} \;
  \cdots \; x_{n-1} \\
  x_1 \; \cdots \; x_{a_1}   & \Big\arrowvert & x_{a_1+2} \; \cdots \;
  x_{a_2}   & \Big\arrowvert \; \cdots \; \Big\arrowvert & x_{a_{l-1}+2} \;
  \cdots \; x_n
\end{array}%
\right)
\]
with appropriate integers $n, l, a_1, \ldots, a_{l-1} \in \mathbb N$
such that,  with $a_0 =-1, a_l = n,$ we have $a_i - a_{i-1} > 1$ for
$i = 1,\ldots,l, n \leq r +1, \dim \tilde{X} = l+ r + 1 -n$ and
$\dim \Sing (\tilde{X}) = r-n.$

Let $ p := (c_0:c_1: \cdots:c_{r+1}).$ As $p \not\in \tilde{X}$
there are two different indices $i, j \in \{0,1,\ldots,n-1\}
\setminus \{a_1, a_2, \ldots,a_{l-1}\}$ such that
\[
\delta := \det \left(
\begin{array}{cc}
  c_i & c_j \\
  c_{i+1} & c_{j+1} \\
\end{array}
\right) \not= 0.
\]
Without loss of generality we may assume that $c_{i+1} \not= 0.$
Define
\[
y_{\alpha} :=
\begin{cases}
    x_{\alpha}, & \hbox{if}\quad  \alpha = i+1, \\
    x_{\alpha} - \frac{c_{\alpha}}{c_{i+1}} x_{i+1}, & \hbox{if}
    \quad
    \alpha \in \{0,\ldots,r+1\}\setminus \{i+1\}.
\end{cases}
\]
Then $S' = k[y_0,y_1,\ldots,y_{r+1}]$ and with respect to the
coordinates $y_0,\ldots,y_{r+1}$ the point $p$ may be written as
$(0:\cdots:0:1:0:\cdots:0)$ with the entry "1" in the $(i+1)$-th
position. Therefore we may assume that the projection $\varrho$ is
induced by the inclusion map
\[
S'' := k[y_0, \ldots, y_i, y_{i+2}, \ldots, y_{r+1}]
\hookrightarrow S'.
\]
We now express the indeterminates which occur in the matrix $M$ in
terms of the variables $y_{\alpha}:$
\[
x_{\alpha} =
\begin{cases}
    y_{\alpha}, & \hbox{for}\quad  \alpha = i+1, \\
    y_{\alpha} + \frac{c_{\alpha}}{c_{i+1}} y_{i+1},, & \hbox{if}
    \quad  \alpha
    \in \{0,\ldots,r+1\}\setminus \{i+1\}.
\end{cases}
\]
Let us first assume that $i < j.$ Then, the $2\times2$-submatrix
$U$ of $M$ which contains $x_i$ and $x_j$ in its first row takes
the form
\[
U = \left(
\begin{array}{cc}
  y_i + \frac{c_i}{c_{i+1}} y_{i+1} & y_j + \frac{c_j}{c_{i+1}}
  y_{i+1} \\
  y_{i+1}  & y_{j+1} + \frac{c_{j+1}}{c_{i+1}} y_{i+1}
\end{array}
\right) .
\]
Now performing sucessively $k$-linear row and column
operations we finally get the following transformed matrix
\[
U' = \left(
\begin{array}{cc}
  y_i &  y_j - \frac{c_i}{c_{i+1}} y_{j+1} - \frac{c_{j+1}}{c_{i+1}}
  y_i  -\frac{\delta}{c_{i+1}^2}y_{i+1}\\
  y_{i+1} & y_{j+1} \\
\end{array}%
\right).
\]
If $i+1 = j,$ then by performing $k$-linear row and column operations, $U$ can be
brought to the form
\[
U' = \left(
\begin{array}{cc}
  y_i &   - \frac{c_i}{c_{i+1}} y_{i+2} - \frac{c_{i+2}}{c_{i+1}}
  y_i -\frac{\delta}{c_{i+1}^2}y_{i+1}\\
  y_{i+1} & y_{i+2} \\
\end{array}%
\right) .
\]

Let $M'$ be the $2\times (n-l+1)$-matrix which is obtained if the
above row and column operations are performed with the whole
matrix $M.$ Observe that the submatrix $U$ of $M$ is transformed
into the submatrix $U'$ of $M',$ which sits in the same columns as
$U.$ Now, as $-\frac{\delta}{c_{i+1}^2} \not= 0$ we may add
appropriate $k$-multiples of the columns of $U'$ to the columns of
$M'$ to remove the indeterminate $y_{i+1}$ from all entries of
$M'$ which do not belong to the two columns of $U'.$ So, we get a
$2 \times (n-l+1)$-matrix $\tilde{M},$ conjugate to $M.$ In
particular, $\tilde{M}$ is $1$-generic (cf Remark \ref{Remark 8.1}
D)) and the entries of $\tilde{M}$ span the same $k$-space as the
entries of $M,$ namely $\sum_{t=0}^n kx_t = \sum_{t=0}^n ky_t.$
Now, let $N$ be the matrix of size $2 \times (n-l-1) = 2 \times
(n-1 -(l+1)+1)$ obtained by deleting the two columns of $U'$ from
$\tilde{M}.$ Then $N$ is $1$-generic (cf Remark \ref{Remark 8.1}
D)) and $y_{i+1}$ does not appear in $N.$ So, the entries of $N$
span a subspace
\[
L \subseteq \Sigma_{t=0, t \not= i+1}^n ky_t \subset
k[y_0,\ldots,y_i,y_{i+1},\ldots,y_n] \subset S'',
\]
whose dimension $m$ is such that $n-m \in \{1,2,3,4\}.$

By Remark \ref{Remark 8.2} C) the ideal $I_2(N) \subset S''$
generated by the $2 \times 2$-minors of $N$ is the vanishing ideal
of a cone $Y \subset \mathbb P^r_k = \Proj(S'')$ over a rational
normal scroll such that $\dim Y = (l+1)+r-(n-1) = \dim \tilde{X}
+1$ and $\dim \Sing(Y) = r-m-1 = \dim \Sing (\tilde{X}) +
(n-m-1).$ As $n-m-1 \in \{0,1,2,3\}$ we get $\dim \Sing
(\tilde{X}) \leq \dim \Sing (Y) \leq \dim \Sing(\tilde{X}) +3.$
According to Remark \ref{Remark 8.2} B) the ideal $I_2(\tilde{M})
\subset S'$ generated by the $2\times 2$-minors of $\tilde{M}$ is
the vanishing ideal of $\tilde{X}.$ Therefore $I_2(\tilde{M}) \cap
S''$ is the vanishing ideal of $\varrho(\tilde{X}).$ As each
$2\times 2$-minor of $N$ is a $2 \times 2$-minor of $\tilde{M},$
it follows $I_2(N) \subseteq I_2(\tilde{M}) \cap S''$ and hence $Y
\supseteq \varrho(\tilde{X}).$

This settles the case $i < j.$ If $j < i$ we first commute the
columns of $U$ and then conclude as above.
\end{proof}

We now apply the previous result to varieties of almost minimal
degree. We still keep the convention that $\mathbb P^{-1}_k =
\emptyset$ and $\dim \emptyset = -1,$ and begin with a preliminary
remark.

\begin{remark} \label{Remark 8.3a} (cf \cite{H}) A) Let $l, n , a_1,
\ldots, a_{l-1} \in \mathbb N$ be as in Remark \ref{Remark 8.1}
and consider the rational $l$-fold scroll $S_{a_1\cdots a_{l-1}}.$
We set $a_0 = -1, a_l = n, d_i = a_i - a_{i-1} -1,$ for $i = 1, \ldots,
l$ and define the linear subspaces
\[
\mathbb E_i = \Proj(S/\Sigma_{j \not\in \{a_{i-1}+1,\ldots,a_i\}}
Sx_j) = \mathbb P^{d_i}_k \subset \mathbb P^n_k, i = 1,\ldots, l.
\]
For each $i \in \{1,\ldots,l\}$ we consider the Veronese
embedding
\[
\nu_i : \mathbb P^1_k \to \mathbb E_i, \; (s:t) \mapsto (s^{d_i} :
s^{d_i-1}t : \ldots : st^{d_i -1} : t^{d_i}),
\]
so that $\nu_i(\mathbb P^1_k) \subset \mathbb E_i = \mathbb
P^{d_i}_k$ is a rational normal curve. Now, for each $q \in
\mathbb P^1_k$ let
\[
\mathbb E(q) = \langle \nu_1(q), \ldots, \nu_l(q) \rangle =
\mathbb P^{l-1}_k \subset \mathbb P^n_k
\]
be the projective space spanned by the $l$ points $\nu_1(q),
\ldots, \nu_l(q).$ Then $q \not= q'$ implies $\mathbb
E(q) \cap \mathbb E(q') = \emptyset$ for all $q, q' \in \mathbb
P^1_k.$ Moreover, $S_{a_1\cdots a_{l-1}} = \cup_{q \in \mathbb
P^1_k} \mathbb E(q).$ \\
B) Keep the above notation. For each $1 \leq i \leq l$ let
$\mathbb V_i = k^{d_i+1} \subset k^{n+1}$ be the affine cone over
$\mathbb E_i,$ fix $q = (s : t) \in \mathbb P^1_k$ and set $v_i =
(s^{d_i}, s^{d_i-1}t, \ldots, t^{d_i}).$ Moreover, for each $1
\leq i \leq l$ let $w_i = (w_{i0},\ldots,w_{id_i}) \in \mathbb
V_i$ such that $(w_{i0} : \ldots :w_{id_i}) \in \mathbb E_i
\setminus \nu_i(q)$ is a point on the tangent of the rational
normal curve $\nu_i(\mathbb P^1_k) \subset \mathbb E_i$ in the
point $\nu_i(q).$ In particular, $v_i$ and $w_i$ are lineraly
independent. Now let $\pi : k^{n+1}\setminus \{0\}
\twoheadrightarrow \mathbb P^n_k = \mathbb P(k^{n+1})$ be the
canonical projection and let
\[
r = \pi(\Sigma_{i=1}^l r_iv_i) \in \mathbb E(q)= \pi(\oplus_{i=1}^l
kv_i \setminus \{0\}) = \mathbb P(\oplus_{i=1}^l kv_i),
\]
where $(r_1,\ldots, r_l) \in k^l \setminus \{0\}.$ Then the
tangent space to the scroll $S_{a_1\cdots a_{l-1}}$ in the point
$r$ is given by
\[
T_r(S_{a_1\cdots a_{l-1}}) = \mathbb P(k(\Sigma_{i=1}^lr_iw_i) +
\oplus_{i=1}^lkv_i) = \langle \mathbb E(q) \cup \pi(\Sigma_{i=1}^l
r_iw_i) \rangle.
\]
From this we easily deduce that $T_r(S_{a_1\cdots a_{l-1}}) \cap
T_{r'}(S_{a_1\cdots a_{l-1}}) = \mathbb E(q)$ for all $r, r' \in
\mathbb E(q)$ with $r \not= r'.$ \\
C) Now, let $s \in \mathbb N$ such that $n < s$ and let $\tilde{X}
\subset \mathbb P^s_k = \Proj(R), R = k[x_0,\ldots,x_n],$ be a cone
over the rational normal $l$-fold scroll
\[
S_{a_1\cdots a_{l-1}} \subset \mathbb P^n_k =
\Proj(R/(x_{n+1},\ldots,x_s)R = k[x_0,\ldots,x_n]).
\]
Then, the vertex of $\tilde{X}$ is given by $\Sing(\tilde{X}) =
\Proj(R/(x_0,\ldots,x_n)) = \mathbb P^{s-n-1}_k \subset \mathbb
P^s_k$ (cf Remark \ref{Remark 8.2} A)). Now, for each $q \in
\mathbb P^1_k$ let
\[
\mathbb F(q) = \langle \mathbb E(q) \cup \Sing(\tilde{X}) \rangle
= \mathbb P^{l+s-n-1}_k = \mathbb P^{\dim \tilde{X}}_k \subset \mathbb P^s_k
\]
be the linear subspace spanned by $\mathbb E(q) = \mathbb
P^{l-1}_k \subset \mathbb P^s_k$ and the vertex $\Sing(\tilde{X})$
of $\tilde{X}.$ Then by part A), $q \not= q'$ implies
that $\mathbb F(q) \cap \mathbb F(q') = \Sing(\tilde{X})$ for all
$q, q' \in \mathbb P^1_k$ and moreover $\tilde{X} = \cup_{q \in
\mathbb P^1_k} \mathbb F(q).$

It also follows easily from part B), that for any $q \in \mathbb
P^1_k$ and any $r \in \mathbb E(q) \setminus \Sing(\tilde{X})$ the
tangent space of $\tilde{X}$ at $r$ is given by $T_r(\tilde{X}) = 
\langle T_{\tilde{r}}(S_{a_1\cdots a_{l-1}})
\cup \Sing(\tilde{X}) \rangle = \mathbb P^d_k,$ where $\tilde{r}
=(r_0: \ldots : r_n)$ is the canonical projection of $r$ from
$\Sing(\tilde{X}).$ As a consequence of the last statement in part
B) we thus get for all $q \in \mathbb P^1_k$ and all $r, r' \in
\mathbb F(q) \setminus \Sing(\tilde{X}):$ If $r \not= r',$ then
$T_r(\tilde{X}) \cap T_{r'}(\tilde{X}) = \mathbb F(q).$
\end{remark}

\begin{theorem} \label{Corollary 8.4} Let $X \subset \mathbb
P^r_k$ be a variety of almost minimal degree which is the
projection of a (cone over a) rational normal scroll $\tilde{X}
\subset \mathbb P^{r+1}_k$ with $\dim \Sing(\tilde{X}) =:h$ from a
point $p \in \mathbb P^{r+1}_k \setminus \tilde{X}.$ Then
\begin{itemize}
\item[(a)] $X$ is contained in a (cone over a) rational normal
scroll $Y \subset \mathbb P^r_k$ such that $\codim_Y (X) = 1$ and
$h \leq \dim \Sing(Y) \leq h+3.$

\item[(b)] $X$ is of arithmetic depth $t \leq h+5.$
\end{itemize}
\end{theorem}

\begin{proof} (a): This is clear by Theorem \ref{Theorem 8.3}.\\
(b): Assume first that $X$ is not arithmetically Cohen-Macaulay.
Then, the non CM-locus $Z$ of $X$ is a linear subspace $\mathbb
P^{t-2}_k$ of $\mathbb P^r_k$ (cf Theorem \ref{5.6 Theorem} (d),
(f)). As $\codim_Y (X) = 1, {\mathcal O}_{X,x}$ is a
Cohen-Macaulay ring for each point $x \in X\setminus \Sing (Y).$
It follows that $\mathbb P^{t-2}_k = Z \subseteq X \cap \Sing(Y)$
and hence $t-2 \leq \dim \Sing(Y) \leq h+3.$

Now, let $X$ be arithmetically Cohen-Macaulay. In this case we
conclude by a geometric argument which in fact also implies in the
previous case. Let $d = \dim X.$ After an appropriate change of
coordinates, we may assume that we are in the situation of Remark
\ref{Remark 8.3a} B) and C). So, we may write $\tilde{X} = \cup_{q
\in \mathbb P^1_k} \mathbb F(q),$ where $\mathbb F(q) = \mathbb
P^{d-1}_k \subset \mathbb P^r_k$ is a linear subspace for all $q
\in \mathbb P^1_k.$ Let $U = \{ (q,q') \in \mathbb P^1_k \times
\mathbb P^1_k | q \not= q'\}.$ Then, according to Remark
\ref{Remark 8.3a} C) we have $\mathbb F(q) \cap \mathbb F(q') =
\Sing(\tilde{X}) = \mathbb P^h_k,$ whenever $(q,q') \in U.$ Now, for each pair $(q,q') \in
U$ consider the linear subspace
\[
\mathbb H(q,q') = \langle \mathbb F(q) \cup \{p\} \rangle \cap
\langle \mathbb F(q') \cup \{p\} \rangle \subset \mathbb
P^{r+1}_k.
\]
Observe that $\Sing(\tilde{X}) \subseteq \mathbb H(q,q')$ and
$\dim \mathbb H(q,q') \leq h+2$ for all $(q,q') \in U.$ Moreover
$\dim \mathbb H(q,q') = h+2$ if and only if $p \in \langle \mathbb
F(q) \cup \mathbb F(q') \rangle.$ Consequently we have $\dim
\mathbb H(q,q') = h+ 2$ if and only if there is a line running
through $p$ and intersecting $\mathbb F(q)$ and $\mathbb F(q').$
Clearly such a line is contained in $\mathbb H(q,q')$ and its
intersection points with $\mathbb F(q)$ and $\mathbb F(q')$ are
different as $p \not\in \mathbb F(q) \cup \mathbb F(q').$

Let $V \subseteq U$ be the closed subset of all pairs $(q,q')$ for which
$\dim \mathbb H(q,q') = h+2.$ It follows that the union of all
proper secant lines of $\tilde{X}$ which run through $p$ is a
subset of $W = \cup_{(q,q')\in V} \mathbb H(q,q').$ Moreover, it
follows from the last statement of Remark \ref{Remark 8.3a} C),
that for each point $q \in \mathbb P^1_k$ there is at most one
point $r(q) \in \mathbb F(q) \setminus \Sing(\tilde{X})$ such that
there is a tangent line $l(q) = \mathbb P^1_k$ of $\tilde{X}$ at
$r(q)$ running through $p.$ Let $T \subset \mathbb P^1_k$ be the
closed subset of all $q \in \mathbb P^1_k$ for which this happens.
Then, all tangents to non-singular points of $\tilde{X}$ through
$p$ are contained in $Y = \cup_{q \in T} l(q).$ Finally observe
that the remaining tangents are the lines running trough $p$ and
$\Sing(\tilde{X}).$ It follows $\mathbb P^{t-1}_k =
\Sec_p(\tilde{X}) \subseteq W \cup Y \cup \langle \Sing(\tilde{X}) \cup \{p\}
\rangle$ and hence $t-1 \leq \max \{ \dim W, \dim Y, h+1\}.$ As
$\{\mathbb H(q,q') | (q,q') \in V\}$ is a family of linear
$(h+2)$-subspaces of $\mathbb P^{r+1}_k,$ it follows $\dim W \leq
h + 2 + \dim V \leq h+4.$ As $\{ l(q) | q \in T\}$ is a family of
lines we have $\dim Y \leq 1 + \dim T \leq 2.$ So we get $t-1 \leq
h+4,$ hence $t \leq h+5.$
\end{proof}

\begin{corollary} \label{Corollary 8.5} Let $X \subset \mathbb
P^r_k$ be a variety of almost minimal degree which is 
a projection of a rational normal scroll $\tilde{X} \subset
\mathbb P^{r+1}_k$ from a point $p \in \mathbb P^{r+1}_k \setminus
\tilde{X}.$ Then $X$ is  of arithmetic depth $t \leq 4.$
\end{corollary}

\begin{proof} Clear from Theorem \ref{Corollary 8.4}.
\end{proof}

As a final comment of this section let us say something about the
exceptional case of projections of the Veronese surface.

\begin{remark} \label{Remark 8.5} (The exceptional case) Let us
recall that the Veronese surface $F \subset \mathbb P^5_k$ is
defined by the $2\times 2$-minors of the matrix
\[
M =
\left(%
\begin{array}{ccc}
  x_0 & x_1 & x_2 \\
  x_1 & x_3 & x_4 \\
  x_2 & x_4 & x_5 \\
\end{array}%
\right).
\]
Let $p \in \mathbb P^5_K \setminus F$ denote a closed point.
Suppose that $\rank M \hspace{-.2cm}\mid_p = 3,$ i.e. the case of
a generic point and remember that $\det M = 0$ defines the secant
variety of $F.$ Then the projection of $F$ from $p$ defines a
surface $X \subset \mathbb P^4_K$ of almost minimal degree and
$\depth A = 1.$

Recall that $\dim_k (I)_2 = 0$ (cf. Corollary \ref{4.4 Corollary}
C)). Therefore the surface $X$ is cut out by  cubics, i.e. it is
not contained in a variety of minimal degree.
\end{remark}

\section{Betti numbers}
\label{8. Betti numbers}

Our next aim is to study the Betti numbers of $A$ if $X$ is of
almost minimal degree, non-arithmetically Cohen-Macaulay and a
projection of a (cone over a) rational normal scroll.

\begin{lemma} \label{9.1 Lemma} Assume that $X \subset \mathbb
P^r_k$ is of almost minimal degree and of arithmetic depth $\leq d
= \dim X.$ Let $B = \Hom_A(K(A), K(A)).$ Then
\[
\Tor_i^S(k,B) \simeq \left\{%
\begin{array}{ll}
    k(0) \oplus k(-1), & \hbox{if }  \quad i = 0,\\
    k^{b_i}(-i-1), & \hbox{if } \quad 0 < i \leq r-d,\\
    0, & \hbox{if }\quad  r-d < i,\\
\end{array}%
\right.
\]
where $b_i = (r+1-d)\binom{r-d}{i} - \binom{r-d}{i+1}$ for $1 \leq
i \leq r-d.$
\end{lemma}

\begin{proof} By Theorem \ref{5.3 Theorem} (a) the $A$-module $B$
is Cohen-Macaulay and hence of depth $d+1$ over $S.$ Therefore
$\Tor_i^S(k,B) = 0$ for all $i > r-d.$ According to Theorem
\ref{5.3 Theorem} (b) there is a short exact sequence of graded
$S$-modules
\begin{equation} \label{9.1}
0 \to A \to B \to C \to 0, \quad C \simeq
(S/(y_{t-2},\ldots,y_r)S)(-1),
\end{equation}
where $y_0,\ldots,y_r$ form a generic set of linear forms of $S.$
In particular, $C$ is of dimension $t-1 < d$ and generated by a
single element of degree 1. This already shows that $\Tor_0^S(k,B)
\simeq k(0)\oplus k(-1).$

Applying cohomology to the above short exact sequence we get an
isomorphism $H^{d+1}(A) \simeq H^{d+1}(B)$ which shows that
$\text{end } H^{d+1}(B) = -d$ (cf Theorem \ref{4.2 Theorem} (b)).
As $\depth B = d+1$ it follows $\reg B = 1.$ Moreover the above
exact sequence yields
\[
\dim_k B_1 = \dim_k A_1 +1 = r+2 = \dim_k (S(0)\oplus S(-1))_1.
\]
So, the graded $S$-module $B$ must have a minimal free resolution
of the form
\[
0 \to S^{b_{r-d}}(-r+d-1) \to \ldots \to S^{b_i}(-i-1) \to \ldots
\to S^{b_1}(-2) \to S\oplus S(-1) \to B \to 0
\]
with $b_1,\ldots,b_{r-d} \in \mathbb N.$

As $B$ is a Cohen-Macaulay module of dimension $d+1,$ regularity 1
and of multiplicity $\deg X = r-d+2$ (cf Theorem \ref{5.3 Theorem}
(c)) its Hilbert series is given by
\[
F(\lambda, B) = \frac{1+ (r+1-d)\lambda}{(1-\lambda)^{d+1}} .
\]
On the use of Betti numbers $b_i$ we also may write
\[
F(\lambda, B) = \frac{1}{(1-\lambda)^{d+1}} (1 + \lambda +
\sum_{i=0}^{r-d} (-1)^i b_i \lambda^{i+1}).
\]
Comparing coefficients we obtain
\[
b_i = (r+1-d)\binom{r-d}{i} -\binom{r-d}{i+1}, \; i = 1, \ldots,
r-d,
\]
as required.
\end{proof}

Next we recall a well known result about the Betti numbers of a
variety of minimal degree.

\begin{lemma} \label{9.2 Lemma} Let $Y \subset \mathbb P^r_k$ be a
variety of minimal degree with $\dim Y = d+1.$ Let $U$ be the
homogeneous coordinate ring of $Y.$ Then
\[
\Tor_i^S(k,U) \simeq \left\{%
\begin{array}{ll}
    k, & \hbox{if } \quad i = 0,\\
    k^{c_i}(-i-1), & \hbox{if } \quad 0 < i < r-d,\\
    0, & \hbox{if }\quad  r-d \leq i,\\
\end{array}%
\right.
\]
where $c_i = i \binom{r-d}{i+1}$ for $1 \leq i < r-d.$
\end{lemma}

\begin{proof} This is well known (cf for instance \cite{E}). In
fact the Eagon-Northcott complex provides a minimal free
resolution of $U$ over $S.$
\end{proof}

\begin{theorem} \label{9.3 Theorem} Let $X \subset \mathbb P^r_k$
be a variety of almost minimal degree which is the projection of a
(cone over a) rational normal scroll $\tilde{X} \subset \mathbb
P^{r+1}_k$ from a point $p \in \mathbb P^{r+1}_k \setminus
\tilde{X}.$ Assume that $t := \depth A \leq d := \dim X.$  Then
\[
\Tor_i^S(k,A) \simeq
\begin{cases}
k, & \text{if} \quad i = 0,\\
k^{u_i}(-i-1) \oplus k^{v_i}(-i-2), & \text{if} \quad 0 < i \leq
r-t+1,\\
0, & \text{if} \quad r-t+1 < i,
\end{cases}
\]
where
\begin{itemize}
\item[(a)]
\begin{itemize}
\item[] $u_1 = t + \binom{r+1-d}{2} -d-2,$

\item[] $i \binom{r-d}{i+1} \leq u_i \leq (r+1-d)\binom{r-d}{i} -
\binom{r-d}{i+1},\quad \text{if } \; 1 < i < r-2d+t-1,$

\item[] $u_i = i \binom{r-d}{i+1},\quad \text{if } \; r-2d+t-1 \leq i
< r-d,$

\item[] $u_i = 0, \quad \text{if } \; r-d\leq i < r-t+1,$

\medskip
\end{itemize}
\item[(b)]
\begin{itemize}
\item[] $\max \{ 0, \binom{r-t+2}{i+1} - (i+2)\binom{r-d}{i+1}\}
\leq v_i \leq \binom{r-t+2}{i+1}, \quad \text{if } \; 1 \leq i <
r-2d+t-2,$

\item[] $v_i = \binom{r-t+2}{i+1} - (i+2)\binom{r-d}{i+1},\quad
\text{if } \; r-2d+t-2 \leq i < r-d,$

\item[] $v_i = \binom{r-t+2}{i+1},\quad \text{if } \; r-d \leq i \leq
r-t+1.$
\end{itemize}
\end{itemize}
Moreover, $v_i - u_{i+1} = \binom{r-t+2}{i+1} -
(r-d+1)\binom{r-d}{i+1} + \binom{r-d}{i+2}$ for all $1 \leq i <
r-d.$
\end{theorem}

\begin{proof} As $\depth A = t$ and $\reg A = 2$ (cf Theorem \ref{4.2
Theorem}) the modules $\Tor_i^S(k,A)$ behave as stated in the main
equality. So let the numbers $u_i, v_i$ be defined according to
this main equality with the convention that $u_i = v_i = 0$ for $i
> r-t+1.$ Moreover let $b_i$ and $c_i$ be as in Lemma \ref{9.1
Lemma} resp. \ref{9.2 Lemma} with the convention that $b_i = 0$
for $i > r-d$ and $c_i = 0$ for $i\geq r-d.$

According to Corollary \ref{Corollary 8.5} there is a (cone over
a) rational normal scroll $Y \subset \mathbb P^r_k$ of dimension
$d+1$ such that $X \subset Y.$ Let $J \subset S$ be the vanishing
ideal of $Y$ and let $U := S/J$ be the homogeneous coordinate ring
of $Y.$ The short exact sequence $0 \to I/J \to U \to A \to 0$
together with Lemma \ref{9.2 Lemma} implies  short exact sequences
\begin{equation} \label{9.2}
0 \to k^{c_i}(-i-1) \to k^{u_i}(-i-1) \oplus k^{v_i}(-i-2) \to
\Tor_{i-1}^S(k,I/J) \to 0
\end{equation}
for all $i \geq 1.$

Keep in mind that $\beg (I/J) = 2, u_1 = \dim_k I_2 =
t+\binom{r+1-d}{2} -d-2$ (cf Corollary \ref{4.4 Corollary} (c))
and $\dim_k J_2 = \binom{r-d}{2}$ (cf Lemma \ref{9.2 Lemma}) so
that
\[
\dim_k (I/J)_2 = r-2d+t-2.
\]
Whence, by Green's Linear Syzygy Theorem (cf \cite[Theorem
7.1]{ES}) we have
\[
\Tor_j^S(k, I/J)_{j+2} = 0\; \text{ for all } j \geq r-2d+t-2.
\]
So, the sequence \ref{9.2} yields that $u_i = c_i$ for all $i \geq
r-2d+t-1.$ This proves statement (a) in the range $i \geq
r-2d+t-1.$ The sequence \ref{9.2} also yields that $c_i \leq u_i$
for all $i \leq r-2d+t-1.$

Next we consider the short exact sequence of graded $S$-modules
\ref{9.1}. In particular we have $\Tor_i^S(k, C) \simeq
k^{\binom{r-t+2}{i}}(-i-1)$ for all $i \in \mathbb N_0.$ So, by
the sequence \ref{9.1} and in view of Lemma \ref{9.2 Lemma} we get
exact sequences
\begin{multline}
\label{9.3} k^{b_{i+1}}(-i-2)  \to k^{\binom{r-t+2}{i+1}}(-i-2)
\to k^{u_i}(-i-1) \oplus k^{v_i}(-i-2)\\ \to k^{b_i}(-i-1)
 \to k^{\binom{r-t+2}{i}}(-i-1) \to k^{u_{i-1}}(-i) \oplus
k^{v_{i-1}}(-i-1) \to k^{b_{i-1}}(-i)
\end{multline}
for all $ i \geq 2.$ Now, we read off that $u_i \leq b_i$ for all
$i \geq 1$ and statement (a) is proved completely.

The sequence \ref{9.3} also yields that
\begin{equation} \label{9.4}
v_i = u_{i+1} - b_{i+1} + \binom{r-t+2}{i+1} \; \text{for all }
i\geq 1.
\end{equation}
Observe that $c_{i+1}-b_{i+1} = -(i+2)\binom{r-d}{i+1}$ for $1
\leq i < r-d.$ If $ 1 \leq i < r-2d+t-1,$ statement (a) gives
$c_{i+1} \leq u_{i+1} \leq b_{i+1}$ so that
\[
c_{i+1} - b_{i+1} + \binom{r-t+2}{i+1} \leq v_i \leq
\binom{r-t+2}{i+1}.
\]
This proves the first estimate of statement (b).

If $r-2d+t-2 \leq i < r-d,$ statement (a) yields $u_{i+1} =
c_{i+1},$ hence $v_i = c_{i+1}-b_{i+1} + \binom{r-t+2}{i+1}. $
This proves the second claim of statement (b). Finally, if $r-d
\leq i < r-t+1$ statement (a) and Lemma \ref{9.2 Lemma} yield that
$u_{i+1} = c_{i+1} = 0.$ Now the last claim of the statement (b)
follows by \ref{9.4}.
\end{proof}

\section{Examples}
\label{6. Examples}

In this final section we present a few examples which illustrate
the previous results. All calculations of the ``graded Betti
numbers'' $u_i$ and $v_i$ (cf Theorem~\ref{9.3 Theorem}) have been
performed by means of the computer algebra system {\sc Singular}
\cite{GrP}. As for rational scrolls and their secant varieties we
refer to \cite{C} and \cite{H}.

First, we present three examples of $3$-folds $X$ of almost
minimal degree in ${\mathbb P}^{11}_k$, one of them being defined
by $32$ quadrics, the second by $32$ quadrics and $1$ cubic, the
third by $32$ quadrics and $3$ cubics. These examples show that,
contrary to the number of defining quadrics (cf Corollary~\ref{4.4
Corollary} (c) ), the number of defining cubics may vary if the
embedding dimension $r$, the dimension $d$ and the arithmetic
$\depth t$ of $X$ are fixed. Notice that each smooth variety $X
\subset {\mathbb P}^{11}_k$ of almost minimal degree which is not
arithmetically Cohen-Macaulay is obtained by projecting a rational
scroll $\tilde {X} \subset {\mathbb P}^{12}_k$ from a point $p \in
{\mathbb P}^{12}_k \setminus \Sec(\tilde {X})$ (cf
Theorem~\ref{5.6 Theorem}).

We first fix some notation. Let $l, n, d_1, \ldots, d_l \in
\mathbb N$ such that $d_1 \leq d_2 leq \ldots d_l$ and
$\sum_{i=1}^l d_i = n-l+1.$ Let $a_i = i-1+ \sum_{j=1}^id_j, i =
1,\ldots, l-1.$ Then, we write $S(d_1,\ldots,d_l)$ for the
rational normal scroll $S_{a_1\cdots, a_{l-1}}$ (cf Remark
\ref{Remark 8.1}).

\begin{example}
\label{6.1 Examples} A) Let $\tilde {X} \subset {\mathbb
P}^{12}_k$ the $3$-scroll $S(2, 2, 6)$, thus the smooth variety of
degree $10$ defined by the $2 \times 2$ minors of the matrix
\[
\begin{pmatrix} x_0 & x_1 \Big\arrowvert x_3 & x_4
                       \Big\arrowvert x_6 & x_7 & x_8 & x_9 & x_{10}
                       & x_{11} \\
                       x_1 & x_2 \Big\arrowvert x_4 & x_5
                       \Big\arrowvert x_7 & x_8 & x_9 & x_{10} & x_{11}
                       & x_{12}
\end{pmatrix} .
\]
Its homogeneous coordinate ring is
\[
B = k[ (s,t)^2 u^5, (s,t)^2 v^5 , (s,t)^6 w ]
       \subset k[s, t, u, v, w] .
\]
Projecting $\tilde {X}$ from the point
\[
p_1 = (0 : 0 : 0 : 0 : 0 : 0 : 0 : 0 : 0 : 1 : 0 : 0 : 0) \in
       {\mathbb P}^{12}_k \setminus \tilde {X}
\]
we get a non-degenerate variety $X \subset {\mathbb P}^{11}_k$ of
dimension $3$ and of degree $\leq 10$, (cf Remark~\ref{3.3 Remark}
A) ). Let $S$ denote a polynomial ring in $12$ indeterminates, let
$I \subset S$ be the homogeneous vanishing ideal and let $A = S /
I$ be the homogeneous coordinate ring of $X$. Also, consider the
only not necessarily vanishing graded Betti numbers
\[
u_i := \dim _k \Tor^S_i(k,A)_{i + 1 }, \; v_i := \dim _k
\Tor^S_i(k,A)_{i + 2 }
\]
of $X$. These numbers present themselves as shown below:
\[
\begin{tabular}{c | r r r r r r r r r r r }
    $i$ & 1 & 2 & 3 & 4 & 5 & 6 & 7 & 8 & 9 & 10 & 11 \\ \hline
    $u_i $ & 32 & 130 & 234 & 234 & 140 & 48 & 7 & 0 & 0
      & 0 & 0 \\
    $v_i $ & 0 & 20 & 155 & 456 & 728 & 728 & 486 & 220 & 66
      & 12 & 1 \\ \hline
\end{tabular}
\]

In particular $t:= \depth A = 1$ so that $X$ cannot be of minimal
degree and hence $\deg X = 10 = 11 - 3 + 2$. Therefore, $X$ is of
almost minimal degree and of arithmetic $\depth 1$. In particular
the projection map $\nu : \tilde {X} \to X$ is an isomorphism (cf
Theorem~\ref{5.7 Corollary}) and so $X$ becomes smooth. Observe,
that $I$ is generated by $32$ quadrics.\\
B) Let $\tilde {X} \subset {\mathbb P}^{12}_k$ be as in part A)
but project $\tilde {X}$ from the point
\[
p_2 = (0 : 0 : 0 : 0 : 0 : 0 : 0 : 0 : 0 : 0 : 1 : 0 : 0) \in
       {\mathbb P}^{12}_k \setminus \tilde {X} .
\]
Again let $X \subset {\mathbb P}^{11}_k$ be the image of $\tilde
{X}$ under this projection and define $S, I, A$ as in part A). Now
the Betti numbers $u_i, v_i$ present themselves as follows:
\[
\begin{tabular}{c | r r r r r r r r r r r }
    $i$ & 1 & 2 & 3 & 4 & 5 & 6 & 7 & 8 & 9 & 10 & 11 \\ \hline
    $u_i $ & 32 & 131 & 234 & 234 & 140 & 48 & 7 & 0 & 0
      & 0 & 0 \\
    $v_i $ & 1 & 20 & 155 & 456 & 728 & 728 & 486 & 220 & 66
      & 12 & 1 \\ \hline
\end{tabular}
\]
So, as above, we see that $X$ is a smooth variety of almost
minimal degree having dimension $3$ and arithmetic $\depth 1$.
Observe, that now $I$ is minimally generated by $32$ quadrics and
$1$ cubic. So, if the same scroll $\tilde {X} = S(2, 2, 6) \subset
{\mathbb P}^{12}_k$ is projected from two different points $p_1,
p_2 \in {\mathbb P} ^{12}_k \setminus \Sec(\tilde {X})$, the
homological nature of the projection $X \subset {\mathbb
P}^{11}_k$ may differ.\\
C) Now, consider the scroll $\tilde {X}:= S(2, 4, 4) \subset
{\mathbb P}^{12}_k$, so that $\tilde {X}$ is the smooth variety of
dimension $3$ and degree $10$ defined by the $2 \times 2$-minors
of the matrix
\[
\begin{pmatrix} x_0 & x_1 \Big\arrowvert x_3 & x_4 & x_5
                       & x_6 \Big\arrowvert x_8 & x_9 & x_{10}
                       & x_{11} \\
                       x_1 & x_2 \Big\arrowvert x_4 & x_5 & x_6 & x_7
                       \Big\arrowvert x_9 & x_{10} & x_{11}
                       & x_{12}
\end{pmatrix} .
\]
Its homogeneous coordinate ring is
\[
B = k [ (s,t)^2 u^3, (s,t)^4 v, (s,t)^4 w ]
       \subset k[s, t, u, v, w] .
\]
Define $X \subset {\mathbb P}^{11}_k$ as the projection of $\tilde
{X}$ from the point $p_2 \in {\mathbb P}^{12}_k \setminus \tilde
{X}$ (cf part B) ). In this case, the Betti numbers $u_i$ and
$v_i$ of $X$ take the values listed in the following table:
\[
\begin{tabular}{c | r r r r r r r r r r r }
    $i$ & 1 & 2 & 3 & 4 & 5 & 6 & 7 & 8 & 9 & 10 & 11 \\ \hline
    $u_i $ & 32 & 133 & 248 & 234 & 140 & 48 & 7 & 0 & 0
      & 0 & 0 \\
    $v_i $ & 3 & 34 & 155 & 456 & 728 & 728 & 486 & 220 &
      66 & 12 & 1 \\ \hline
\end{tabular}
\]
So, again $X \subset {\mathbb P}^{11}_k$ is a smooth variety of
almost minimal degree having dimension $3$ and arithmetic $\depth
1$. But this time, besides $32$ quadrics three cubics are needed
to generate the homogeneous vanishing ideal $I$ of $X$.
\end{example}

The previous example where all of arithmetic $\depth 1$ and of
dimension $3$. By projecting rational $3$-scrolls in ${\mathbb
P}^{12}_k$ from appropriate points we also may obtain
$3$-dimensional varieties $X \subset {\mathbb P}^{11}_k$ of almost
minimal degree and of arithmetic $\depth $ not equal to 1. We
present two examples to illustrate this.

\begin{example}
\label{6.2 Examples} A) Next consider the $3$-scroll $\tilde {X}:=
S(3, 3, 4) \subset {\mathbb P}^{12}_k$ defined by the $2 \times
2$-minors of the matrix
\[
\begin{pmatrix} x_0 & x_1 & x_2 \Big\arrowvert  x_4 & x_5
                       & x_6 \Big\arrowvert x_8 & x_9 & x_{10}
                       & x_{11} \\
                       x_1 & x_2 & x_3 \Big\arrowvert x_5 & x_6 & x_7
                       \Big\arrowvert x_9 & x_{10} & x_{11}
                       & x_{12}
\end{pmatrix} .
\]
$\tilde {X}$ has the homogeneous coordinate ring
\[
B = k[ (s,t)^3 u^2, (s,t)^3 v^2 , (s,t)^4 w]
       \subset k[s, t, u, v, w] .
\]
We project $\tilde{X}$ from the point
\[
p_3 = (0 : 0 : 0 : 0 : 0 : 0 : 1 : 0 : 0 : 0 : 0 : 0 : 0) \in
       {\mathbb P}^{12}_k \setminus \tilde {X} .
\]
Like above we get a non-degenerate variety $X \subset {\mathbb
P}^{11}_k$ of degree $10  = \codim X + 2$ and Betti
numbers:
\[
\begin{tabular}{c | r r r r r r r r r r  }
    $i$ & 1 & 2 & 3 & 4 & 5 & 6 & 7 & 8 & 9 & 10  \\ \hline
    $u_i $ & 33 & 142 & 278 & 284 & 155 & 48 & 7 & 0 & 0
      & 0  \\
    $v_i $ & 1 & 9 & 40 & 141 & 266 & 266 & 156 & 55 &
      11 & 1  \\ \hline
\end{tabular}
\]
So, $X$ is of arithmetic $\depth 2.$

The tangent line of the curve
\[
\sigma : k \to \tilde {X} ; \ s \mapsto \sigma (s):=
       (0 : 0 : 0 : 0 : s^3 : s^2 : s : 1 : 0 : 0 : 0 : 0 : 0)
\]
in the point $\sigma (0)$ contains $p_3$. So, the secant cone
$\Sec_{p_3}(\tilde {X})$ -- which must be a line according to
Theorem~\ref{5.6 Theorem} (d) -- is just the line $\ell $ which
joins $p_3$ and $\sigma (0)$. The projection of $\ell $ from $p_3$
to ${\mathbb P}^{11}_k$ is the point
\[
q:= (0 : 0 : 0 : 0 : 0 : 0 : 1 : 0 : 0 : 0 : 0 : 0) \in X .
\]
So, in the notation of Theorem~\ref{5.6 Theorem}, we have
$\Sing(\nu ) = \varrho (\ell \setminus \{ p_2 \} ) = \{ q \}$.
Now, let $a:= \frac {u^2}{v^2} , b:= \frac {s}{t}, c:= \frac {tw}
{v^2}$. An easy calculation shows that there is an isomorphism
\[
\epsilon : X_{t^3 v^2} \overset {\simeq }{\longrightarrow }
       Y:= \Spec ( k[a, ab, b^2, b^3, b^2c, bc, c])
\]
with $\epsilon (q) = \underline {0}$, where $X_{t^3 v^2} \subset
X$ is the affine open neighborhood of $q$ defined by $t^3 v^2
\not= 0$. It is easy to verify, that ${\mathcal O}_{Y, \underline
{0}}$ is a $G$-ring and hence that $q \in X$ is a $G$-point, as
predicted by Theorem~\ref{5.6 Theorem}. \\
B) Next, we project the $3$-scroll $\tilde {X}:= S(2, 4, 4)
\subset {\mathbb P}^{12}_k$ of Example~\ref{6.1 Examples} C) from
the point
\[
p_4 = (0 : 1 : 0 : 0 : 0 : 0 : 0 : 0 : 0 : 0 : 0 : 0 : 0) \in
       {\mathbb P}^{12}_k \setminus \tilde {X}.
\]
We get a $3$-dimensional variety $X \subset {\mathbb P}^{11}_k$ of
degree 10 whose non-vanishing Betti numbers are:
\[
\begin{tabular}{c | r r r r r r r r r  }
    $i$ & 1 & 2 & 3 & 4 & 5 & 6 & 7 & 8 & 9  \\ \hline
    $u_i  $ & 34 & 151 & 314 & 364 & 230 & 69 & 7 & 0 & 0 \\
    $v_i  $ & 0 & 0 & 0 & 6 & 35 & 56 & 36 & 10 & 1  \\ \hline
\end{tabular}
\]
Now, $X$ is of arithmetic $\depth 3 = \dim (X).$ For each pair
$(s,t) \in k^2 \setminus \{ (0,0) \} $ consider the point
\[
\pi (s,t) := (s^2 : st : t^2 : 0 : 0 : 0 : 0 : 0 : 0 : 0 : 0 : 0 :
       0) \in \tilde {X} .
\]
Whenever $st \not= 0$, the two points $\pi (s,t)$ and $\pi (-
\sqrt{-1} s, \sqrt{-1} t)$ are different and the line joining them
contains $p_4$ and hence belongs to the secant cone
$\Sec_{p_4}(\tilde {X})$. Moreover the tangent line of the curve
\[
\tau : k \to \tilde {X} ; \ s \mapsto \pi (s,1)
\]
in the point $\tau (0) = \pi (0,1)$ runs through $p_4$ and thus
belongs to $\Sec_{p_4}(\tilde {X})$. Altogether this shows (cf
Theorem~\ref{5.6 Theorem} (d) ) that $\Sec_{p_4}(\tilde {X})$
coincides with the $2$-plane
\[
\{ (a : b : c : 0 : 0 : 0 : 0 : 0 : 0 : 0 : 0 : 0 : 0)
       \big\arrowvert (a : b : c) \in {\mathbb P}^2_k \} \subset
       {\mathbb P}^{12}_k .
\]
Projecting this plane from $p_4$ we obtain the line
\[
h:= \{ (a : c : 0 : 0 : 0 : 0 : 0 : 0 : 0 : 0 : 0 : 0)
       \big\arrowvert (a : c) \in {\mathbb P}^1_k \} \subset X.
\]
So, in the notation of Theorem~\ref{5.6 Theorem} we have $h =
\Sing(\nu )$. Let $a:= \frac{t}{s}, b:= \frac {s^2 v}{u^3}, c =
\frac{s^2 w}{u^3}$ and let $X_{s^2 u^3} \subset X$ be the affine
open set defined by $s^2 u^2 \not= 0$. It is easy to verify, that
there is an isomorphism $\varphi: X_{s^2 u^3} \overset {\simeq }
{\to } Y:= \Spec (k [a^2, b, ab, c, ac])$ such that $P:= (b, ab,
c, ac) \in Y$ is the generic point of $\varphi (h \cap X_{s^2
u^3})$. An easy calculation shows that ${\mathcal O}_{Y,P}$ is a
$G$-ring and hence, that the generic point of $h$ in $X$ is again
a $G$-point.
\end{example}

We now present a class of non-normal Del Pezzo varieties. Note
that these varieties are in fact arithmetically Gorenstein.

\begin{example}
\label{6.3 Example} A) Let $r \geq 4$ and let $\tilde {X} \subset
{\mathbb P}^{r + 1}_k$ be the rational surface scroll $S(2, r -
1)$, hence the variety which is defined by the $2 \times 2$-minors
of the matrix
\[
\begin{pmatrix} x_0 & x_1 \Big\arrowvert x_3 & x_4 & \cdots
                       & x_r  \\
                       x_1 & x_2 \Big\arrowvert x_4 & x_5 & \cdots
                       & x_{r + 1}
\end{pmatrix} .
\]
$\tilde{X}$  has the homogeneous coordinate ring
\[
B:= k [(s,t)^2 u^{r - 2} , (s,t)^{r - 2} v^2] \subset k
       [s, t, u, v] .
\]
Now, let $\varrho : {\mathbb P}^{r + 1}_k \setminus \{ p \}
\twoheadrightarrow {\mathbb P}^r_k, (x_0 : x_1 : x_2 : \cdots :
x_{r + 1}) \mapsto (x_0 : x_2 : x_3 : \cdots : x_{r + 1})$ be the
projection from the point $p = (0 : 1 : 0 : \cdots : 0) \in
{\mathbb P}^{r + 1}_k \setminus \tilde {X}$ and let $X:= \varrho
(\tilde {X}) \subset {\mathbb P}^r_k.$  Then $X$ is a surface and
has the homogeneous coordinate ring
\[
A:= k[s^2 u^{r - 2}, t^2 u^{r - 2}, (s, t)^{r - 2} v^2] \subset B.
\]
As $B$ is a birational extension of $A$, the morphism $\nu =
\varrho \upharpoonright  : \tilde {X} \to X$ is birational, so
that $\deg X = \deg \tilde{X} = r$ and $X \subset {\mathbb P}^r_k$
is a surface of almost minimal degree. Moreover, as $\tilde{X}$ is
smooth, $\nu = \varrho \upharpoonright : \tilde {X} \to X$ is a
normalization of $X$ and $\Sing(\nu ) = \varrho (\Sec_p(\tilde
{X}) \setminus \{ p \} )$ is the non-normal locus of $X$.

Similar as in example~\ref{6.2 Examples} B) we can check that the
secant cone of $\tilde {X}$at $p$ satisfies
\[
\Sec_p(\tilde {X}) = \{ (a : b : c : 0 : \cdots : 0)
       \big\arrowvert (a : b : c) \in {\mathbb P}^2_k \}
\]
and hence is a $2$-plane. So, by Theorem~\ref{5.6 Theorem}, $X$
cannot be of arithmetic $\depth \leq 2 = \dim X$ and hence is
arithmetically Cohen-Macaulay.

Moreover
\[
h := X \setminus \Nor(X) = \varrho (\Sec_p(\tilde {X}) \setminus
\{ p \}) =
       \{ (a : c : 0 : \cdots : 0) \in {\mathbb P}^r_k
       \big\arrowvert (a : c) \in {\mathbb P}^1_k \} .
\]
So, the non-normal locus $h$ of $X$ is a line.

Now consider the affine open set $X_{s^2 u^{r - 2}} \subset X$
defined by $s^2 u^{r - 2} \not= 0$ and let $a:= \frac{s^{r - 5}
v^2 t} {u^{r - 2}}$ and $b:= \frac {s^{r - 4}v^2}{u^{r - 2}}$.
Then, an easy calculation shows that there is an isomorphism
\[
\varphi : X_{s^2 u^{r - 2}} \overset {\simeq }{\longrightarrow }
       Y:= \Spec(k[a, b, \frac{a^2}{b^2}]) = \Spec
       (k[a, b, c]/(cb^2 - a^2))
\]
which maps $h \cap X_{s^2 u^{r - 2}}$ to the singular line $a = b
= 0$ of the surface $Y$. The pinch point $\underline {0}$ of $Y$
can be written as $\varphi (\varrho (\ell \setminus \{ p \} ))$,
where $\ell $ is the tangent line to $\tilde {X}$ at the point $(1
: 0 : \cdots : 0)$ which contains $p$.

The same arguments apply to the affine open set $X_{t^2 u^{r - 2}}
\subset X$. This allows to conclude that the open neighborhood
$X_{s^2 u^{r - 2}} \cup X_{t^2 u^{r - 2}}$ of the singular line
$h$ of $X$ is isomorphic to the blow-up $\Proj (k[a, b][a^2T,
b^2T])$ of the affine plane ${\mathbb A}^2_k = \Spec(k[a, b])$
with respect to the polynomials $a^2$ and $b^2$. \\
B) Let $r \geq 5$ and let $\tilde {X} \subset {\mathbb P}^{r +
1}_k$ be the rational normal 3-scroll $S(1,1, r - 3)$, hence the
variety which is defined by the $2 \times 2$-minors of the matrix
\[
\begin{pmatrix} x_0 \Big\arrowvert x_2 \Big\arrowvert x_4 & x_5 &
                       \cdots & x_r  \\
                       x_1 \Big\arrowvert  x_3 \Big\arrowvert x_5 &
                       x_6 & \cdots
                       & x_{r + 1}
\end{pmatrix}.
\]
$\tilde{X}$ has the homogeneous coordinate ring
\[
B:= k [(s,t)u^{r - 4} , (s,t)v^{r -4}, (s,t)^{r - 4} w] \subset k
       [s, t, u, v, w] .
\]
Now, let $\varrho : {\mathbb P}^{r + 1}_k \setminus \{ p \}
\twoheadrightarrow {\mathbb P}^r_k, (x_0 : x_1 : x_2 : \cdots :
x_{r + 1}) \mapsto (x_0 : x_1 - x_2 : x_3 : \cdots : x_{r + 1})$
be the projection from the point $p = (0 : 1 : 1 : 0 : \cdots : 0)
\in {\mathbb P}^{r + 1}_k \setminus \tilde {X}$ and let $X:=
\varrho (\tilde {X}) \subset {\mathbb P}^r_k.$ Then $X$ is of
dimension 3 and has the homogeneous coordinate ring
\[
A:= k[s u^{r - 4}, tu^{r-4} - sv^{r-4}, t v^{r - 4}, (s, t)^{r -
4} w] \subset B.
\]
As $B$ is a birational extension of $A$, the morphism $\nu =
\varrho \upharpoonright  : \tilde {X} \to X$ is birational, so
that $\deg X = \deg \tilde{X} = r$ and $X \subset {\mathbb P}^r_k$
has dimension $3$ and is of almost minimal degree. Moreover, as
$\tilde{X}$ is smooth, $\nu = \varrho \upharpoonright : \tilde {X}
\to X$ is a normalization of $X$ and $\Sing(\nu ) = \varrho
(\Sec_p(\tilde {X}) \setminus \{ p \} )$ is the non-normal locus
of $X$.

Similar as in example~  A) above we can check that the
secant cone of $\tilde {X}$at $p$ satisfies
\[
\Sec_p(\tilde {X}) = \{ (a : b : c : d : 0 : \cdots : 0)
       \big\arrowvert (a : b : c: d) \in {\mathbb P}^3_k \}
\]
and hence is a $3$-plane. So, by Theorem~\ref{5.6 Theorem} the
variety $X$ cannot be of arithmetic $\depth \leq 3 = \dim X$ and
hence is arithmetically Cohen-Macaulay, that is a non-normal Del
Pezzo variety of dimension 3.

Moreover
\[
\varrho (\Sec_p(\tilde {X}) \setminus \{ p \}) =
       \{ (a : b : d : 0 : \cdots : 0) \in {\mathbb P}^r_k
       \big\arrowvert (a : b : d) \in {\mathbb P}^2_k \} .
\]
So, the non-normal locus of $X$ is a plane, in accordance with
Proposition~\ref{5.8 Proposition} and Corollary~\ref{6.10
Corollary}.
\end{example}

Finally observe that $X$ in \ref{6.3 Example} A) is a divisor on
the variety of minimal degree $\varrho (Z) \subset {\mathbb
P}^r_k$, where $Z \subset {\mathbb P}^{r + 1}_k$ is the variety
defined by the $2 \times 2$ minors of the matrix
\[
\begin{pmatrix} x_3 & x_4 & \cdots & x_r  \\
                       x_4 & x_5 & \cdots & x_{r + 1}
       \end{pmatrix} .
\]

In the previous example we have met arithmetically
Cohen-Macaulay varieties of almost minimal degree which occur as a
subvariety of codimension one on a variety of minimal degree. We
now present an example of a normal Del Pezzo variety which does
not have this property.

\begin{example}
\label{6.4 Example} Let $X \subset {\mathbb P}^9_k$ be the smooth
$6$-dimensional arithmetically Gorenstein variety of degree $5$
defined by the $4 \times 4$ Pfaffian quadrics $F_1, F_2, F_3, F_4,
F_5$ of the skew symmetric matrix (cf \cite{BuE})
\[
M =   \begin{pmatrix} 0 & x_0 & x_1 & x_2 & x_3&   \\
                       - x_0 & 0 & x_4 & x_5 & x_6 \\
                       - x_1 & - x_4 & 0 & x_7 & x_8 \\
                       - x_2 & - x_5 & - x_7 & 0 & x_9 \\
                       - x_3 & - x_6 & - x_8 & - x_9 & 0
       \end{pmatrix} .
\]

According to \cite{BuE} the columns of $M$ provide a minimal
system of generators for the first syzygy module of the
homogeneous vanishing ideal $I \subset S = k[x_0, x_2, \cdots ,
x_9]$ of $X$. Assume now that there is a variety $W
\subset{\mathbb P}^9_k$ of minimal degree with $\dim W = 7$ and $X
\subset W$. So, $W$ is arithmetically Cohen-Macaulay and of
codimension $2$ and by the Theorem of Hilbert-Burch the
homogeneous vanishing ideal $J \subset S$ of $W$ is generated by
the three $2 \times 2$-minors $G_1, G_2, G_3 \in S_2$ of a $2
\times 3$-matrix with linearly independent entries in $S_1$ (cf
\cite{E}). So, after an eventual renumbering of the generators
$F_i$, we may assume that $G_1, G_2, G_3, F_4, F_5 \in I_2$ is a
minimal system of generators of $I$. As $J$ admits two independent
syzygies
\[
\lambda _{i1} G_1 + \lambda _{i2} G_2 + \lambda _{i3}
       G_3 = 0,  \; \lambda _{ij} \in S_1,  i = 1, 2,  j = 1, 2, 3,
\]
a minimal system of generators for the first syzygy module of $I$
would be given by the matrix of the form
\[
N = \begin{pmatrix} 0 & 0 & \ast & \ast & \ast \\
                       0 & 0 & \ast & \ast & \ast \\
                       \ast & \ast & \ast & \ast & \ast \\
                       \ast & \ast & \ast & \ast & \ast \\
                       \ast & \ast & \ast & \ast & \ast
       \end{pmatrix}
       \in S_1^{ 5 \times 5}.
\]
On the other hand there should be a $k$-linear transformation
which converts $N$ into $M$ -- a contradiction.
\end{example}

\begin{remark}
\label{7.5 Remark} A) The variety $X \subset \mathbb P^9_k$ of
Example \ref{6.4 Example} is normal and Dell Pezzo and hence not a
projection of a variety $\tilde{X} \subset \mathbb P^{10}_k$ of
minimal degree. The non-existence of the above variety $W \subset
\mathbb P^9_k$ of minimal degree thus is in accordance with
Theorem \ref{Theorem 8.3}. By Remark \ref{Remark 8.5} the
projection $X \subset \mathbb P^4_k$ of the Veronese surface $F
\subset \mathbb P^5_k$ is not contained in a variety $Y \subset
\mathbb P^4_k$ of minimal degree either, according to the fact,
that $F$ is not a scroll. So Remark \ref{Remark 8.5} and Example
\ref{6.4 Example} illustrate that the hypotheses of Theorem
\ref{Theorem 8.3} cannot be weakened.\\
B) The examples of this section (with the execption of the last
one) are all of relatively big codimension. It turns out, that the
structure of varieties of almost minimal degree and small
codimension is fairly fixed and cannot vary very much. We study
these varieties more extensively in \cite{BS3}.
\end{remark}

\subsection*{Acknowledgment} 
We thank the referee for his valuable hints concerning further 
investigations on birational projections from varieties of minimal and 
almost minimal degree. In fact, a part of his ideas come up already in 
an on-going joint investigation of the authors on projective surfaces $X$ 
of degree $r+1$ in $P^r_k.$ We thank also Euisung Park for his comment 
concerning Proposition \ref{3.4 Proposition}.

\bibliographystyle{plain}

\vskip 2 cm

\author{
 \begin{tabular}{llll}
        M. Brodmann                      &&P. Schenzel \\
        Institut f\"ur Mathematik        &&Martin-Luther-Univ.
                 Halle-Wittenberg \\
        Universit\"at Z\"urich           &&Fachbereich Mathematik und
                 Informatik \\
        Winterthurerstrasse 190          && Von-Seckendorff-Platz 1\\
        CH-8057 Z\"urich, Schwitzerland  && D-06120 Halle (Saale), Germany \\
        & \\
       {\it email: } brodmann@math.unizh.ch
        &&{\it email: } schenzel@informatik.uni-halle.de

 \end{tabular}
}

\end{document}